\documentclass[11pt,twoside,draft]{article}
\usepackage{citesort}
\usepackage{amssymb}
\usepackage{amsmath}
\usepackage[mathscr]{eucal}
\usepackage{amsfonts}
\usepackage{latexsym}
\usepackage{amsxtra}
\usepackage{amsbsy}
\usepackage{amsthm}
\usepackage{amscd}
\usepackage{amsopn}
\usepackage{amstext}
\newtheorem{ddddd}{Theorem}[section]

\newtheorem{hhhhh}{Conjecture}[section]

\newtheorem{proposition}{Proposition}[section]
\newtheorem{ccc}{Lemma}[section]
\newtheorem{ddd}{Theorem}[section]

\newtheorem{hhh}{Conjecture}[section]

\newtheorem{proposition:s}{Proposition}[subsection]
\newtheorem{cccc}{Lemma}[subsection]
\newtheorem{dddd}{Theorem}[subsection]
\newtheorem{ffff}{Corollary}[subsection]

\theoremstyle{definition}
\newtheorem{eee}{Remark}[section]

\theoremstyle{definition}
\newtheorem{eeee}{Remark}[subsection]

\theoremstyle{definition}

\newcommand{\me}{\mathrm{e}}
\newcommand{\mi}{\mathrm{i}}
\newcommand{\md}{\mathrm{d}}
\renewcommand{\Im}{\mathrm{Im}}
\renewcommand{\Re}{\mathrm{Re}}
\allowdisplaybreaks[4]
\begin{document}
\baselineskip=12pt
\frenchspacing
\title{Connection Formulae for Asymptotics of Solutions of the Degenerate 
Third Painlev\'{e} Equation. I}
\author{A.~V.~Kitaev\thanks{\texttt{E-mail: kitaev@pdmi.ras.ru}. Present 
address: School of Mathematics and Statistics F07, University of Sydney, NSW 
2006, Australia, \texttt{e-mail: kitaev@maths.usyd.edu.au}} \\
Steklov Mathematical Institute \\
Fontanka 27 \\
St. Petersburg 191023 \\
Russia \and 
A.~H.~Vartanian\thanks{\texttt{E-mail: arthur@math.duke.edu}. On leave of 
absence from: Department of Mathematics, Winthrop University, Rock Hill, South 
Carolina 29733, U.~S.~A., \texttt{e-mail: vartaniana@winthrop.edu}} \\
Department of Mathematics \\
Duke University \\
Durham, North Carolina 27708 \\
U.~S.~A.}
\date{27 November 2003}
\maketitle
\begin{abstract}
\noindent
The degenerate third Painlev\'{e} equation, $u^{\prime \prime} \! = \! \frac{
(u^{\prime})^{2}}{u} \! - \! \frac{u^{\prime}}{\tau} \! + \! \frac{1}{\tau}(-
8 \varepsilon u^{2} \! + \! 2ab) \! + \! \frac{b^{2}}{u}$, where $\varepsilon,
b \! \in \! \mathbb{R}$, and $a \! \in \! \mathbb{C}$, and the associated 
tau-function are studied via the Isomonodromy Deformation Method. Connection 
formulae for asymptotics of the general as $\tau \! \to \! \pm 0$ and $\pm 
\mi 0$ solution and general regular as $\tau \! \to \! \pm \infty$ and $\pm 
\mi \infty$ solution are obtained.

\vspace{1.35cm}
\! \! \! {\bf 2000 Mathematics Subject Classification.} 33E17, 34M40, 34M50, 
34M55, 34M60

\vspace{0.50cm}
{\bf Abbreviated Title.} Degenerate Third Painlev\'{e} Equation

\vspace{0.50cm}

tau-function, Schlesinger transformations, WKB method, Stokes phenomena
\end{abstract}
\clearpage
\section{Introduction}
In this paper we study the degenerate third Painlev\'{e} equation,
\begin{equation}
\label{eq:dp3}
u^{\prime \prime} \! = \! \dfrac{(u^{\prime})^{2}}{u} \! - \! \dfrac{u^{
\prime}}{\tau} \! + \! \dfrac{1}{\tau}(-8 \varepsilon u^{2} \! + \! 2 ab) \! 
+ \! \dfrac{b^{2}}{u},
\end{equation}
where $u \! = \! u(\tau)$, the primes denote differentiation with respect to 
$\tau$ (or $t$: see below), $a$ and $b$ $(\not= \! 0)$ are $\mathbb{C}$-valued 
parameters, and $\varepsilon \! = \! \pm 1$. Originally, 
Equation~(\ref{eq:dp3}) appeared as a special case of, and in the same context 
as, the complete third Painlev\'{e} equation,
\begin{equation}
\label{eq:p3}
y^{\prime \prime} \! = \! \dfrac{(y^{\prime})^{2}}{y} \! - \! \dfrac{y^{
\prime}}{t} \! + \! \dfrac{1}{t} \! \left(\alpha y^{2} \! + \! \beta \right) 
\! + \! \gamma y^{3} \! + \! \dfrac{\delta}{y},
\end{equation}
in the studies of Painlev\'{e}; however, its particular significance was not 
distinguished at that time.

There is another relation between Equations~(\ref{eq:dp3}) and~(\ref{eq:p3}) 
(which is why Equation~(\ref{eq:dp3}) is called degenerate), namely, the 
following \emph{double-scaling limit},
\begin{gather}
\begin{split}
\label{eq:limit}
t &= \epsilon \tau, \qquad \quad y(t) \! \underset{\underset{\tau =\mathcal{O}
(1)}{\epsilon \to 0}}{=} \! \epsilon u(\tau) \! + \! o(\epsilon), \\
\alpha \! = \! -\dfrac{8 \varepsilon}{\epsilon^{2}}&, \qquad \beta \! = \! 
2ab, \qquad \gamma \! = \! o(\epsilon^{-4}), \qquad \delta \! = \! b^{2}.
\end{split}
\end{gather}
The $o(\epsilon)$ term in the formula for $y(t)$ above also depends on $\tau$: 
it is supposed that the derivatives with respect to $\tau$ of this term are 
also of the order $o(\epsilon)$. Only a few of the formal properties of 
Equation~(\ref{eq:dp3}), namely, the Hamiltonian structure and B\"{a}cklund 
transformations, can be derived {}from those for Equation~(\ref{eq:p3}) as a 
straightforward consequence of the double-scaling limit (\ref{eq:limit}), 
whereas its analytic and asymptotic properties require separate 
considerations. Recently, Equation~(\ref{eq:dp3}) has appeared in a number 
of physical \cite{a1,a2,a3,a4} and geometrical applications \cite{a5} (in 
contexts independent of Equation~(\ref{eq:p3})) where knowledge of asymptotic 
properties of its solutions is of special importance.

The B\"{a}cklund transformation for Equation~(\ref{eq:dp3}) was obtained by 
Gromak \cite{G1}: he also proved \cite{G2} that the only algebraic solutions 
of the complete third Painlev\'{e} equation are rational functions of 
$\tau^{1/3}$. Actually, these functions are solutions of 
Equation~(\ref{eq:dp3}) for $a \mi \! = \! n \! \in \! \mathbb{Z}$ and 
arbitrary, non-vanishing values of $\varepsilon$ and $b$. For fixed values of 
$a \mi \! = \! n$, $b$ and $\varepsilon$, there is exactly one algebraic 
solution of Equation~(\ref{eq:dp3}) which is a multi-valued function with 
three branches. This solution can be obtained by applying $\vert n \vert$ 
B\"{a}cklund transformations to the simplest solution of 
Equation~(\ref{eq:dp3}), namely, $u\! = \! b^{2/3} \tau^{1/3}/2 \varepsilon$ 
(for $a \! = \! 0)$. All other solutions are non-classical in the sense of 
Darboux-Umemura \cite{GD,U}, that is, the absence of invariant algebraic 
curves of the corresponding Hamiltonian vector field: this was recently 
proved by Ohyama \cite{YO}.

The most efficacious approach for studying the asymptotic behaviour of general 
solutions of the Painlev\'{e} equations and, especially, connection formulae 
for their asymptotics, is the method of isomonodromic deformations 
\cite{a10,a11,a12,a13,a9}, or a closely related technique based on a 
steepest descent-type analysis of the associated Riemann-Hilbert problem 
\cite{a14}. Equation~(\ref{eq:dp3}), in the case $a \! = \! 0$, was studied by 
the isomonodromy deformation method in \cite{a15}, where asymptotics as $\tau 
\! \to \! 0$ and $\infty$, as well as the corresponding connection formulae, 
were obtained: that work was based on the study of isomonodromic deformations 
of a $3 \! \times \! 3$ matrix linear ODE with two irregular singular points. 
At the same time, as noted in \cite{a16}, there is another matrix linear ODE 
in terms of $2 \! \times \! 2$ matrices whose isomonodromy deformations are 
described by Equation~(\ref{eq:dp3}) (with arbitrary $a)$. In this work, the 
latter ODE is used for the analysis of Equation~(\ref{eq:dp3}).
\begin{proposition}[\textrm{\cite{a16}}]
The necessary and sufficient condition for the compatibility of the 
linear system
\begin{equation}
\label{eq:UV}
\partial_{\lambda} \Phi (\lambda,\tau) \! = \! \mathscr{U}(\lambda,\tau) 
\Phi (\lambda,\tau), \qquad \quad \partial_{\tau} \Phi (\lambda,\tau) \! = 
\! \mathscr{V}(\lambda,\tau) \Phi (\lambda,\tau),
\end{equation}
with
\begin{gather*}
\mathscr{U}(\lambda,\tau) \! = \! \tau \! \left(\! -\mi \sigma_{3} \! - \! 
\dfrac{a \mi}{2 \tau \lambda} \sigma_{3} \! - \! \dfrac{1}{\lambda} \! 
\begin{pmatrix}
0 & C(\tau) \\
D(\tau) & 0
\end{pmatrix} \! + \! \dfrac{\mi}{2 \lambda^{2}} \! 
\begin{pmatrix}
\sqrt{\smash[b]{-A(\tau)B(\tau)}} & A(\tau) \\
B(\tau) & -\sqrt{\smash[b]{-A(\tau)B(\tau)}}
\end{pmatrix} \right), \\
\mathscr{V}(\lambda,\tau) \! = \! \left(\! -\mi \lambda \sigma_{3} \! + \! 
\dfrac{a \mi}{2 \tau} \sigma_{3} \! - \! 
\begin{pmatrix}
0 & C(\tau) \\
D(\tau) & 0
\end{pmatrix} \! - \! \dfrac{\mi}{2 \lambda} \! 
\begin{pmatrix}
\sqrt{\smash[b]{-A(\tau)B(\tau)}} & A(\tau) \\
B(\tau) & -\sqrt{\smash[b]{-A(\tau)B(\tau)}}
\end{pmatrix} \right),
\end{gather*}
$\sigma_{3} \! = \! 
\left(
\begin{smallmatrix}
1 & 0 \\
0 & -1
\end{smallmatrix}
\right)$, $\partial_{z} \! := \! \tfrac{\partial}{\partial z}$, and 
differentiable, scalar-valued functions $A(\tau)$, $B(\tau)$, $C(\tau)$, and 
$D(\tau)$, is that $A(\tau)$, $B(\tau)$, $C(\tau)$, and $D(\tau)$ satisfy the 
following system of isomonodromy deformations:
\begin{equation}
 \label{sys:ABCD}
\begin{gathered}
A^{\prime}(\tau) \! = \! 4C(\tau) \sqrt{\smash[b]{-A(\tau)B(\tau)}}, \qquad 
\quad B^{\prime}(\tau) \! = \! -4D(\tau) \sqrt{\smash[b]{-A(\tau)B(\tau)}}, \\
(\tau C(\tau))^{\prime} \! = \! 2a \mi C(\tau) \! - \! 2 \tau A(\tau), \qquad 
\quad (\tau D(\tau))^{\prime} \! = \! -2a \mi D(\tau) \! + \! 2 \tau B(\tau), 
\\
\left(\sqrt{\smash[b]{-A(\tau)B(\tau)}} \right)^{\prime} \! = \! 2(A(\tau)
D(\tau) \! - \! B(\tau)C(\tau)).
\end{gathered}
\end{equation}
\end{proposition}
\emph{Proof.} Follows {}from the Frobenius compatibility condition, $\partial_{
\tau} \partial_{\lambda} \Phi(\lambda,\! \tau) \! = \!\partial_{\lambda} 
\partial_{\tau} \Phi(\lambda,\! \tau)$. \hfill $\square$
\begin{eee}
System~(\ref{sys:ABCD}) implies that when 
$A(\tau) \! \equiv \! 0 \! \Rightarrow \! B(\tau)
C(\tau) \! = \! 0$, and, when $B(\tau) \! \equiv \! 0 \! \Rightarrow \! 
A(\tau) D(\tau) \! = \! 0$: for these cases, System~(\ref{sys:ABCD}) can be 
solved exactly; hence, one excludes these cases {}from further 
consideration. Hereafter, all explicit $\tau$ dependencies are suppressed, 
except where confusion may arise. \hfill $\blacksquare$
\end{eee}
\begin{proposition}
Let $u \! = \! u(\tau)$ and $\varphi \! = \! \varphi (\tau)$ solve the 
following system,
\begin{gather}
\begin{split}
u^{\prime \prime} \! = \! \frac{(u^{\prime})^{2}}{u} \! - \! \frac{u^{\prime}}
{\tau} \! &+ \! \frac{1}{\tau}(-8 \varepsilon u^{2} \! + \! 2ab) \! + \! 
\frac{b^{2}}{u}, \\
\varphi^{\prime} &= \! \frac{2a}\tau \! + \! \frac{b}{u},
\end{split}
\end{gather}
where $\varepsilon \! = \! \pm 1$, and $a,b \! \in \! \mathbb{C}$ are 
independent of $\tau$. Then
\begin{equation*}
A(\tau) \! := \! \dfrac{u(\tau)}{\tau} \me^{\mi \varphi (\tau)}, \quad  
B(\tau) \! := \! -\dfrac{u(\tau)}{\tau} \me^{-\mi \varphi (\tau)}, \quad 
C(\tau) \! := \! \dfrac{\varepsilon \tau}{4u(\tau)} A^{\prime}(\tau), \quad 
D(\tau) \! := \! - \frac{\varepsilon \tau}{4u(\tau)}B^{\prime}(\tau),
\end{equation*}
solve System~{\rm (\ref{sys:ABCD})}. Conversely, let 
$A(\tau) \! \not\equiv \! 0$, $B(\tau) 
\! \not\equiv \! 0$, $C(\tau)$, and $D(\tau)$ solve 
System~{\rm(\ref{sys:ABCD})}, and define
\begin{equation*}
u(\tau) \! := \! \varepsilon \tau \sqrt{\smash[b]{-A(\tau)B(\tau)}}, \quad \, 
\, \varphi (\tau) \! := \! -\tfrac{\mi}{2} \ln (-A(\tau)/B(\tau)), \quad 
\mathrm{and} \quad b \! := \! u(\tau)(\varphi^{\prime}(\tau) \! - \! 2a/\tau).
\end{equation*}
Then $b$ is independent of $\tau$, and $u(\tau)$ solves 
Equation~{\rm(\ref{eq:dp3})}.
\end{proposition}

\emph{Proof.} Straightforward verification. \hfill $\square$
\begin{eee}
For $b \! = \! 0$, System~(6) yields $u(\tau) \! = \! \tfrac{\varepsilon 
\widetilde{c}_{1}}{16 \tau} \! - \! \tfrac{\varepsilon \widetilde{c}_{1} 
\left( 1 + \widetilde{c}_{2} \tau^{\sqrt{\widetilde{c}_{1}}} \right)^{2}}{16 
\tau \left( 1 - \widetilde{c}_{2} \tau^{\sqrt{\widetilde{c}_{1}}} \right)^{
2}}$, where $\widetilde{c}_{i} \! \in \! \mathbb{C}$, $i \! = \! 1,2$, and 
$\varphi (\tau) \! = \! 2a \ln \tau \! + \! \varphi_{o}$, with $\varphi_{o} 
\! \in \! \mathbb{C}$; therefore, hereafter, only the case $b \! \not=\! 
0$ will be considered. In the latter case, one can rescale the parameters 
$\varepsilon, b$ {}from Equation~(\ref{eq:dp3}) via the following sequence of 
scaling transformations: 
(1) $\tau \! \to \! \tau /\sqrt{\smash[b]{\varepsilon b}}$ 
and $u \! \to \! u \sqrt{\smash[b]{\varepsilon b}}$; followed by 
(2) $\tau \! \to \! \varepsilon \tau$ and $u \! \to \! u$. As these algebraic 
rescalings change $\arg \tau$, which is important for asymptotics, it is 
convenient to keep the parameters $\varepsilon, b$ 
explicitly. \hfill $\blacksquare$
\end{eee}

The Hamiltonian structure of Equation~(\ref{eq:dp3}) can be derived {}from 
the corresponding structure of Equation~(\ref{eq:p3}) \cite{a11,a17} by 
applying the double-scaling limit~(\ref{eq:limit}).
\begin{proposition}
Let
\begin{equation}
\mathcal{H}_{\varepsilon_{1}}(p,q;\tau) \! := \! \frac{p^{2}q^{2}}{\tau} 
\! - \! \frac{2 \varepsilon_{1}pq(a \mi \! +\! 1/2)}{\tau} \! + \! 4 
\varepsilon q \! + \! \mi bp \! + \! \frac{(a \mi \! +\! 1/2)^{2}}{2 
\tau},
\end{equation}
where $\varepsilon^{2} \! = \! \varepsilon_{1}^{2} \! = \! 1$. Then 
Hamilton's equations,
\begin{equation}
\frac{\md p}{\md \tau} \! = \! -\frac{\partial \mathcal{H}_{\varepsilon_{1}}
(p,q;\tau)}{\partial q} \qquad \mathrm{and} \qquad \frac{\md q}{\md \tau} \! = 
\! \frac{\partial \mathcal{H}_{\varepsilon_{1}}(p,q;\tau)}{\partial p},
\end{equation}
are equivalent to either one of the following degenerate third Painlev\'{e} 
equations,
\begin{align}
p^{\prime \prime} \! &= \! \dfrac{(p^{\prime})^{2}}{p} \! - \! \dfrac{p^{
\prime}}{\tau} \! + \! \dfrac{1}{\tau} \! \left(-2 \mi bp^{2} \! + \! 8 
\varepsilon \! \left(a \mi \varepsilon_{1} \! + \! \tfrac{(\varepsilon_{1}-
1)}{2} \right) \right) \! - \! \dfrac{16}{p}, \\
q^{\prime \prime} \! &= \! \dfrac{(q^{\prime})^{2}}{q} \! - \! \dfrac{q^{
\prime}}{\tau} \! + \! \dfrac{1}{\tau} \! \left(-8 \varepsilon q^{2} \! - \! 
b(2a \varepsilon_{1} \! - \! \mi (1 \! + \! \varepsilon_{1})) \right) \! + \! 
\dfrac{b^{2}}{q}.
\end{align}
\end{proposition}

\emph{Proof.} The Hamiltonian System~(8) can be rewritten as
\begin{equation}
p \! = \! \dfrac{\tau (q^{\prime} \! - \! \mi b)}{2q^{2}} \! + \! \dfrac{(a 
\mi \! + \! 1/2) \varepsilon_{1}}{q}, \qquad \quad q \! = \! -\dfrac{\tau 
(p^{\prime} \! + \! 4 \varepsilon)}{2p^{2}} \! + \! \dfrac{(a \mi \! + \! 1/2) 
\varepsilon_{1}}{p}.
\end{equation}
The proposition follows {}from Equation~(11) by straightforward calculations. 
\hfill $\square$
\begin{eee}
Either expression in Equation~(11) is equivalent to the B\"{a}cklund 
transformation \cite{G2} for Equation~(\ref{eq:dp3}). The connection of 
$\mathcal{H}_{\varepsilon_{1}}(p,q;\tau)$ with the isomonodromy deformations 
is given in Remark~2.2 of Section~2. \hfill $\blacksquare$
\end{eee}

It is worth noting two other differential equations related to the degenerate 
third Painlev\'{e} equation:
\begin{gather*}
\tau^{2} \! \left(f^{\prime \prime} \! + \! 4 \mi \varepsilon b \right)^{2} 
\! = \! \left(4f \! - \! \varepsilon_{1}(2 \mi a \! + \! 1) \right)^{2} \! 
\left((f^{\prime})^{2} \! + \! 8 \mi \varepsilon b f \right), \\
\left(\tau \sigma^{\prime \prime} \! - \! \sigma^{\prime} \right)^{2} \! = 
\! 2 \! \left(2 \sigma \! - \! \tau \sigma^{\prime} \right) \! \left(\sigma^{
\prime} \right)^{2} \! - \! 32 \mi \varepsilon b \tau \! \left(\! \left(
\tfrac{1-\varepsilon_{1}}{2} \! - \! a \mi \varepsilon_{1} \right) \! \sigma^{
\prime} \! + \! 2 \mi \varepsilon b \tau \right),
\end{gather*}
where
\begin{equation*}
f(\tau) := \dfrac{p(\tau)q(\tau)}{2},
\end{equation*}
and
\begin{eqnarray*}
\sigma (\tau) \! \! \! &:=& \! \! \! \tau\mathcal{H}_{\varepsilon_{1}}(p(\tau),
q(\tau);\tau) \! + \! 2f(\tau) \! + \! \tfrac{1}{2} \! \left(a \mi \! + \! 
\tfrac{1}{2} \right)^{2} \! - \! \varepsilon_{1} \! \left(a \mi \! + \! 
\tfrac{1}{2}\right) \! + \! \tfrac{1}{4} \\
&=& \! \! \! \left(p(\tau)q(\tau) \! - \! \varepsilon_1 \! \left(a \mi \! + \! 
\tfrac{1}{2} \! - \! \tfrac{\varepsilon_1}{2} \right) \right)^2 \! + \! \tau 
(4 \varepsilon q(\tau) \! + \! \mi bp(\tau)).
\end{eqnarray*}
The functions $u(\tau)$ and  $f(\tau)$ solve some well-known 
integrable differential-difference and difference equations which are given 
in Section~\ref{sec:6}. In this work, asymptotics of Equation~(\ref{eq:dp3}) 
as $\tau \! \to \! \pm 0$, $\pm \mi 0$, $\pm \infty$, and $\pm \mi \infty$ 
are parametrised in terms of the monodromy data of System~(\ref{eq:UV}). 
This parametrisation is equivalent to finding the corresponding connection 
formulae: indeed, given asymptotics of some solution as $\tau \! \to \! \pm 
0$ (resp., $\pm \mi 0)$ or $\pm \infty$ (resp., $\pm \mi \infty)$, it is 
straightforward to determine the corresponding monodromy data, and therefore 
to obtain asymptotics of the same solution as $\tau \! \to \! \pm \infty$ 
(resp., $\pm \mi \infty)$ or $\pm 0$ (resp., $\pm \mi 0)$.

There are different technical means by which isomonodromy deformations can be 
used to obtain the desired parametrisation of asymptotics. For the first time, 
isomonodromy deformations for the solution of the connection problem for the 
six Painlev\'{e} equations were used by Jimbo \cite{a13}; however, that work 
only dealt with asymptotics of solutions of the Painlev\'{e} equations in the 
neighbourhood of the regular singular points. At about the same time, Its 
\emph{et al.} \cite{IP} used isomonodromy deformations for the parametrisation 
of asymptotics of a special solution to the third Painlev\'{e} equation in the 
neighbourhood of its essential singular point by means of the corresponding 
monodromy data; however, this was done in the context of the solution of an 
asymptotic problem for special (rapidly decaying) solutions of the sine-Gordon 
equation rather than the connection problem for the third Painlev\'{e} 
equation. Novokshenov \cite{N} used similar methods to solve the connection 
problem for a particular solution of the third Painlev\'{e} equation. These 
asymptotic methods for the Painlev\'{e} equations (see, also, \cite{a18}), 
based on the asymptotic analysis of associated linear ODEs, were developed 
further in a number of works; for a recent example, see \cite{a19}. Here, we 
also follow these ideas.

Primary emphasis is given to the global asymptotic properties of $\Phi 
(\lambda,\tau)$, that is, the possibility of matching different local 
asymptotic expansions of $\Phi (\lambda,\tau)$ at singular and turning points, 
which culminates in the asymptotic expansion of solutions of the Painlev\'{e} 
equations. The justification scheme for this asymptotic method was suggested 
in \cite{a20}. Further details concerning its application to particular 
Painlev\'{e} equations can be found in \cite{a21,a22}.

This paper is organized as follows. In Section~\ref{sec:2}, the principal 
object of study, that is, the manifold of monodromy data, is introduced. In 
Section~3, the results of this work, namely, asymptotics of general solutions 
of Equation~(\ref{eq:dp3}) parametrised in terms of the monodromy data (points 
of the monodromy manifold), are given. In Section~4 (resp., Section~5), the 
$\tau \! \to \! +\infty$ (resp., $\tau \! \to \! +0)$ asymptotic analysis 
for the general solutions of Equation~(\ref{eq:dp3}) is presented. In 
Section~\ref{sec:6}, the B\"{a}cklund transformations for 
Equation~(\ref{eq:dp3}) are derived as a consequence of the Schlesinger 
transformations (cf.~\cite{a16}). This allows one to remove the restriction 
on the parameter (of formal monodromy) $a$ which is imposed in Section~5, 
and to extend the corresponding asymptotic (and connection) results as $\tau 
\! \to \! +0$ to any complex value of the parameter $a$. In the same section 
we also find the actions of the Lie-point symmetries for 
System~(\ref{sys:ABCD}) on the manifold of the monodromy data: the latter is 
used to extend the connection results found in the previous sections for
asymptotics on the positive semi-axis to asymptotics on the negative 
semi-axis, and, in the Appendix, for asymptotics on the imaginary axis.

Many modern applications actually require knowledge of the function $\mathcal{
H}(\tau)$ rather than the original Painlev\'{e} function $u(\tau)$. Therefore, 
we present the corresponding asymptotic results for the function $\mathcal{H}
(\tau)$. Actually, our asymptotics for the function $\mathcal{H}(\tau)$ as 
$\tau \! \to \! \pm \infty$ is ``incomplete'' as the leading term that 
directly follows from our calculations is parametrised by only one monodromy 
parameter, while there should be two. To obtain a complete parametrisation, 
one has to find one more term in the asymptotic expansion of $\mathcal{H}
(\tau)$: the latter requires the calculation of two further terms in the 
asymptotic expansion for the function $u(\tau)$. This calculation is not 
related to the isomonodromy deformation technique considered here, but can 
be obtained by a direct substitution of the corresponding ansatz into 
Equation~(\ref{eq:dp3}). Our asymptotic results as $\tau \! \to \! \pm \infty$ 
are also incomplete in another sense: we haven't found asymptotics for all 
domains of the monodromy manifold; in particular, asymptotics as $\tau \! \to 
\! \pm \infty$ of the singular real solutions of Equation~(\ref{eq:dp3}). We 
are planning on closing both ``gaps'' mentioned above in a subsequent work.
\section{The Manifold of Monodromy Data}
\label{sec:2}
The leading terms, as $\lambda \! \to \! 0$, of the matrices $\mathscr{U}
(\lambda,\tau)$ and $\mathscr{V}(\lambda,\tau)$ in System~(\ref{eq:UV}) are 
degenerate. Therefore, we begin with a Fabry-type transformation (see 
\cite{V}) to obtain, instead of System (\ref{eq:UV}), one with non-degenerate 
leading terms at irregular singular points.
\begin{proposition}
Define variables $\mu$ and $\Psi \! = \! \Psi (\mu,\tau):$
\begin{eqnarray*}
&\lambda \! = \! \mu^{2}, \qquad \Phi (\lambda,\tau) \! := \! \sqrt{\mu} \! 
\left(\!
\begin{pmatrix}
1 & 0 \\
0 & 0
\end{pmatrix} \! + \! \dfrac{1}{\mu} \! 
\begin{pmatrix}
0 & - \frac{A}{\sqrt{\smash[b]{-AB}}} \\
0 & 1
\end{pmatrix} \! \right) \! \Psi (\mu,\tau),&
\end{eqnarray*}
where $\Phi(\lambda,\tau)$ is a fundamental solution of 
System~{\rm(\ref{eq:UV})}. Then
\begin{equation}
\label{eq:UVwave}
\partial_{\mu} \Psi  \! = \! \widetilde{\mathscr{U}}(\mu,\tau) \Psi, \qquad 
\quad \partial_{\tau} \Psi \! = \! \widetilde{\mathscr{V}}(\mu,\tau) \Psi,
\end{equation}
where
\begin{gather*}
\widetilde{\mathscr{U}}(\mu,\tau) = -2 \mi \tau \mu \sigma_{3} \! + \! 2 \tau 
\begin{pmatrix}
0 & \frac{2 \mi A}{\sqrt{\smash[b]{-AB}}} \\
-D & 0
\end{pmatrix} \! - \! \dfrac{1}{\mu} \! \left(a \mi \! + \! \frac{2 \tau AD}
{\sqrt{\smash[b]{-AB}}} \! + \! \frac{1}{2} \right) \! \sigma_{3} \! + \! 
\frac{1}{\mu^{2}} 
\begin{pmatrix}
0 & \widetilde{\alpha} \\
\mi \tau B & 0
\end{pmatrix}, \\
\widetilde{\mathscr{V}}(\mu,\tau) = -\mi \mu^{2} \sigma_3 \! + \! \mu 
\begin{pmatrix}
0 & \frac{2 \mi A}{\sqrt{\smash[b]{-AB}}} \\
-D & 0
\end{pmatrix} \! + \! 
\left(\frac{a \mi}{2 \tau} \! - \! \frac{AD}{\sqrt{\smash[b]{-AB}}} \right) \! 
\sigma_{3} \! - \! \frac{1}{2 \tau \mu}
\begin{pmatrix}
0 & \widetilde{\alpha} \\
\mi \tau B & 0
\end{pmatrix},
\end{gather*}
with
\begin{equation}
\widetilde{\alpha} \! := \! -\frac{2}{B} \! \left(a\mi\sqrt{-AB}+\!\tau(AD \! 
+ \! BC) \right).
\end{equation}
\end{proposition}

\emph{Proof.} Follows {}from the change of independent and dependent variables 
given in the Proposition and System~(\ref{eq:UV}). \hfill $\square$
\begin{ccc}
\label{lemma:2.1}
\begin{eqnarray}
&\dfrac{1}{\tau} \det 
\begin{pmatrix}
0 & \widetilde{\alpha} \\
\mi \tau B & 0
\end{pmatrix} \! = \! -\mi \widetilde{\alpha} B \! = \! \varepsilon b, \quad 
\varepsilon \! = \! \pm 1,&
\end{eqnarray}
with $\widetilde{\alpha}$ defined by Equation~{\rm (13)}.
\end{ccc}

\emph{Proof.} Substituting the parametrisation for $A$, $B$, $C$, and $D$ in 
terms of $u(\tau)$ and $\varphi (\tau)$ (cf. Proposition~1.2), one deduces 
that $\widetilde{\alpha} B \! = \! \mi \varepsilon u (\varphi^{\prime} \! - 
\! 2a/ \tau)$: the result stated in the Lemma now follows {}from Equation~(6). 
\hfill $\square$

System~(12) has two irregular singular points, $\mu \! = \! \infty$ and $\mu 
\! = \! 0$. For $\delta \! > \! 0$ and $k \! \in \! \mathbb{Z}$, define the 
(sectorial) neighbourhoods $\Omega_{k}^{\infty}$ and $\Omega_{k}^{0}$, 
respectively, of these points:
\begin{align*}
\Omega_{k}^{\infty} \! &:= \! \left\{\mathstrut \mu; \, \vert \mu \vert \! > 
\! \delta^{-1}, \, -\tfrac{\pi}{2} \! + \! \tfrac{\pi k}{2} \! < \! \arg \mu 
\! + \! \tfrac{1}{2} \arg \tau \! < \! \tfrac{\pi}{2} \! + \! \tfrac{\pi k}{2} 
\right\}, \\
\Omega_{k}^{0} \! &:= \! \left\{\mathstrut \mu; \, \vert \mu \vert \! < \! 
\delta, \, -\pi \! + \! \pi k \! < \! \arg \mu \! - \! \tfrac{1}{2} \arg \tau 
\! - \! \tfrac{1}{2} \arg (\varepsilon b) \! < \! \pi \! + \! \pi k \right\}.
\end{align*}
The following Proposition is a direct consequence of general asymptotic 
results for linear ODEs \cite{W,F} and Lemma~\ref{lemma:2.1}.
\begin{proposition}
For $k \! \in \! \mathbb{Z}$, there exist solutions $Y_{k}^{\infty}(\mu)$ and 
$X_{k}^{0}(\mu)$ of System~{\rm(\ref{eq:UVwave})} which are uniquely defined 
by the following asymptotic expansions: 
\begin{eqnarray}
Y_{k}^{\infty}(\mu) \! \! &\underset{\underset{\mu \in \Omega_{k}^{\infty}}{
\mu \to \infty}}{:=}& \! \! \left(\mathrm{I} \! + \! \frac{1}{\mu} \Psi^{(1)} 
\! + \! \frac{1}{\mu^{2}} \Psi^{(2)} \! + \! \dotsb \right) \! \exp \! \left(
-\mi \! \left(\tau \mu^{2} \! + \! (a \! - \! \mi/2) \ln \mu \right) \! 
\sigma_{3} \right), \\
X_{k}^{0}(\mu) \! \! &\underset{\underset{\mu \in \Omega_{k}^{0}}{\mu \to 0}}{
:=}& \! \! \Psi_{0} \! \left(\mathrm{I} \! + \! \mathcal{Z}_{1} \mu \! + \! 
\dotsb \right) \exp \! \left(\! -\frac{\mi \sqrt{\tau \varepsilon b}}{\mu} \, 
\sigma_{3} \right),
\end{eqnarray}
where $\ln \mu \! := \! \ln \vert \mu \vert \! + \! \mi \arg \mu$,
\begin{eqnarray}
&\mathrm{I} \! = \! 
\begin{pmatrix}
1 & 0 \\
0 & 1
\end{pmatrix}, \qquad
\Psi^{(1)} \! = \! 
\left(
\begin{smallmatrix}
0 & \frac{A}{\sqrt{\smash[b]{-AB}}} \\
\frac{D}{2 \mi} & 0
\end{smallmatrix}
\right), \qquad \Psi^{(2)} \! = \!
\left(
\begin{smallmatrix}
\psi^{(2)}_{11} & 0 \\
0 & \psi^{(2)}_{22}
\end{smallmatrix}
\right),& \\
&\psi^{(2)}_{11} \! = \! -\frac{\mi}{2} \! \left(\tau \sqrt{\smash[b]{-AB}}+ 
\! \tau DC \! + \! \frac{AD}{\sqrt{\smash[b]{-AB}}} \right), \qquad \quad 
\psi^{(2)}_{22} \! = \! \frac{\mi \tau}{2} \! \left(\sqrt{\smash[b]{-AB}}+ \! 
CD \right),& \nonumber \\
&\Psi_{0} \! = \! \frac{\mi}{\sqrt{\smash[b]{2}}} \! \left(\! \frac{
(\varepsilon b)^{1/4}}{\tau^{1/4} \sqrt{\smash[b]{B}}} \! \right)^{\sigma_{3}} 
\left(\sigma_{1} \! + \! \sigma_{3} \right), \qquad \quad \mathcal{Z}_{1} \! = 
\! \left(
\begin{smallmatrix}
z_{1}^{(11)} & z_{1}^{(12)} \\
-z_{1}^{(12)} & - z_{1}^{(11)}
\end{smallmatrix}
\right),& \\
&z_{1}^{(11)} \! =\frac{\left(a \mi + \frac{1}{2} + \frac{2 \tau AD}{\sqrt{
\smash[b]{-AB}}} \right)^{2}}{2 \mi \sqrt{\smash[b]{\tau \varepsilon b}}} \! - 
\! \frac{2 \mi \tau^{3/2} \sqrt{\smash[b]{-AB}}}{\sqrt{\smash[b]{\varepsilon 
b}}} \! - \! \frac{D \sqrt{\smash[b]{\tau \varepsilon b}}}{B}, \qquad \quad 
z_{1}^{(12)} \! = \! \frac{\left(a \mi +\frac{1}{2}+\frac{2 \tau A D}{\sqrt{
\smash[b]{-AB}}} \right)}{2 \mi \sqrt{\tau \varepsilon b}}.& \nonumber
\end{eqnarray}
\end{proposition}
\noindent
In Proposition~2.2 and below, we use the Pauli matrices
\begin{equation*}
\sigma_{1} \! = \! 
\begin{pmatrix}
0 & 1 \\
1 & 0
\end{pmatrix}, \qquad
\sigma_{2} \! = \! 
\begin{pmatrix}
0 & -\mi \\
\mi & 0
\end{pmatrix}, \qquad \text{and} \qquad 
\sigma_{3} \! = \! 
\begin{pmatrix}
1 & 0 \\
0 & -1
\end{pmatrix}.
\end{equation*}

The \emph{canonical solutions}, $Y_{k}^{\infty}(\mu)$ and $X_{k}^{0}(\mu)$, 
enable one to define the \emph{Stokes matrices} $S_{k}^{\infty}$ and $S_{k}^{
0}$:
\begin{eqnarray}
&Y_{k+1}^{\infty}(\mu) \! = \! Y_{k}^{\infty}(\mu)S_{k}^{\infty},& \\
&X_{k+1}^{0}(\mu) \! = \! X_{k}^{0}(\mu)S_{k}^{0}.&
\end{eqnarray}
The Stokes matrices are independent of the parameters $\mu$ and $\tau$, and 
have the following structures:
\begin{equation}
S_{2k}^{\infty} \! = \!
\begin{pmatrix}
1 & 0 \\
s_{2k}^{\infty} & 1
\end{pmatrix}, \quad S_{2k+1}^{\infty} \! = \! 
\begin{pmatrix}
1 & s_{2k+1}^{\infty} \\
0 & 1
\end{pmatrix}, \quad S_{2k}^{0} \! = \! 
\begin{pmatrix}
1 & s_{2k}^{0} \\
0 & 1
\end{pmatrix}, \quad S_{2k+1}^{0} \! = \! 
\begin{pmatrix}
1 & 0 \\
s_{2k+1}^{0} & 1
\end{pmatrix}.
\end{equation}
The parameters $s_{n}^{\infty}$ and $s_{n}^{0}$, $n \! \in \! \mathbb{Z}$, are 
called the \emph{Stokes multipliers}. Using Equations~(15) and~(16), one shows 
that
\begin{equation}
Y_{k+4}^{\infty}(\mu \me^{2 \pi \mi}) \! = \! Y_{k}^{\infty}(\mu) \me^{2 \pi 
(a-\mi/2) \sigma_{3}}, \qquad \quad X_{k+2}^{0}(\mu \me^{2 \pi \mi}) \!= \! 
X_{k}^{0}(\mu);
\end{equation}
hence, {}from Equations~(19), (20), and~(22), one arrives at
\begin{equation}
S_{k+4}^{\infty} \! = \! \me^{-2 \pi (a-\mi/2) \sigma_{3}}S_{k}^{\infty} 
\me^{2 \pi (a-\mi/2) \sigma_{3}}, \qquad \quad S_{k+2}^{0} \! = \! S_{k}^{0}.
\end{equation}
Equations~(23) show that the number of independent Stokes multipliers does not 
exceed six, e.g., $s_{0}^{0}$, $s_{1}^0$, $s_{0}^{\infty}$, $s_{1}^{\infty}$, 
$s_{2}^{\infty}$, and $s_{3}^{\infty}$. Furthermore, due to the special 
structure of System~(12), namely, the coefficient matrices of odd (resp., 
even) powers of $\mu$ in $\widetilde{\mathscr{U}}(\mu,\tau)$ are diagonal 
(resp., off-diagonal) and \emph{vice-versa} for $\widetilde{\mathscr{V}}(\mu,
\tau)$, one can deduce further relations between the Stokes multipliers. More 
precisely, as a result of the above-mentioned special structure, one obtains 
the following symmetry reductions for the canonical solutions:
\begin{equation}
Y_{k}^{\infty}(\mu) \! = \! \sigma_{3} Y_{k+2}^{\infty}(\mu \me^{\pi \mi}) 
\sigma_{3} \me^{-\pi (a-\mi/2) \sigma_{3}}, \qquad \quad X_{k}^{0}(\mu) \! = 
\! \sigma_{3}X_{k+1}^{0}(\mu \me^{\pi \mi}) \sigma_{1}.
\end{equation}
Equations~(24) imply the following relations for the Stokes matrices 
(multipliers):
\begin{equation}
S_{k+2}^{\infty} \! = \! \sigma_{3} \me^{-\pi(a-\mi/2) \sigma_{3}}S_{k}^{
\infty} \me^{\pi(a-\mi/2) \sigma_{3}} \sigma_{3}, \qquad \quad S_{k}^{0} = \! 
\sigma_{1}S_{k+1}^{0} \sigma_{1}.
\end{equation}
Equations~(25) reduce the number of independent Stokes multipliers by a factor 
of 2; in particular, all Stokes multipliers can be expressed in terms of $s_{
0}^{0}$, $s_{0}^{\infty}$, $s_{1}^{\infty}$, and the parameter of formal 
monodromy, $a$.

There is one more relation between the Stokes multipliers which follows {}from 
the so-called cyclic relation. Define the monodromy matrices at infinity, 
$M^{\infty}$, and at zero, $M^{0}$, by the following relations:
\begin{gather}
Y_{0}^{\infty}(\mu \me^{-2 \pi \mi}) \! := \! Y_{0}^{\infty}(\mu)M^{\infty}, \\
X_{0}^{0}(\mu \me^{-2 \pi \mi}) \! := \! X_{0}^{0}(\mu)M^{0}.
\end{gather}
Since $Y_{0}^{\infty}(\mu)$ and $X_{0}^{0}(\mu)$ are solutions of System~(12), 
they differ by a right-hand factor $G$,
\begin{equation}
Y_{0}^{\infty}(\mu) \! := \! X_{0}^{0}(\mu)G,
\end{equation}
where $G$ is called the \emph{connection matrix}. As matrices relating 
fundamental solutions of System (12), the monodromy, connection, and Stokes 
matrices are independent of $\mu$ and $\tau$. Furthermore, since $\mathrm{tr}
(\widetilde{\mathscr{U}}(\mu,\tau)) \! = \! \mathrm{tr}(\widetilde{\mathscr{
V}}(\mu,\tau)) \! = \! 0$ (cf. Proposition~2.1, System~(12)), it follows that
\begin{equation}
\det (M^{0}) \! = \! \det (M^{\infty}) \! = \! \det (G) \! = \! 1.
\end{equation}
{}From the definition of the monodromy and connection matrices, one deduces 
the following \emph{cyclic relation}:
\begin{equation}
GM^{\infty} \! = \! M^{0}G.
\end{equation}
The monodromy matrices can be expressed in terms of the Stokes matrices:
\begin{equation}
M^{\infty} \! = \! S_{0}^{\infty}S_{1}^{\infty}S_{2}^{\infty}S_{3}^{\infty} 
\me^{-2 \pi (a-\mi/2) \sigma_{3}}, \qquad \quad M^{0} \! = \! S_{0}^{0}
S_{1}^{0}.
\end{equation}
Due to the symmetry reductions~(24), one can derive, using Equations~(19), 
(20), and~(28), the \emph{semi-cyclic} relation:
\begin{equation}
\label{eq:semicycle}
G^{-1}S_{0}^{0} \sigma_{1}G \! = \! S_{0}^{\infty}S_{1}^{\infty} \sigma_{3} 
\me^{-\pi (a-\mi/2) \sigma_{3}}.
\end{equation}
Using the first of Equations~(25), one shows that $M^{\infty} \! = \! (S_{0}^{
\infty}S_{1}^{\infty} \sigma_{3} \me^{-\pi (a-\mi/2) \sigma_{3}})^{2}$; hence, 
the semi-cyclic relation (Equation~(32)) implies the cyclic one 
(Equation~(30)). The Stokes multipliers, $s_{0}^{0}$, $s_{0}^{\infty}$, and 
$s_{1}^{\infty}$, the elements of the connection matrix, $(G)_{ij} \! =: \! 
g_{ij}$, $i,j \! = \! 1,2$, and the parameter of formal monodromy, $a$, are 
called the \emph{monodromy data}. Consider $\mathbb{C}^{8}$ with co-ordinates 
$(a,s_{0}^{0},s_{0}^{\infty},s_{1}^{\infty},g_{11},g_{12},g_{21},g_{22})$. The 
algebraic variety defined by $\det (G) \! = \! 1$ and Equation~(32) is called 
the \emph{manifold of monodromy data}, $\mathscr{M}$. Since only three of the 
four equations in (\ref{eq:semicycle}) are independent, it is clear that 
$\mathrm{dim}_{\mathbb{C}}(\mathscr{M}) \! = \! 4$. More precisely, the 
equations defining $\mathscr{M}$ are as follows:
\begin{eqnarray}
&s_{0}^{\infty}s_{1}^{\infty} \! = \! -1 \! - \! \me^{-2 \pi a} \! - \! \mi 
s_{0}^{0} \me^{-\pi a}, \quad \qquad g_{22}g_{21} \! - \! g_{11} g_{12} \! + 
\! s_{0}^{0}g_{11}g_{22} \! = \! \mi \me^{-\pi a},& \\
&g_{11}^{2} \! - \! g_{21}^{2} \! - \! s_{0}^{0} g_{11} g_{21} \! = \! \mi 
\me^{-\pi a} s_{0}^{\infty}, \quad \, g_{22}^{2} \! - \! g_{12}^{2} \! + \! 
s_{0}^{0} g_{12} g_{22} \! = \! \mi \me^{\pi a} s_{1}^{\infty}, \quad \, 
g_{11} g_{22} \! - \! g_{12} g_{21} \! = \! 1.& \nonumber
\end{eqnarray}
These equations are equivalent to one of the following (three) systems:
\begin{enumerate}
\item[({\bf 1})] $g_{11}g_{22} \! \not= \! 0 \, \Rightarrow$
\begin{equation*}
s_{0}^{\infty} \! = \! - \dfrac{(g_{21} \! + \! \mi g_{11} \me^{\pi a})}{g_{
22}}, \qquad s_{1}^{\infty} \! = \! \dfrac{(g_{12} \! - \! \mi g_{22} \me^{-
\pi a})}{g_{11}}, \qquad s_{0}^{0} \! = \! \dfrac{\mi \me^{-\pi a}}{g_{11}
g_{22}} \! + \! \dfrac{g_{12}}{g_{22}} \! - \! \dfrac{g_{21}}{g_{11}};
\end{equation*}
\item[({\bf 2})] $g_{11} \! = \! 0 \! \Rightarrow \! g_{22} \! \not= \! 0$, 
$\mathrm{parameters} \! = \! (s_{0}^{0},g_{22})$,
\begin{equation*}
g_{12} \! = \! \mi g_{22} \me^{\pi a}, \quad g_{21} \! = \! \dfrac{\mi \me^{
-\pi a}}{g_{22}}, \quad s_{0}^{\infty} \! = \! -\dfrac{\mi \me^{-\pi a}}{g_{
22}^{2}}, \quad s_{1}^{\infty} \! = \! -\mi g_{22}^{2} \! \left(1 \! + \! 
\me^{2 \pi a} \! + \! \mi s_{0}^{0} \me^{\pi a} \right) \! \me^{-\pi a};
\end{equation*}
\item[({\bf 3})] $g_{11} \! \not= \! 0 \! \Rightarrow \! g_{22} \! = \! 0$, 
$\mathrm{parameters} \! = \! (s_{0}^{0},g_{11})$,
\begin{equation*}
g_{12} \! = \! -\dfrac{\mi \me^{-\pi a}}{g_{11}}, \quad g_{21} \! = \! -\mi 
\me^{\pi a} g_{11}, \quad s_{1}^{\infty} \! = \! -\dfrac{\mi \me^{-3 \pi a}}
{g_{11}^{2}}, \quad s_{0}^{\infty} \! = \! -\mi g_{11}^{2} \! \left(1 \! + 
\! \me^{2 \pi a} \! + \! \mi s_{0}^{0} \me^{\pi a} \right) \! \me^{\pi a}.
\end{equation*}
\end{enumerate}
\begin{eee}
In view of the rescalings mentioned in Remark~1.2, one can renormalise the 
canonical solutions, namely, $\tau \! \to \! \tau /\sqrt{\smash[b]{\varepsilon 
b}}$, $\mu \! \to \! \mu (\varepsilon b)^{1/4}$, and $Y_{k}^{\infty} \! \to 
\! Y_{k}^{\infty} \exp \! \left(-\tfrac{1}{4}(a \! - \! \tfrac{\mi}{2}) \ln 
(\varepsilon b) \sigma_{3} \right)$, such that the monodromy data will be 
independent of the parameters $\varepsilon, b$. Thus, the dependence of the 
solution, $u(\tau)$, on the parameter (product) $\varepsilon b$ is similar to 
its dependence on $\tau$; in particular, $u(\tau)$ has singular points at 
$\varepsilon b \! = \! 0$ and $\infty$. \hfill $\blacksquare$
\end{eee}
\begin{proposition}
Define the functions
\begin{equation}
\mathcal{H}_{0}(\tau) \! := \! \mathrm{tr} \! \left(\tfrac{\mi \sqrt{
\smash[b]{\varepsilon b}}}{2 \sqrt{\smash[b]{\tau}}} \mathcal{Z}_{1} \sigma_{
3} \right) \qquad \mathrm{and} \qquad \mathcal{H}_{\infty}(\tau) \! := \! 
\mathrm{tr} \! \left(\mi \! \left(2 \Psi^{(2)} \! - \! (\Psi^{(1)})^{2} 
\right) \! \sigma_{3} \right),
\end{equation}
where $\Psi^{(1)}$, $\Psi^{(2)}$, and $\mathcal{Z}_{1}$ are the coefficients 
of the asymptotic expansions~{\rm (15)} and~{\rm (16)}. Let
\begin{equation}
\mathcal{H}(\tau) \! := \! \mathcal{H}_{0}(\tau) \! + \! \mathcal{H}_{\infty}
(\tau).
\end{equation}
Then
\begin{equation}
\mathcal{H}_{0}(\tau) \! - \! \mathcal{H}_{\infty}(\tau) \! = \! -\tfrac{1}{2 
\tau} \! \left(a \! - \! \tfrac{\mi}{2} \right)^{2},
\end{equation}
\begin{align}
\begin{split}
\mathcal{H}(\tau) \! &= \! \tfrac{\left(a \mi +\frac{1}{2}+\frac{2 \tau AD}{
\sqrt{\smash[b]{-AB}}} \right)^{2}}{2 \tau} \! + \! 4 \tau \sqrt{\smash[b]{-
AB}} \! - \! \tfrac{\mi \varepsilon bD}{B} \! + \! 2 \tau CD \! + \! \tfrac{A
D}{\sqrt{\smash[b]{-AB}}} \\
 &= \! \left(a \! - \! \tfrac{\mi}{2} \right) \! \tfrac{b}{u(\tau)} \! + \! 
\tfrac{1}{2 \tau} \! \left(a \! - \! \tfrac{\mi}{2} \right)^{2} \! + \! 
\tfrac{\tau}{4u^{2}(\tau)} \! \left((u'(\tau))^2 \! + \! b^{2} \right) \! + \! 
4 \varepsilon u(\tau).
\end{split}
\end{align}
\end{proposition}

\emph{Proof.} Using the explicit expressions for $\Psi^{(1)}$, $\Psi^{(2)}$, 
and $\mathcal{Z}_{1}$ given in Equations~(17) and~(18), one obtains the first 
equation of~(37) and
\begin{equation*}
\mathcal{H}_{0}(\tau) \! - \! \mathcal{H}_{\infty}(\tau) \! = \! \tfrac{\left(
a \mi +\frac{1}{2}+\frac{2 \tau AD}{\sqrt{\smash[b]{-AB}}} \right)^{2}}{2 
\tau} \! - \! \tfrac{\mi \varepsilon bD}{B} \! - \! 2 \tau CD \! - \! \tfrac{A
D}{\sqrt{\smash[b]{-AB}}}.
\end{equation*}
Now, {}from the formulae for $A(\tau)$, $B(\tau)$, $C(\tau)$, and $D(\tau)$ in 
terms of $u(\tau)$ and $\varphi (\tau)$ given in Proposition~1.2, one arrives 
at Equation~(36) and the second of~(37). \hfill $\square$
\begin{eee}
Note that $\mathcal{H}_{\varepsilon_{1}}(p(\tau),u(\tau);\tau) \vert_{
\varepsilon_{1}=-1} \! = \! \mathcal{H}(\tau)$, where $\mathcal{H}_{
\varepsilon_{1}}(p,q;\tau)$ is defined in Equation (7), $\mathcal{H}(\tau)$ 
is given by the second of~(37), and $p(\tau)$ is calculated {}from the first 
equation of~(11) with $q \! = \! u(\tau)$. \hfill $\blacksquare$
\end{eee}

An important object in the theory of Painlev\'e equations is the $\pmb{\pmb{
\boldsymbol{\tau}}}$-function \cite{a11,a12,a23}. In this work, it is defined 
as follows:
\begin{equation*}
\mathcal{H}(\tau) \! := \! \partial_{\tau} \ln (\pmb{\pmb{\boldsymbol{\tau}}}
(\tau)).
\end{equation*}
\section{Summary of Results}
In order to present our asymptotic results, it is convenient to introduce the 
following auxiliary mapping $\mathscr{F}_{\varepsilon_{1},\varepsilon_{2}} 
\colon \mathscr{M} \to \mathscr{M}$, $(a,s^{0}_{0},s^{\infty}_{0},s^{\infty}_{
1},g_{11},g_{12},g_{21},g_{22}) \to ((-1)^{\varepsilon_{2}}
a,s^{0}_{0},s^{\infty}_{0}(\varepsilon_{1},\varepsilon_{2}),s^{\infty}_{1}
(\varepsilon_{1},\linebreak[4]
\varepsilon_{2}),g_{11}(\varepsilon_{1},\varepsilon_{2}),g_{12}(\varepsilon_{
1},\varepsilon_{2}),g_{21}(\varepsilon_{1},\varepsilon_{2}),g_{22}
(\varepsilon_{1},\varepsilon_{2}))$, $\varepsilon_{1},\varepsilon_{2} \! = \! 
0,\pm 1$. Define:
\begin{enumerate}
\item[(1)] $\mathscr{F}_{0,0}$ as the identity mapping: $s^{\infty}_{0}(0,0) 
\! = \! s^{\infty}_{0}$, $s^{\infty}_{1}(0,0) \! = \! s^{\infty}_{1}$, and 
$g_{ij}(0,0) \! = \! g_{ij}$, $i,j \! = \! 1,2$;
\item[(2)] $\mathscr{F}_{0,-1}$ as: $s^{\infty}_{0}(0,-1) \! = \! s^{\infty}_{
1} \me^{-\pi a}$, $s^{\infty}_{1}(0,-1) \! = \! s^{\infty}_{0} \me^{-\pi a}$, 
$g_{11}(0,-1) \! = \! -g_{22} \me^{-\frac{\pi a}{2}}$, $g_{12}(0,-1) \! = \! 
-(g_{21} \! + \! s^{\infty}_{0}g_{22}) \me^{\frac{\pi a}{2}}$, $g_{21}(0,-1) 
\! = \! -(g_{12} \! - \! s^{0}_{0}g_{22}) \me^{-\frac{\pi a}{2}}$, and $g_{22}
(0,-1) \! = \! -(g_{11} \! - \! s^{0}_{0}g_{21} \! + \! (g_{12} \! - \! s^{
0}_{0}g_{22})s^{\infty}_{0}) \me^{\frac{\pi a}{2}}$;
\item[(3)] $\mathscr{F}_{0,1}$ as: $s^{\infty}_{0}(0,1) \! = \! s^{\infty}_{1} 
\me^{-\pi a}$, $s^{\infty}_{1}(0,1) \! = \! s^{\infty}_{0} \me^{-\pi a}$, 
$g_{11}(0,1) \! = \! -\mi g_{12} \me^{-\frac{\pi a}{2}}$, $g_{12}(0,1) \! = \! 
-\mi (g_{11} \! + \! s^{\infty}_{0}g_{12}) \me^{\frac{\pi a}{2}}$, $g_{21}(0,
1) \! = \! -\mi g_{22} \me^{-\frac{\pi a}{2}}$, and $g_{22}(0,1) \! = \! -\mi 
(g_{21} \! + \! s^{\infty}_{0}g_{22}) \me^{\frac{\pi a}{2}}$;
\item[(4)] $\mathscr{F}_{-1,0}$ as: $s^{\infty}_{0}(-1,0) \! = \! -s^{\infty}_{
0} \me^{-\pi a}$, $s^{\infty}_{1}(-1,0) \! = \! -s^{\infty}_{1} \me^{\pi a}$, 
$g_{11}(-1,0) \! = \! g_{21} \me^{-\frac{\pi a}{2}}$, $g_{12}(-1,0) \! = \! 
-g_{22} \me^{\frac{\pi a}{2}}$, $g_{21}(-1,0) \! = \! (g_{11} \! - \! s^{0}_{0}
g_{21}) \me^{-\frac{\pi a}{2}}$, and $g_{22}(-1,0) \! = \! -(g_{12} \! - \! 
s^{0}_{0}g_{22}) \me^{\frac{\pi a}{2}}$;
\item[(5)] $\mathscr{F}_{-1,-1}$ as: $s^{\infty}_{0}(-1,-1) \! = \! 
-s^{\infty}_{1}$, $s^{\infty}_{1}(-1,-1) \! = \! -s^{\infty}_{0} \me^{-2 \pi 
a}$, $g_{11}(-1,-1) \! = \! g_{12} \! - \! s^{0}_{0}g_{22}$, $g_{12}(-1,-1) \! 
= \! -g_{11} \! + \! s^{0}_{0}g_{21} \! - \! (g_{12} \! - \! s^{0}_{0}g_{22})
s^{\infty}_{0}$, $g_{21}(-1,-1) \! = \! g_{22} \! - \! (g_{12} \! - \! s^{0}_{
0}g_{22})s^{0}_{0}$, and $g_{22}(-1,-1) \! = \! -g_{21} \! + \! (g_{11} \! - 
\! s^{0}_{0}g_{21})s^{0}_{0} \! - \! (g_{22} \! - \! (g_{12} \! - \! s^{0}_{0}
g_{22})s^{0}_{0})s^{\infty}_{0}$;
\item[(6)] $\mathscr{F}_{-1,1}$ as: $s^{\infty}_{0}(-1,1) \! = \! -s^{\infty}_{
1}$, $s^{\infty}_{1}(-1,1) \! = \! -s^{\infty}_{0} \me^{-2 \pi a}$, $g_{11}(-1,
1) \! = \! \mi g_{22}$, $g_{12}(-1,1) \! = \! -\mi (g_{21} \! + \! s^{\infty}_{
0}g_{22})$, $g_{21}(-1,1) \! = \! \mi (g_{12} \! - \! s^{0}_{0}g_{22})$, and 
$g_{22}(-1,1) \! = \! -\mi (g_{11} \! - \! s^{0}_{0}g_{21} \! + \! (g_{12} \! 
- \! s^{0}_{0}g_{22})s^{\infty}_{0})$;
\item[(7)] $\mathscr{F}_{1,0}$ as: $s^{\infty}_{0}(1,0) \! = \! -s^{\infty}_{
0} \me^{\pi a}$, $s^{\infty}_{1}(1,0) \! = \! -s^{\infty}_{1} \me^{-\pi a}$, 
$g_{11}(1,0) \! = \! (g_{21} \!+ \! s^{0}_{0}g_{11}) \me^{\frac{\pi a}{2}}$, 
$g_{12}(1,0) \! = \! -(g_{22} \! + \! s^{0}_{0}g_{12}) \me^{-\frac{\pi a}{
2}}$, $g_{21}(1,0) \! = \! g_{11} \me^{\frac{\pi a}{2}}$, and $g_{22}(1,0) \! 
= \! -g_{12} \me^{-\frac{\pi a}{2}}$;
\item[(8)] $\mathscr{F}_{1,-1}$ as: $s^{\infty}_{0}(1,-1) \! = \! -s^{\infty}_{
1} \me^{-2 \pi a}$, $s^{\infty}_{1}(1,-1) \! = \! -s^{\infty}_{0}$, $g_{11}(1,
-1) \! = \! g_{12} \me^{-\pi a}$, $g_{12}(1,-1) \! = \! -(g_{11} \! + \! 
s^{\infty}_{0}g_{12}) \me^{\pi a}$, $g_{21}(1,-1) \! = \! g_{22} \me^{-\pi 
a}$, and $g_{22}(1,-1) \! = \! -(g_{21} \! + \! s^{\infty}_{0}g_{22}) \me^{\pi 
a}$; and
\item[(9)] $\mathscr{F}_{1,1}$ as: $s^{\infty}_{0}(1,1) \! = \! -s^{\infty}_{
1} \me^{-2 \pi a}$, $s^{\infty}_{1}(1,1) \! = \! -s^{\infty}_{0}$, $g_{11}(1,
1) \! = \! \mi (g_{22} \! + \! s^{0}_{0}g_{12}) \me^{-\pi a}$, $g_{12}(1,1) \! 
= \! -\mi (g_{21} \! + \! s^{0}_{0}g_{11} \! + \! (g_{22} \! + \! s^{0}_{0}
g_{12})s^{\infty}_{0}) \me^{\pi a}$, and $g_{22}(1,1) \! = \! -\mi (g_{11} \! 
+ \! s^{\infty}_{0}g_{12}) \me^{\pi a}$.
\end{enumerate}
\begin{eee}
Throughout this work, $\epsilon$ and $\delta$ (assumed sufficiently small), 
with and without subscripts, superscripts, etc., denote positive real numbers, 
and the context(s) in which they appear should make clear the parameter(s), 
if any, on which they depend. The roots of positive quantities are assumed 
positive, whilst the branches of the roots of complex quantities can be taken 
arbitrarily, unless stated otherwise. For negative values of a real variable 
$\varkappa$ $(\varkappa \! = \! \tau$ or $\varepsilon b$, say), $\varkappa^{
1/3} \! := \! -\vert \varkappa \vert^{1/3}$. \hfill $\blacksquare$
\end{eee}
\begin{ddd}
Let $\varepsilon_{1},\varepsilon_{2} \! = \! 0,\pm 1$, $\varepsilon b \! = \! 
\vert \varepsilon b \vert \me^{\mi \pi \varepsilon_{2}}$, and $u(\tau)$ be a 
solution of Equation~{\rm (\ref{eq:dp3})} corresponding to the monodromy 
data $(a,s^{0}_{0},s^{\infty}_{0},
s^{\infty}_{1},g_{11},g_{12},g_{21},g_{22})$. Suppose that
\begin{equation}
g_{11}(\varepsilon_{1},\varepsilon_{2})g_{12}(\varepsilon_{1},\varepsilon_{2})
g_{21}(\varepsilon_{1},\varepsilon_{2})g_{22}(\varepsilon_{1},\varepsilon_{2}) 
\! \not= \! 0, \quad \left \vert \Re \! \left(\tfrac{\mi}{2 \pi} \ln (g_{11}
(\varepsilon_{1},\varepsilon_{2})g_{22}(\varepsilon_{1},\varepsilon_{2})) 
\right) \right \vert \! < \! \dfrac{1}{6}.
\end{equation}
Then $\exists \, \, \delta \! > \! 0$ such that $u(\tau)$ has the asymptotic 
expansion
\begin{align}
u(\tau) \underset{\tau \to \infty \me^{\mi \pi \varepsilon_{1}}}{=}& \, \dfrac{
(-1)^{\varepsilon_{1}} \varepsilon \sqrt{\smash[b]{\vert \varepsilon b \vert}}
}{3^{1/4}} \! \left(\! \sqrt{\dfrac{\vartheta (\tau)}{12}}+ \! \sqrt{
\widetilde{\nu}(\varepsilon_{1},\varepsilon_{2}) \! + \! 1} \, \me^{\frac{3 
\pi \mi}{4}} \cosh \! \left(\mi \vartheta (\tau) \! + \! (\widetilde{\nu}
(\varepsilon_{1},\varepsilon_{2}) \! + \! 1) \ln \vartheta (\tau) \right. 
\right. \nonumber \\
+&\left. \left. \, z(\varepsilon_{1},\varepsilon_{2}) \! + \! o \! \left(
\tau^{-\delta} \right) \right) \right),
\end{align}
where
\begin{equation*}
\vartheta (\tau) \! := \! 3 \sqrt{\smash[b]{3}} \, \vert \varepsilon b \vert^{
1/3} \vert \tau \vert^{2/3}, \qquad \quad \widetilde{\nu}(\varepsilon_{1},
\varepsilon_{2}) \! + \! 1 \! := \! \dfrac{\mi}{2 \pi} \ln (g_{11}
(\varepsilon_{1},\varepsilon_{2})g_{22}(\varepsilon_{1},\varepsilon_{2})),
\end{equation*} 
\begin{align*}
z(\varepsilon_{1},\varepsilon_{2}) :=& \, \dfrac{1}{2} \ln (2 \pi) \! - \! 
\dfrac{\pi \mi}{2} \! - \! \dfrac{3 \pi \mi}{2} \! \left(\widetilde{\nu}
(\varepsilon_{1},\varepsilon_{2}) \! + \! 1 \right) \! + \! (-1)^{\varepsilon_{
2}} \mi a \ln \! \left(2 \! + \! \sqrt{\smash[b]3} \right) \! + \! \left(
\widetilde{\nu}(\varepsilon_{1},\varepsilon_{2}) \! + \! 1 \right) \! \ln 12 \\
-& \, \ln \! \left(\omega (\varepsilon_{1},\varepsilon_{2}) \sqrt{\smash[b]{
\widetilde{\nu}(\varepsilon_{1},\varepsilon_{2}) \! + \! 1}} \, \Gamma 
(\widetilde{\nu}(\varepsilon_{1},\varepsilon_{2}) \! + \! 1) \right),
\end{align*}
with
\begin{equation*}
\omega (\varepsilon_{1},\varepsilon_{2}) \! := \! \dfrac{g_{12}(\varepsilon_{
1},\varepsilon_{2})}{g_{22}(\varepsilon_{1},\varepsilon_{2})},
\end{equation*}
and $\Gamma (\cdot)$ is the gamma function {\rm \cite{a24}}.

Let $\mathcal{H}(\tau)$ be the Hamiltonian function defined in 
Equation~{\rm (35)} corresponding to the function $u(\tau)$ given above. Then
\begin{align}
\mathcal{H}(\tau) \underset{\tau \to \infty \me^{\mi \pi \varepsilon_{1}}}{=}& 
\, 3(\varepsilon b)^{2/3} \tau^{1/3} \! + \! 2 \vert \varepsilon b \vert^{1/3} 
\tau^{-1/3} \! \left(\! \left(a \! - \! (-1)^{\varepsilon_{2}} \mi/2 \right) 
\! - \! 2 \sqrt{\smash[b]{3}} \, \mi \! \left(\widetilde{\nu}(\varepsilon_{1},
\varepsilon_{2}) \! + \! 1 \right) \! + \! o \! \left(\tau^{-\delta} \right) 
\right) \nonumber \\
+& \, \dfrac{(a \! - \! (-1)^{\varepsilon_{2}} \mi/2)^{2}}{2 \tau}.
\end{align}
\end{ddd}
\begin{eee}
If conditions~(38) are valid for two distinct pairs of values of 
$(\varepsilon_{1},\varepsilon_{2})$, $(\varepsilon_{1}^{\flat},\varepsilon_{
2}^{\flat})$ and $(\varepsilon_{1}^{\natural},\varepsilon_{2}^{\natural})$, 
say, then Theorem~3.1 yields connection formulae for asymptotics of $u(\tau)$ 
as $\tau \! \to \! \infty \me^{\mi \pi \varepsilon_{1}^{\flat}}$ and $\tau 
\! \to \! \infty \me^{\mi \pi \varepsilon_{1}^{\natural}}$; for example, the 
connection formulae in terms of the parameters $\omega (\pm 1,0)$ and 
$\widetilde{\nu}(\pm 1,0)$ read:
\begin{gather*}
\omega (1,0) \! = \! \omega \! + \! \left(\mi \me^{-\pi a} \! + \! \dfrac{1}{
\omega} \right) \! \me^{2 \pi \mi (\widetilde{\nu}+1)}, \\
\me^{-2 \pi \mi (\widetilde{\nu}(1,0)+1)} \! = \! -\omega \! \left(\omega 
\me^{-2 \pi \mi (\widetilde{\nu}+1)} \! + \! \mi \me^{-\pi a} \right), \\
\omega (-1,0) \! = \! \dfrac{\omega \me^{2 \pi \mi (\widetilde{\nu}+1)}}{\me^{
2 \pi \mi (\widetilde{\nu}+1)} \! - \! 1 \! - \! \mi \omega \me^{-\pi a}}, \\
\me^{-2 \pi \mi (\widetilde{\nu}(-1,0)+1)} \! = \! \dfrac{1}{\omega^{2}} \! 
\left(1 \! - \! \me^{2 \pi \mi (\widetilde{\nu}+1)} \right) \! \left(1 \! - \! 
\me^{-2 \pi \mi (\widetilde{\nu}+1)} \! + \! \mi \omega \me^{-\pi a} \right),
\end{gather*}
where $\omega \! := \! \omega (0,0)$, and $\widetilde{\nu} \! := \! \widetilde{
\nu}(0,0)$. \hfill $\blacksquare$
\end{eee}
\begin{ddd}
Let $\varepsilon_{1},\varepsilon_{2} \! = \! 0,\pm 1$, $\varepsilon b \! = \! 
\vert \varepsilon b \vert \me^{\mi \pi \varepsilon_{2}}$, and $u(\tau)$ be a 
solution of Equation~{\rm (\ref{eq:dp3})} corresponding to the monodromy 
data $(a,s^{0}_{0},
s^{\infty}_{0},s^{\infty}_{1},g_{11},g_{12},g_{21},g_{22})$. Suppose that
\begin{equation*}
g_{21}(\varepsilon_{1},\varepsilon_{2}) \! = \! 0, \qquad \quad g_{11}
(\varepsilon_{1},\varepsilon_{2})g_{22}(\varepsilon_{1},\varepsilon_{2}) \! = 
\! 1.
\end{equation*}
Then $\exists \, \, \delta \! > \! 0$ such that $u(\tau)$ has the asymptotic 
expansion
\begin{align}
u(\tau) \underset{\tau \to \infty \me^{\mi \pi \varepsilon_{1}}}{=}& \, 
\dfrac{\varepsilon (\varepsilon b)^{2/3}}{2} \tau^{1/3} \! + \! \dfrac{(-1)^{
\varepsilon_{1}} \varepsilon \sqrt{\smash[b]{\vert \varepsilon b \vert}} \, 
(s_{0}^{0} \! - \! \mi \me^{(-1)^{\varepsilon_{2}+1} \pi a})}{2^{3/2}3^{1/4} 
\sqrt{\smash[b]{\pi}}} \! \left(\! \dfrac{\sqrt{\smash[b]{3}}- \! 1}{\sqrt{
\smash[b]{3}}+ \! 1} \right)^{(-1)^{\varepsilon_{2}} \mi a} \nonumber \\
\times& \, \exp \! \left(-\mi \! \left(3 \sqrt{\smash[b]{3}} \, \vert 
\varepsilon b \vert^{1/3} \vert \tau \vert^{2/3} \! - \! \dfrac{\pi}{4} 
\right) \right) \! \left(1 \! + \! o \! \left(\tau^{-\delta} \right) \right).
\end{align}
\end{ddd}
\begin{ddd}
Let $\varepsilon_{1},\varepsilon_{2} \! = \! 0,\pm 1$, $\varepsilon b \! = \! 
\vert \varepsilon b \vert \me^{\mi \pi \varepsilon_{2}}$, and $u(\tau)$ be a 
solution of Equation~{\rm (\ref{eq:dp3})} corresponding to the monodromy 
data $(a,s^{0}_{0},s^{\infty}_{0},
s^{\infty}_{1},g_{11},g_{12},g_{21},g_{22})$. Suppose that
\begin{equation*}
g_{12}(\varepsilon_{1},\varepsilon_{2}) \! = \! 0, \qquad \quad g_{11}
(\varepsilon_{1},\varepsilon_{2})g_{22}(\varepsilon_{1},\varepsilon_{2}) \! = 
\! 1.
\end{equation*}
Then $\exists \, \, \delta \! > \! 0$ such that $u(\tau)$ has the asymptotic 
expansion
\begin{align}
u(\tau) \underset{\tau \to \infty \me^{\mi \pi \varepsilon_{1}}}{=}& \, 
\dfrac{\varepsilon (\varepsilon b)^{2/3}}{2} \tau^{1/3} \! + \! \dfrac{(-1)^{
\varepsilon_{1}} \varepsilon \sqrt{\smash[b]{\vert \varepsilon b \vert}} \, 
(s_{0}^{0} \! - \! \mi \me^{(-1)^{\varepsilon_{2}+1} \pi a})}{2^{3/2}3^{1/4} 
\sqrt{\smash[b]{\pi}}} \! \left(\! \dfrac{\sqrt{\smash[b]{3}}+ \! 1}{\sqrt{
\smash[b]{3}}- \! 1} \right)^{(-1)^{\varepsilon_{2}} \mi a} \nonumber \\
\times& \, \exp \! \left(\mi \! \left(3 \sqrt{\smash[b]{3}} \, \vert 
\varepsilon b \vert^{1/3} \vert \tau \vert^{2/3} \! + \! \dfrac{3 \pi}{4} 
\right) \right) \! \left(1 \! + \! o \! \left(\tau^{-\delta} \right) \right).
\end{align}
\end{ddd}
\begin{eee}
The function $\mathcal{H}(\tau)$ also has asymptotics~(40) for the conditions 
on the monodromy data given in Theorems~3.2 and~3.3. \hfill $\blacksquare$
\end{eee}
\begin{ddd}
Let $\varepsilon_{1},\varepsilon_{2} \! = \! 0,\pm 1$, $\varepsilon b \! = \! 
\vert \varepsilon b \vert \me^{\mi \pi \varepsilon_{2}}$, and $u(\tau)$ be a 
solution of Equation~{\rm (\ref{eq:dp3})} corresponding to the monodromy 
data $(a,s^{0}_{0},s^{\infty}_{0},
s^{\infty}_{1},g_{11},g_{12},g_{21},g_{22})$. Suppose that
\begin{equation}
\vert \Im (a) \vert \! < \! 1, \qquad \quad g_{11}(\varepsilon_{1},
\varepsilon_{2})g_{22}(\varepsilon_{1},\varepsilon_{2}) \! \not= \! 0, \qquad 
\quad \rho \! \not= \! 0, \qquad \quad \vert \Re (\rho) \vert \! < \! \dfrac{
1}{2},
\end{equation}
where
\begin{eqnarray}
\cos (2 \pi \rho) \! := \! -\dfrac{\mi s^{0}_{0}}{2} \! = \cosh (\pi a) \! + 
\! \dfrac{1}{2} s_{0}^{\infty}s_{1}^{\infty} \me^{\pi a}.
\end{eqnarray}
Then $\exists \, \, \delta \! > \! 0$ such that $u(\tau)$ has the asymptotic 
expansion
\begin{align}
u(\tau) \underset{\tau \to 0 \me^{\mi \pi \varepsilon_{1}}}{=}& \, \dfrac{
(-1)^{\varepsilon_{2}} \tau b}{16 \pi} \exp \! \left((-1)^{\varepsilon_{2}} 
\dfrac{\pi a}{2} \right) \! \left(\mathfrak{p}((-1)^{\varepsilon_{2}}a,\rho) 
\chi_{1}(\vec{g}(\varepsilon_{1},\varepsilon_{2});\rho) \vert \tau \vert^{2 
\rho} \! + \! \mathfrak{p}((-1)^{\varepsilon_{2}}a,-\rho) \right. \nonumber \\
\times&\left. \, \chi_{1}(\vec{g}(\varepsilon_{1},\varepsilon_{2});-\rho) 
\vert \tau \vert^{-2 \rho} \right) \! \left(\mathfrak{p}((-1)^{\varepsilon_{2}
+1}a,\rho) \me^{-\mi \pi \rho} \chi_{2}(\vec{g}(\varepsilon_{1},\varepsilon_{
2});\rho) \vert \tau \vert^{2 \rho} \right. \nonumber \\
+&\left. \, \mathfrak{p}((-1)^{\varepsilon_{2}+1}a,-\rho) \me^{\mi \pi \rho} 
\chi_{2}(\vec{g}(\varepsilon_{1},\varepsilon_{2});-\rho) \vert \tau \vert^{-
2 \rho} \right) \! \left(1 \! + \! \mathcal{O} \! \left(\tau^{\delta} \right) 
\right),
\end{align}
where
\begin{equation}
\mathfrak{p}(z_{1},z_{2}) \! := \! \left(\dfrac{\vert \varepsilon b \vert}{32} 
\me^{\frac{\mi \pi}{2}} \right)^{z_{2}} \! \left(\dfrac{\Gamma (\frac{1}{2} \! 
- \! z_{2})}{\Gamma (1 \! + \! z_{2})} \right)^{2} \dfrac{\Gamma (1 \! + \! 
z_{2} \! + \! \frac{\mi z_{1}}{2})}{\tan (\pi z_{2})},
\end{equation}
\begin{equation}
\begin{gathered}
\chi_{1}(\vec{g}(\varepsilon_{1},\varepsilon_{2});z_{3}) \! := \! g_{11}
(\varepsilon_{1},\varepsilon_{2}) \me^{\mi \pi z_{3}} \me^{\frac{\mi \pi}{4}} 
\! + \! g_{21}(\varepsilon_{1},\varepsilon_{2}) \me^{-\mi \pi z_{3}} \me^{-
\frac{\mi\pi}{4}}, \\
\chi_{2}(\vec{g}(\varepsilon_{1},\varepsilon_{2});z_{4}) \! := \! g_{12}
(\varepsilon_{1},\varepsilon_{2}) \me^{\mi \pi z_{4}} \me^{\frac{\mi \pi}{4}} 
\! + \! g_{22}(\varepsilon_{1},\varepsilon_{2}) \me^{-\mi \pi z_{4}} \me^{-
\frac{\mi \pi}{4}}.
\end{gathered}
\end{equation}

Let $\mathcal{H}(\tau)$ be the Hamiltonian function defined in 
Equation~{\rm (35)} corresponding to the function $u(\tau)$ given above. Then
\begin{align}
\mathcal{H}(\tau) \underset{\tau \to 0 \me^{\mi \pi \varepsilon_{1}}}{=}& \, 
\dfrac{2 \rho}{\tau} \dfrac{\left(\mathfrak{p}((-1)^{\varepsilon_{2}}a,\rho) 
\chi_{1}(\vec{g}(\varepsilon_{1},\varepsilon_{2});\rho) \vert \tau \vert^{2 
\rho} \! - \! \mathfrak{p}((-1)^{\varepsilon_{2}}a,-\rho) \chi_{1}(\vec{g}
(\varepsilon_{1},\varepsilon_{2});-\rho) \vert \tau \vert^{-2 \rho} \right)}{
\left(\mathfrak{p}((-1)^{\varepsilon_{2}},\rho) \chi_{1}(\vec{g}(\varepsilon_{
1},\varepsilon_{2});\rho) \vert \tau \vert^{2 \rho} \! + \! \mathfrak{p}((-1)^{
\varepsilon_{2}}a,-\rho) \chi_{1}(\vec{g}(\varepsilon_{1},\varepsilon_{2});-
\rho) \vert \tau \vert^{-2 \rho} \right)} \nonumber \\
+& \, \dfrac{1}{2 \tau} \! \left(a \! \left(a \! - \! (-1)^{\varepsilon_{2}} 
\mi \right) \! + \! \dfrac{1}{4} \! + \! 8 \rho^{2} \right) \! + \! o \! 
\left(\dfrac{1}{\tau} 
\right).
\end{align}
\end{ddd}
\begin{ddd}
Let $\varepsilon_{1},\varepsilon_{2} \! = \! 0,\pm 1$, $\varepsilon b \! = \! 
\vert \varepsilon b \vert \me^{\mi \pi \varepsilon_{2}}$, and $u(\tau)$ be a 
solution of Equation~{\rm (\ref{eq:dp3})} corresponding to the monodromy 
data $(a,s^{0}_{0},s^{\infty}_{0},
s^{\infty}_{1},g_{11},g_{12},g_{21},g_{22})$. Suppose that
\begin{equation}
\vert \Im (a) \vert \! < \! 1, \qquad \quad g_{11}(\varepsilon_{1},
\varepsilon_{2})g_{22}(\varepsilon_{1},\varepsilon_{2}) \! \not= \! 0, \qquad 
\quad s^{0}_{0} \! = \! 2 \mi.
\end{equation}
Then $\exists \, \, \delta \! > \! 0$ such that $u(\tau)$ has the asymptotic 
expansion
\begin{align}
u(\tau) \underset{\tau \to 0 \me^{\mi \pi \varepsilon_{1}}}{=}& \, \dfrac{
(-1)^{\varepsilon_{2}} \tau b \exp ((-1)^{\varepsilon_{2}} \frac{\pi a}{2})}{
2a \sinh (\frac{\pi a}{2})} \! \left(\chi_{1}(\vec{g}(\varepsilon_{1},
\varepsilon_{2});0) \! \left(1 \! - \! \tfrac{(-1)^{\varepsilon_{2}} \mi a}{2}
Q((-1)^{\varepsilon_{2}}a) \right) \! + \! \dfrac{(-1)^{\varepsilon_{2}} \pi 
a}{4} \right. \nonumber \\
\times&\left. \, (g_{21}(\varepsilon_{1},\varepsilon_{2}) \me^{-\frac{\mi \pi}{
4}} \! - \! 3g_{11}(\varepsilon_{1},\varepsilon_{2}) \me^{\frac{\mi \pi}{4}}) 
\! + \! (-1)^{\varepsilon_{2}} \mi a \chi_{1}(\vec{g}(\varepsilon_{1},
\varepsilon_{2});0) \ln \vert \tau \vert \right) \nonumber \\
\times&\, \left(\chi_{2}(\vec{g}(\varepsilon_{1},\varepsilon_{2});0) \! \left(
1 \! + \! \tfrac{(-1)^{\varepsilon_{2}} \mi a}{2}Q((-1)^{\varepsilon_{2}+1}a) 
\right) \! + \! \dfrac{(-1)^{\varepsilon_{2}} \pi a}{4}(g_{12}(\varepsilon_{1},
\varepsilon_{2}) \me^{\frac{\mi \pi}{4}} \right. \nonumber \\
-&\left. \, 3g_{22}(\varepsilon_{1},\varepsilon_{2}) \me^{-\frac{\mi \pi}{4}}) 
\! - \! (-1)^{\varepsilon_{2}} \mi a \chi_{2}(\vec{g}(\varepsilon_{1},
\varepsilon_{2});0) \ln \vert \tau \vert \right) \! \left(1 \! + \! 
\mathcal{O} \! \left(\tau^{\delta} \right) \right),
\end{align}
where $\chi_{j}(\vec{g}(\varepsilon_{1},\varepsilon_{2});\boldsymbol{\cdot})$, 
$j \! = \! 1,2$, are defined in Theorem~{\rm 3.4}, Equations~{\rm (47)},
\begin{equation}
Q(z) \! := \! 4 \psi (1) \! - \! \psi (\mi z/2) \! + \! \ln 2 \! - \! \ln 
(\vert \varepsilon b \vert),
\end{equation}
$\psi (x) \! := \! \tfrac{\md}{\md x} \ln \Gamma (x)$ is the psi function, and 
$\psi (1) \! = \! -0.57721566490 \dotsc$ {\rm \cite{a24}}.

Let $\mathcal{H}(\tau)$ be the Hamiltonian function defined in 
Equation~{\rm (35)} corresponding to the function $u(\tau)$ given above. 
Then
\begin{equation}
\mathcal{H}(\tau) \! \underset{\tau \to 0 \me^{\mi \pi \varepsilon_{1}}}{=} 
\! \dfrac{1}{2 \tau} \! \left(a \! \left(a \! - \! (-1)^{\varepsilon_{2}} \mi 
\right) \! + \! \dfrac{1}{4} \right) \! + \! \dfrac{b_{2}(\varepsilon_{1},
\varepsilon_{2})}{\tau (a_{2}(\varepsilon_{1},\varepsilon_{2}) \! + \! b_{2}
(\varepsilon_{1},\varepsilon_{2}) \ln \vert \tau \vert)} \! + \! o \! \left(
\dfrac{1}{\tau} \right),
\end{equation}
where
\begin{align*}
a_{2}(\varepsilon_{1},\varepsilon_{2}) :=& \, \chi_{1}(\vec{g}(\varepsilon_{
1},\varepsilon_{2});0) \! \left(1 \! - \! \tfrac{(-1)^{\varepsilon_{2}} \mi a}{
2}Q((-1)^{\varepsilon_{2}}a) \right) \! + \! \dfrac{(-1)^{\varepsilon_{2}} 
\pi a}{4}(g_{21}(\varepsilon_{1},\varepsilon_{2}) \me^{-\frac{\mi \pi}{4}} 
\nonumber \\
-& \, 3g_{11}(\varepsilon_{1},\varepsilon_{2}) \me^{\frac{\mi \pi}{4}}),
\end{align*}
\begin{equation*}
b_{2}(\varepsilon_{1},\varepsilon_{2}) \! := \! (-1)^{\varepsilon_{2}} \mi a 
\chi_{1}(\vec{g}(\varepsilon_{1},\varepsilon_{2});0).
\end{equation*}
\end{ddd}
\begin{hhh}
For the conditions stated in Theorem~{\rm 3.4},
\begin{align}
\pmb{\pmb{\boldsymbol{\tau}}}(\tau) \underset{\tau \to 0 \me^{\mi \pi 
\varepsilon_{1}}}{=}& \, \mathrm{const.} \, \tau^{\frac{1}{2} \! \left(a(a-
(-1)^{\varepsilon_{2}} \mi)+\frac{1}{4}+8\rho^{2} \right)} \! \left(
\mathfrak{p}((-1)^{\varepsilon_{2}}a,\rho) \chi_{1}(\vec{g}(\varepsilon_{1},
\varepsilon_{2});\rho) \vert \tau \vert^{2 \rho} \right. \nonumber \\
+&\left. \, \mathfrak{p}((-1)^{\varepsilon_{2}}a,-\rho) \chi_{1}(\vec{g}
(\varepsilon_{1},\varepsilon_{2});-\rho) \vert \tau \vert^{-2 \rho} \right) \! 
\left(1 \! + \! o \! \left(\tau^{\delta} \right) \right).
\end{align}
\end{hhh}
\begin{hhh}
For the conditions stated in Theorem~{\rm 3.5},
\begin{equation}
\pmb{\pmb{\boldsymbol{\tau}}}(\tau) \! \underset{\tau \to 0 \me^{\mi \pi 
\varepsilon_{1}}}{=} \! \mathrm{const.} \, \tau^{\frac{1}{2} \! \left(a(a-
(-1)^{\varepsilon_{2}} \mi)+\frac{1}{4} \right)} \! \left(a_{2}(\varepsilon_{
1},\varepsilon_{2}) \! + \! b_{2}(\varepsilon_{1},\varepsilon_{2}) \ln \vert 
\tau \vert \right) \! \left(1 \! + \! o \! \left(\tau^{\delta} \right) \right).
\end{equation}
\end{hhh}
\begin{eee}
Conjectures~3.1 and~3.2 can be proved under the assumption that the error 
terms in Equations~(48) and~(52) behave as power-like functions under 
integration. One can also obtain asymptotics for $\pmb{\pmb{\boldsymbol{\tau}}}
(\tau)$ as $\tau \! \to \! \infty \me^{\mi \pi \varepsilon_{1}}$, 
$\varepsilon_{1} \! = \! 0,\pm 1$, up to $\exp (\mathcal{O}(\tau^{2/3}))$, 
under similar assumptions on the error term in Equation~(40). \hfill 
$\blacksquare$
\end{eee}
\begin{eee}
The asymptotic results as $\tau \! \to \! +\infty$ for the functions 
$A(\tau)$, $B(\tau)$, $C(\tau)$, and $D(\tau)$, which solve 
System~(\ref{sys:ABCD}), can be found, via Equations~(56), in 
Proposition~4.3.1. The asymptotics as $\tau \! \to \! +0$ for these 
functions for $\rho \! \not= \! 0$ (resp., $\rho 
\! = \! 0)$ are given in Proposition~5.5 (resp., Proposition~5.7) via 
Proposition~5.6. Asymptotics of these functions for negative and (pure) 
imaginary values of $\tau$ can be obtained by applying the transformations 
given in Section~6, as was done for $u(\tau)$ (see Subsection~6.2 and the 
Appendix). \hfill $\blacksquare$
\end{eee}
\section{Calculation of the Monodromy Data as $\tau \! \to \! +\infty$}
In this section the results stated in Theorem~3.1 for $\tau \! \to \! +\infty$ 
and $\varepsilon b \! > \! 0$ are derived. The corresponding results stated in 
Theorems~3.2 and~3.3 can be obtained analogously. The derivation is based on 
the WKB analysis of the $\mu$-part of System~(12):
\begin{eqnarray}
&\partial_{\mu} \Psi \! = \! \widetilde{\mathscr{U}}(\mu,\tau) \Psi.&
\end{eqnarray}
\subsection{WKB Analysis}
This subsection is devoted to the WKB analysis of Equation~(55) as $\tau \! 
\to \! +\infty$. In order to put Equation~(55) into a form suitable for 
WKB analysis, it is convenient to introduce the notation given in 
Proposition~4.1.1 below.
\begin{proposition:s}
Let
\begin{equation}
\begin{gathered}
A(\tau) \! = \! a(\tau) \tau^{-2/3}, \qquad B(\tau) \! = \! b(\tau) \tau^{-
2/3}, \qquad C(\tau) \! = \! c(\tau) \tau^{-1/3}, \qquad D(\tau) \! = \! 
d(\tau) \tau^{-1/3}, \\
\widetilde{\mu} \! = \! \mu \tau^{1/6}, \qquad \widetilde{\Psi}(\widetilde{
\mu}) \! = \! \tau^{-(1/12) \sigma_{3}} \Psi(\widetilde{\mu} \tau^{-1/6}).
\end{gathered}
\end{equation}
Then
\begin{eqnarray}
&\partial_{\widetilde{\mu}} \widetilde{\Psi}(\widetilde{\mu}) \! = \! \tau^{
2/3}{\cal A}(\widetilde{\mu},\tau) \widetilde{\Psi}(\widetilde{\mu}),&
\end{eqnarray}
where
\begin{equation}
{\cal A}(\widetilde{\mu},\tau) \! := \! -2 \mi \widetilde{\mu} \sigma_{3} \! 
+ \! 
\begin{pmatrix}
0 & -\frac{4 \mi \sqrt{\smash[b]{-a(\tau)b(\tau)}}}{b(\tau)} \\
-2d(\tau) & 0
\end{pmatrix} \! - \! \dfrac{\mi r(\tau)(\varepsilon b)^{1/3}}{2 
\widetilde{\mu}} \sigma_{3} \! + \! \dfrac{1}{\widetilde{\mu}^{2}} 
\begin{pmatrix}
0 & \frac{\mi \varepsilon b}{b(\tau)} \\
\mi b(\tau) & 0
\end{pmatrix},
\end{equation}
with
\begin{eqnarray}
\dfrac{\mi r(\tau)(\varepsilon b)^{1/3}}{2} \! := \! \left(a \mi \! + \! 
\frac{1}{2} \right) \! \tau^{-2/3} \! + \! \dfrac{2a(\tau)d(\tau)}{\sqrt{
\smash[b]{-a(\tau)b(\tau)}}}.
\end{eqnarray}
\end{proposition:s}

\emph{Proof.} Equations~(57)--(59) are obtained by the direct substitution of 
Equations~(56) into Equation~(55), upon taking into account Lemma~2.1. \hfill 
$\square$

The WKB analysis for Equation~(57) is carried out under the following 
conditions:
\begin{gather}
a(\tau)d(\tau) \! + \! b(\tau)c(\tau) \! + \! a \mi \sqrt{\smash[b]{-a(\tau)
b(\tau)}} \, \tau^{-2/3} \! = \! -\dfrac{\mi \varepsilon b}{2}, \\
r(\tau) \! = \! -2 \! + \! r_{0}(\tau) \tau^{-1/3}, \\
\sqrt{\smash[b]{-a(\tau)b(\tau)}} + \! c(\tau)d(\tau) \! + \! \dfrac{a(\tau)
d(\tau) \tau^{-2/3}}{2 \sqrt{\smash[b]{-a(\tau)b(\tau)}}} \! - \! \dfrac{(a \! 
- \! \mi/2)^{2} \tau^{-4/3}}{4} \! = \! \dfrac{3(\varepsilon b)^{2/3}}{4} \! - 
\! h_{0}(\tau) \tau^{-2/3}, \\
\dfrac{2 \sqrt{\smash[b]{-a(\tau)b(\tau)}}}{(\varepsilon b)^{2/3}} \! = \! 1 
\! + \! u_{0}(\tau) \tau^{-1/3},
\end{gather}
where
\begin{equation}
h_{0}(\tau) \! \underset{\tau \to +\infty}{=} \! \mathcal{O}(1), \quad \, 
\vert r_{0}(\tau) \vert \! \underset{\tau \to +\infty}{\leqslant} \! \mathcal{
O} \! \left(\tau^{\delta_{\diamondsuit}} \right), \, \quad \vert u_{0}(\tau) 
\vert \! \underset{\tau \to +\infty}{\leqslant}\! \mathcal{O} \! \left(\tau^{
\delta_{\diamondsuit}} \right), \quad \, 0 \! < \! \delta_{\diamondsuit} \! < 
\! \dfrac{1}{9}.
\end{equation}
\begin{eeee}
Equation~(60) is the integral of motion for System~(\ref{sys:ABCD}) in terms 
of $a(\tau)$, $b(\tau)$, $c(\tau)$, and $d(\tau)$. Actually, in 
conditions~(61)--(63) the functions $r_{0}(\tau)$, $h_{0}(\tau)$, and 
$u_{0}(\tau)$, with the power-like (growth) behaviours given in 
Equation~(64), are introduced. Conditions~(61)--(63), although somewhat 
artificial looking at this stage, will be clarified as they appear in the 
following (asymptotic) analysis. It is worth noting that 
conditions~(60)--(63) are self-consistent; in fact, they are equivalent to
\begin{align}
\dfrac{a(\tau)d(\tau)}{\sqrt{\smash[b]{-a(\tau)b(\tau)}}} &= -\dfrac{\mi 
(\varepsilon b)^{1/3}}{2} \! + \! \dfrac{\mi (\varepsilon b)^{1/3}r_{0}(\tau)
}{4} \tau^{-1/3} \! - \! \dfrac{\mi}{2} \! \left(a \! - \! \dfrac{\mi}{2} 
\right) \tau^{-2/3}, \\
\dfrac{b(\tau)c(\tau)}{\sqrt{\smash[b]{-a(\tau)b(\tau)}}} &= -\dfrac{\mi 
(\varepsilon b)^{1/3}}{2} \! + \! \mi (\varepsilon b)^{1/3} \! \left(\! \dfrac{
u_{0}(\tau)}{1 \! + \! u_{0}(\tau) \tau^{-1/3}} \! - \! \dfrac{r_{0}(\tau)}{4} 
\right) \! \tau^{-1/3} \! - \! \dfrac{\mi}{2} \! \left(a \! + \! \dfrac{\mi}{
2} \right) \! \tau^{-2/3}, \\
-h_{0}(\tau) &= \dfrac{(\varepsilon b)^{2/3}}{2} \! \left(\! \dfrac{u_{0}^{2}
(\tau)}{1 \! + \! u_{0}(\tau) \tau^{-1/3}} \! + \! \dfrac{r_{0}(\tau)}{2} \! 
\left(\! \dfrac{u_{0}(\tau)}{1 \! + \! u_{0}(\tau) \tau^{-1/3}} \! - \! 
\dfrac{r_{0}(\tau)}{4} \right) \right) \nonumber \\
&+ \dfrac{(\varepsilon b)^{1/3}}{2} \! \left(a \! - \! \dfrac{\mi}{2} \right) 
\! - \! \dfrac{(\varepsilon b)^{1/3}(a \! - \! \mi/2)u_{0}(\tau)}{2(1 \! + 
\! u_{0}(\tau) \tau^{-1/3})} \tau^{-1/3}. \qquad \qquad \qquad \qquad \quad 
\blacksquare
\end{align}
\end{eeee}

In certain domains (see below) of the complex $\widetilde{\mu}$-plane, the 
leading term of asymptotics (as $\tau \! \to \! +\infty)$ of a fundamental 
solution, $\widetilde{\Psi}(\widetilde{\mu})$, of Equation~(57) is given by 
the following WKB formula \cite{F},
\begin{equation}
\widetilde{\Psi}_{\textrm{WKB}}(\widetilde{\mu}) \! = \! T(\widetilde{\mu}) 
\exp \! \left( \! -\mi \tau^{2/3} \sigma_{3} \int_{}^{\widetilde{\mu}}l(\xi) 
\, \md \xi \! - \! \int_{}^{\widetilde{\mu}}{\rm diag\/} \!\left( (T(\xi))^{-
1} \partial_{\xi} T(\xi) \right) \md \xi \right),
\end{equation}
where
\begin{equation}
l(\widetilde{\mu}) \! = \! l(\widetilde{\mu},\tau) \! := \! (\det (\mathcal{A}
(\widetilde{\mu})))^{1/2},
\end{equation}
and $T(\widetilde{\mu})$ diagonalises $\mathcal{A}(\widetilde{\mu})$, that is, 
$(T(\widetilde{\mu}))^{-1} \mathcal{A}(\widetilde{\mu})T(\widetilde{\mu}) \! = 
\! -\mi l(\widetilde{\mu}) \sigma_{3}$ (the $\tau$ dependencies of $\widetilde{
\Psi}_{\textrm{WKB}}(\widetilde{\mu})$, $l(\widetilde{\mu})$, and 
$T(\widetilde{\mu})$ have been suppressed in order to simplify the notation). 
It is convenient to choose $T(\widetilde{\mu})$ as follows,
\begin{equation}
T(\widetilde{\mu}) \! = \! \dfrac{\mi}{\sqrt{\smash[b]{2 \mi l(\widetilde{\mu})
({\cal A}_{11}(\widetilde{\mu}) \! - \! \mi l(\widetilde{\mu}))}}} \! \left(
{\cal A}(\widetilde{\mu}) \! - \! \mi l(\widetilde{\mu}) \sigma_{3}\right) \! 
\sigma_{3},
\end{equation}
where $\mathcal{A}_{11}(\widetilde{\mu})$ is the corresponding element of 
matrix $\mathcal{A}(\widetilde{\mu})$. Note that $\det (T(\widetilde{\mu})) \! 
= \! 1$ and $T(\widetilde{\mu}) \! =_{\widetilde{\mu} \to \infty} \! \mathrm{
I} \! + \! \mathcal{O}(\digamma (\tau) \widetilde{\mu}^{-1})$; therefore, 
$\det (\widetilde{\Psi}_{\textrm{WKB}}(\widetilde{\mu})) \! = \! 1$ (since 
$\mathrm{tr}((T(\widetilde{\mu}))^{-1} \partial_{\widetilde{\mu}}T(\widetilde{
\mu})) \! = \! 0)$. The domains in the complex $\widetilde{\mu}$-plane where 
Equation~(68) gives the asymptotic approximation of solutions of Equation~(57) 
are defined in terms of the \emph{Stokes graph}. The vertices of the Stokes 
graph are the singular points of Equation~(57), namely, $\widetilde{\mu} \! = 
\! 0$ and $\infty$, and the \emph{turning points}, which are the roots of the 
equation $l^{2}(\widetilde{\mu}) \! = \! 0$. The edges of the Stokes graph 
are the \emph{Stokes curves}, defined as $\Re (\int_{\widetilde{\mu}^{\mathrm{
TP}}}^{\widetilde{\mu}}l(\xi) \, \md \xi) \! = \! 0$, where $\widetilde{\mu}^{
\mathrm{TP}}$ denotes a turning point. \emph{Canonical domains} are those 
domains in the complex $\widetilde{\mu}$-plane containing one, and only one, 
Stokes curve and bounded by two adjacent Stokes curves. Note that the 
restriction of any branch of $l(\widetilde{\mu})$ to a canonical domain is a 
single-valued function. In each canonical domain, for any choice of the branch 
of $l(\widetilde{\mu})$, there exists a fundamental solution of Equation~(57) 
which has asymptotics whose leading term as $\tau \! \to \! +\infty$ is 
given by Equation~(68). Recalling the definition of $l(\widetilde{\mu})$ 
(Equation~(69)), one shows that
\begin{equation}
l^{2}(\widetilde{\mu}) \! = \! \dfrac{4}{\widetilde{\mu}^{4}} \! \left(\! 
\left(\widetilde{\mu}^{2} \! - \! \frac{1}{2}(\varepsilon b)^{1/3} \right)^{2} 
\! \left(\widetilde{\mu}^{2} \! + \! (\varepsilon b)^{1/3} \right) \! + \! 
\left(\widetilde{\mu}^{4} \! \left(a \! - \! \frac{\mi}{2} \right) \! + \! 
\widetilde{\mu}^{2}h_{0}(\tau) \right) \! \tau^{-2/3} \right).
\end{equation}
It follows {}from Equation~(71) that there are six turning points: two turning 
points coalesce (as $\tau \! \to \! +\infty)$ at $\tfrac{(\varepsilon b)^{1/6}
}{\sqrt{2}} \! + \! \mathcal{O}(\tau^{-1/3})$, another pair coalesce (as $\tau 
\! \to \! +\infty)$ at $-\tfrac{(\varepsilon b)^{1/6}}{\sqrt{2}} \! + \! 
\mathcal{O}(\tau^{-1/3})$, and the remaining two turning points approach (as 
$\tau \! \to \! +\infty)$, respectively, $\pm \mi (\varepsilon b)^{1/6} \! + 
\! \mathcal{O}(\tau^{-4/9})$. Denote by $\widetilde{\mu}_{1}$ the turning 
point in the first quadrant of the complex $\widetilde{\mu}$-plane which 
approaches $\tfrac{(\varepsilon b)^{1/6}}{\sqrt{2}} \! + \! \mathcal{O}
(\tau^{-1/3})$ as $\tau \! \to \! +\infty$, and by $\widetilde{\mu}_{2}$ the 
pure imaginary turning point which approaches $\mi (\varepsilon b)^{1/6} \! + 
\! \mathcal{O}(\tau^{-4/9})$ as $\tau \! \to \! +\infty$. Denote by $\mathscr{
G}_{1}$ the part of the Stokes graph in quadrant one of the complex 
$\widetilde{\mu}$-plane which consists of the vertices $0$, $\infty$, 
$\widetilde{\mu}_{1}$, and $\widetilde{\mu}_{2}$, and the edges $(+\mi \infty,
\widetilde{\mu}_{2})$, $(0,\widetilde{\mu}_{1})$, $(\widetilde{\mu}_{2},
\widetilde{\mu}_{1})$, and $(\widetilde{\mu}_{1},+\infty)$: the complete 
Stokes graph is the union of the mirror images of $\mathscr{G}_{1}$ with 
respect to the real and imaginary axes of the complex $\widetilde{\mu}$-plane.
\begin{proposition:s}
\begin{align}
\int_{\widetilde{\mu}^{\mathrm{TP}}}^{\widetilde{\mu}}l(\xi) \, \md \xi 
\underset{\tau \to +\infty}{=}& \, \left(\! \left(\xi \! + \! \dfrac{
(\varepsilon b)^{1/3}}{\xi} \right) \! \sqrt{\smash[b]{\xi^{2} \! + \! 
(\varepsilon b)^{1/3}}}+ \! \tau^{-2/3} \! \left(\! \left(a \! - \! \dfrac{
\mi}{2} \right) \! \ln \! \left(2 \sqrt{\smash[b]{\xi^{2} \! + \! (\varepsilon 
b)^{1/3}}}+ \! 2 \xi \right) \right. \right. \nonumber \\
+&\left. \left. \dfrac{1}{2 \sqrt{\smash[b]{3}}} \! \left(\! \left(a \! - 
\! \dfrac{\mi}{2} \right) \! + \! \dfrac{2h_{0}(\tau)}{(\varepsilon b)^{1/3}} 
\right) \! \ln \! \left(\! \tfrac{(-\xi +\sqrt{\smash[b]{2}} \, (\varepsilon 
b)^{1/6}+\sqrt{\smash[b]{3}} \, \sqrt{\smash[b]{\xi^{2}+(\varepsilon b)^{1/3}
}} \,)(\xi -\frac{(\varepsilon b)^{1/6}}{\sqrt{\smash[b]{2}}})}{(\xi +\sqrt{
\smash[b]{2}} \, (\varepsilon b)^{1/6}+\sqrt{\smash[b]{3}} \, \sqrt{\smash[b]{
\xi^{2}+(\varepsilon b)^{1/3}}} \,)(\xi +\frac{(\varepsilon b)^{1/6}}{\sqrt{
\smash[b]{2}}})} \right) \! \right) \right. \nonumber \\
+&\left. \left. \mathcal{O} \! \left(\dfrac{\tau^{-4/3}}{\left(\xi \! \pm \! 
\frac{1}{\sqrt{\smash[b]{2}}}(\varepsilon b)^{1/6} \right)^{2}} \right) \! + 
\! \mathcal{O} \! \left(\dfrac{\tau^{-4/3}}{\sqrt{\smash[b]{\xi \! \pm \! 
\mi (\varepsilon b)^{1/6}}}} \right) \right) \right\vert_{\widetilde{\mu}^{
\mathrm{TP}}}^{\widetilde{\mu}}.
\end{align}
\end{proposition:s}

\emph{Proof.} Presenting $l(\widetilde{\mu}) \! = \! l(\widetilde{\mu},\tau) 
\! = \! l(\widetilde{\mu},+\infty)(1 \! + \! \Delta)^{1/2}$, where $\Delta \! 
:= \! \tfrac{l^{2}(\widetilde{\mu},\tau)-l^{2}(\widetilde{\mu},+\infty)}{l^{2}
(\widetilde{\mu},+\infty)}$, approximating $(1 \! + \! \Delta)^{1/2}$ as $1 \! 
+ \! \tfrac{1}{2} \Delta \! + \! \mathcal{O}(\Delta^{2})$, integrating by 
parts, and using the following integrals \cite{a24}, $\int \! \! \sqrt{x^{2} 
\! + \! a_{o}^{2}} \, \md x \! = \! \tfrac{x \sqrt{\smash[b]{x^{2}+a_{o}^{2}}
}}{2} \! + \! \tfrac{a_{o}^{2}}{2} \ln (x \! + \! \sqrt{\smash[b]{x^{2} \! + 
\! a_{o}^{2}}})$, $\int \! \! \tfrac{\md x}{\sqrt{\smash[b]{a_{o}x^{2}+b_{o}x+
c_{o}}}} \! = \! \tfrac{1}{\sqrt{\smash[b]{a_{o}}}} \ln (2 \sqrt{\smash[b]{a_{
o}}} \, \sqrt{\smash[b]{a_{o}x^{2} \! + \! b_{o}x \! + \! c_{o}}} 
\linebreak[4] 
+ \! 2a_{o}x \! + \! b_{o})$, and $\int \! \! \tfrac{\md x}{x \sqrt{\smash[b]{
a_{o}x^{2}+b_{o}x+c_{o}}}} \! = \! -\tfrac{1}{\sqrt{\smash[b]{c_{o}}}} \ln \! 
\left(\tfrac{2 \sqrt{\smash[b]{c_{o}}}\, \sqrt{\smash[b]{a_{o}x^{2}+b_{o}x+c_{
o}}}+b_{o}x+2c_{o}}{x} \right)$, one arrives at the result stated in 
Equation~(72). \hfill $\square$
\begin{ffff}
For $\widetilde{\mu}^{\mathrm{TP}} \! = \! \tfrac{(\varepsilon b)^{1/6}}{\sqrt{
\smash[b]{2}}} \! + \! \tau^{-1/3} \widetilde{\Lambda}$, where $\widetilde{
\Lambda} \! =_{\tau \to +\infty} \! \mathcal{O}(\tau^{\epsilon_{0}})$, $0 \! < 
\! \epsilon_{0} \! < \! \tfrac{1}{9}$,
\begin{align}
\int_{\widetilde{\mu}^{\mathrm{TP}}}^{\widetilde{\mu}}l(\xi) \, \md \xi 
\underset{\tau \to +\infty}{=}& \, \left(\widetilde{\mu} \! + \! \dfrac{
(\varepsilon b)^{1/3}}{\widetilde{\mu}} \right) \! \sqrt{\smash[b]{\widetilde{
\mu}^{2} \! + \! (\varepsilon b)^{1/3}}}+ \! \tau^{-2/3} \! \left(\! \left(a 
\! - \! \dfrac{\mi}{2} \right) \! \ln \! \left(2 \sqrt{\smash[b]{\widetilde{
\mu}^{2} \! + \! (\varepsilon b)^{1/3}}}+ \! 2 \widetilde{\mu} \right) \right. 
\nonumber \\
+&\left. \dfrac{1}{2 \sqrt{\smash[b]{3}}} \! \left(\! \left(a \! - \! \dfrac{
\mi}{2} \right) \! + \! \dfrac{2h_{0}(\tau)}{(\varepsilon b)^{1/3}} \right) \! 
\ln \! \left(\! \tfrac{(-\widetilde{\mu}+\sqrt{\smash[b]{2}} \, (\varepsilon 
b)^{1/6}+\sqrt{\smash[b]{3}} \, \sqrt{\smash[b]{\widetilde{\mu}^{2}+
(\varepsilon b)^{1/3}}} \,)(\widetilde{\mu}-\frac{(\varepsilon b)^{1/6}}{\sqrt{
\smash[b]{2}}})}{(\widetilde{\mu}+\sqrt{\smash[b]{2}} \, (\varepsilon b)^{1/6}
+\sqrt{\smash[b]{3}} \, \sqrt{\smash[b]{\widetilde{\mu}^{2}+(\varepsilon b)^{
1/3}}} \,)(\widetilde{\mu}+\frac{(\varepsilon b)^{1/6}}{\sqrt{\smash[b]{2}}})} 
\right) \! \right) \nonumber \\
-& \, \dfrac{3 \sqrt{3} \, (\varepsilon b)^{1/3}}{2} \! - \! \left(a \! - \! 
\dfrac{\mi}{2} \right) \! \tau^{-2/3} \ln \! \left(\! \sqrt{2} \, (\sqrt{3} \! 
+ \! 1)(\varepsilon b)^{1/6} \right) \! - \! 2 \sqrt{3} \, \tau^{-2/3} 
\widetilde{\Lambda}^{2} \nonumber \\
+& \, \dfrac{\tau^{-2/3}}{2 \sqrt{3}} \! \left(\! \left(a \! - \! \dfrac{\mi}{
2} \right) \! + \! \dfrac{2h_{0}(\tau)}{(\varepsilon b)^{1/3}} \right) \! 
\left(\ln \! \left(\tfrac{3(\varepsilon b)^{1/6}}{\sqrt{2}} \right) \! - \! 
\ln \widetilde{\Lambda} \! + \! \dfrac{1}{3} \ln \tau \right) \! + \! o(\tau^{
-2/3}).
\end{align}
\end{ffff}

\emph{Proof.} Substituting into Equation~(72) the expression for $\widetilde{
\mu}^{\mathrm{TP}}$ given in the Corollary and expanding the result into a 
power series in $\widetilde{\Lambda}$, one gets for the error estimate 
$\mathcal{O}(\tau^{-1} \widetilde{\Lambda}^{3}) \! + \! \mathcal{O}(\tau^{-
2/3} \widetilde{\Lambda}^{-2}) \! + \! \mathcal{O}(\tau^{-1} \widetilde{
\Lambda})$. \hfill $\square$
\begin{ffff}
For the choice of $\widetilde{\mu}^{\mathrm{TP}}$ given in 
Corollary~{\rm 4.1.1},
\begin{align}
-\mi \tau^{2/3} \int_{\widetilde{\mu}^{\mathrm{TP}}}^{\widetilde{\mu}}l(\xi) 
\, \md \xi \underset{\underset{\underset{\arg \widetilde{\mu}=0}{\widetilde{
\mu} \to \infty}}{\tau \to +\infty}}{=}& \, -\mi \! \left(\tau^{2/3} 
\widetilde{\mu}^{2} \! + \! \left(a \! - \! \dfrac{\mi}{2} \right) \! \ln 
\widetilde{\mu} \right) \! + \! \dfrac{\mi \tau^{2/3}(\varepsilon b)^{1/3}
(3 \sqrt{3} \! - \! 2)}{2} \! + \! 2 \sqrt{3} \, \mi \widetilde{\Lambda}^{2} 
\nonumber \\
+& \, \dfrac{\mi}{2 \sqrt{3}} \! \left(\! \left(a \! - \! \dfrac{\mi}{2} 
\right) \! + \! \dfrac{2h_{0}(\tau)}{(\varepsilon b)^{1/3}} \right) \! \left(
\ln \widetilde{\Lambda} \! - \! \dfrac{1}{3} \ln \tau \right) \! + \! C_{
\infty}^{\mathrm{WKB}} \! + \! o(1),
\end{align}
where
\begin{equation}
C_{\infty}^{\mathrm{WKB}} \! := \! \mi \! \left(a \! - \! \dfrac{\mi}{2} 
\right) \! \ln \! \left(\dfrac{(\sqrt{\smash[b]{3}}+ \! 1)(\varepsilon b)^{
1/6}}{2 \sqrt{\smash[b]{2}}} \right) \! - \! \dfrac{\mi}{2 \sqrt{\smash[b]{
3}}} \! \left(\! \left(a \! - \! \dfrac{\mi}{2} \right) \! + \! \dfrac{2h_{
0}(\tau)}{(\varepsilon b)^{1/3}} \right) \! \ln \! \left(\dfrac{3 \sqrt{
\smash[b]{2}} \, (\varepsilon b)^{1/6}}{(\sqrt{\smash[b]{3}}+ \! 1)^{2}} 
\right),
\end{equation}
and
\begin{align}
-\mi \tau^{2/3} \int_{\widetilde{\mu}^{\mathrm{TP}}}^{\widetilde{\mu}}l(\xi) 
\, \md \xi \underset{\underset{\widetilde{\mu} \to+0}{\tau \to +\infty}}{=}& 
\, \dfrac{\mi \tau^{2/3} \sqrt{\smash[b]{\varepsilon b}}}{\widetilde{\mu}} \! 
- \! \dfrac{3 \sqrt{\smash[b]{3}} \, \mi (\varepsilon b)^{1/3} \tau^{2/3}}{2} 
\! - \! 2 \sqrt{3} \, \mi \widetilde{\Lambda}^{2} \! - \! \dfrac{\mi}{2 \sqrt{
3}} \! \left(\! \left(a \! - \! \dfrac{\mi}{2} \right) \right. \nonumber \\
+&\left. \dfrac{2h_{0}(\tau)}{(\varepsilon b)^{1/3}} \right) \! \left(\ln 
\widetilde{\Lambda} \! - \! \dfrac{1}{3} \ln \tau \right) \! + \! C_{0}^{
\mathrm{WKB}} \! + \! o(1),
\end{align}
where
\begin{equation}
C_{0}^{\mathrm{WKB}} \! := \! \mi \! \left(a \! - \! \dfrac{\mi}{2} \right) \! 
\ln \! \left(\! \dfrac{\sqrt{\smash[b]{2}}}{\sqrt{\smash[b]{3}}+ \! 1} \right) 
\! + \! \dfrac{\mi}{2 \sqrt{\smash[b]{3}}} \! \left(\! \left(a \! - \! \dfrac{
\mi}{2} \right) \! + \! \dfrac{2h_{0}(\tau)}{(\varepsilon b)^{1/3}} \right) \! 
\ln \! \left(\dfrac{3(\varepsilon b)^{1/6}}{\sqrt{\smash[b]{2}}} \right).
\end{equation}
\end{ffff}

\emph{Proof.} Equation~(74) is the large-$\widetilde{\mu}$ asymptotics of 
Equation~(73). Equation~(76) is the small-$\widetilde{\mu}$ asymptotics of 
Equation~(73), but the branch of $l(\xi)$ in Equation~(76) is chosen 
differently than that in Equation~(73). \hfill $\square$
\begin{ffff}
For $\widetilde{\mu}^{\mathrm{TP}} \! = \! \mi (\varepsilon b)^{1/6} \! + 
\! \tau^{-4/9} \widehat{\Lambda}$, where $\widehat{\Lambda} \! =_{\tau \to 
+\infty} \! \mathcal{O}(\tau^{\epsilon_{0}})$, $0 \! < \! \epsilon_{0} \! < 
\! \tfrac{4}{9}$,
\begin{align}
\int_{\widetilde{\mu}^{\mathrm{TP}}}^{\widetilde{\mu}}l(\xi) \, \md \xi 
\underset{\tau \to +\infty}{=}& \, \left(\widetilde{\mu} \! + \! \dfrac{
(\varepsilon b)^{1/3}}{\widetilde{\mu}} \right) \! \sqrt{\smash[b]{\widetilde{
\mu}^{2} \! + \! (\varepsilon b)^{1/3}}}+ \! \tau^{-2/3} \! \left(\! \left(a 
\! - \! \dfrac{\mi}{2} \right) \! \ln \! \left(2 \sqrt{\smash[b]{\widetilde{
\mu}^{2} \! + \! (\varepsilon b)^{1/3}}}+ \! 2 \widetilde{\mu} \right) \right. 
\nonumber \\
+&\left. \dfrac{1}{2 \sqrt{\smash[b]{3}}} \! \left(\! \left(a \! - \! \dfrac{
\mi}{2} \right) \! + \! \dfrac{2h_{0}(\tau)}{(\varepsilon b)^{1/3}} \right) \! 
\ln \! \left(\! \tfrac{(-\widetilde{\mu}+\sqrt{\smash[b]{2}} \, (\varepsilon 
b)^{1/6}+\sqrt{\smash[b]{3}} \, \sqrt{\smash[b]{\widetilde{\mu}^{2}+
(\varepsilon b)^{1/3}}} \,)(\widetilde{\mu}-\frac{(\varepsilon b)^{1/6}}{\sqrt{
\smash[b]{2}}})}{(\widetilde{\mu}+\sqrt{\smash[b]{2}} \, (\varepsilon b)^{1/6}
+\sqrt{\smash[b]{3}} \, \sqrt{\smash[b]{\widetilde{\mu}^{2}+(\varepsilon b)^{
1/3}}} \,)(\widetilde{\mu}+\frac{(\varepsilon b)^{1/6}}{\sqrt{\smash[b]{2}}})} 
\right) \! \right) \nonumber \\
+& \, 2 \sqrt{2} \, \mi^{-3/2}(\varepsilon b)^{1/12} \tau^{-2/3} \widehat{
\Lambda}^{3/2} \! - \! \tau^{-2/3} \! \left(a \! - \! \tfrac{\mi}{2} \right) 
\ln \! \left(2 \mi (\varepsilon b)^{1/6} \right) \! + \! o(\tau^{-2/3}).
\end{align}
\end{ffff}

\emph{Proof.} Analogous to the proof of Corollary~4.1.1; but in this case, for 
the error estimate, one obtains $\mathcal{O}(\tau^{-10/9} \widehat{\Lambda}^{
-1/2}) \! + \! \mathcal{O}(\tau^{-8/9} \widehat{\Lambda}^{1/2})$. \hfill 
$\square$
\begin{ffff}
For the choice of $\widetilde{\mu}^{\mathrm{TP}}$ given in 
Corollary~{\rm 4.1.3},
\begin{align}
-\mi \tau^{2/3} \int_{\widetilde{\mu}^{\mathrm{TP}}}^{\widetilde{\mu}}l(\xi) 
\, \md \xi \underset{\underset{\underset{\arg \widetilde{\mu}=\frac{\pi}{2}}{
\widetilde{\mu} \to \infty}}{\tau \to +\infty}}{=}& \, -\mi \! \left(\tau^{
2/3} \widetilde{\mu}^{2} \! + \! (a \! - \! \tfrac{\mi}{2}) \ln \widetilde{
\mu} \right) \! - \! \mi \tau^{2/3}(\varepsilon b)^{1/3} \! - \! 2 \sqrt{2} 
\, \mi^{-1/2}(\varepsilon b)^{1/12} \widehat{\Lambda}^{3/2} \nonumber \\
+& \, \mi \! \left(a \! - \! \tfrac{\mi}{2} \right) \! \left(\ln \! \left(
\tfrac{(\varepsilon b)^{1/6}}{2} \right) \! + \! \tfrac{\pi \mi}{2} \right) \! 
+ \! \left(\tfrac{\mi}{2 \sqrt{\smash[b]{3}}} \! \left(a \! - \! \tfrac{\mi}{
2} \right) \! + \! \tfrac{2h_{0}(\tau)}{(\varepsilon b)^{1/3}} \right) \! \ln 
\! \left(\tfrac{\sqrt{\smash[b]{3}}+1}{\sqrt{\smash[b]{3}}-1} \right) \! + \! 
o(1).
\end{align}
\end{ffff}

\emph{Proof.} Equation~(79) is the large-$\widetilde{\mu}$ asymptotics of 
Equation~(78). \hfill $\square$
\begin{proposition:s}
\begin{gather}
T(\widetilde{\mu}) \! \underset{\underset{\arg \widetilde{\mu}=0}{\widetilde{
\mu} \to \infty}}{=} \! (b(\tau))^{-\frac{1}{2} \sigma_{3}} \! \left(\mathrm{
I} \! + \! \mathcal{O} \! \left(\widetilde{\mu}^{-1} \right) \right) \! (b
(\tau))^{\frac{1}{2} \sigma_{3}}, \\
T(\widetilde{\mu}) \! \underset{\widetilde{\mu} \to +0}{=} \! (b(\tau))^{-
\frac{1}{2} \sigma_{3}} \! \left(\frac{1}{\sqrt{2}} \! 
\begin{pmatrix}
1 & (\varepsilon b)^{1/2} \\
-(\varepsilon b)^{-1/2} & 1
\end{pmatrix} \! + \! \mathcal{O}(\widetilde{\mu}) \right) \! (b(\tau))^{
\frac{1}{2} \sigma_{3}}.
\end{gather}
\end{proposition:s}

\emph{Proof.} Asymptotics~(80) and~(81) are obtained {}from Equation~(70) upon 
choosing the same branches of $l(\widetilde{\mu})$ as in Equations~(74) 
and~(76), respectively. \hfill $\square$
\begin{eeee}
The elements of the matrices $\mathcal{O}(\widetilde{\mu}^{-1})$ and $\mathcal{
O}(\widetilde{\mu})$ in Equations~(80) and~(81), respectively, as functions of 
$\tau$, behave like $\mathcal{O}(1)$ as $\tau \! \to \! +\infty$. \hfill 
$\blacksquare$
\end{eeee}
\begin{proposition:s}
For $0 \! < \! \delta_{\diamondsuit}\! < \! \epsilon_{0} \! < \! \tfrac{1}{
9}$, set $\widetilde{\mu} \! = \! \tfrac{(\varepsilon b)^{1/6}}{\sqrt{2}} \! 
+ \! \tau^{-1/3} \widetilde{\Lambda}$, where $\widetilde{\Lambda} =_{\tau \to 
+\infty} \! \mathcal{O}(\tau^{\epsilon_{0}})$, and $\arg \widetilde{\Lambda} 
\! = \! \tfrac{1}{2}(1 \! - \! \widetilde{\varepsilon}_{1}) \pi$, $\widetilde{
\varepsilon}_{1} \! = \! \pm 1$. Then,
\begin{equation*}
T(\widetilde{\mu}) \! \underset{\tau \to +\infty}{=} \! \dfrac{(b(\tau))^{-
\frac{1}{2} \sigma_{3}}}{\sqrt{\smash[b]{2 \sqrt{3} \, (\sqrt{3}+\widetilde{
\varepsilon}_{1})}}} \! \left(\! 
\begin{pmatrix}
\sqrt{3}+\widetilde{\varepsilon}_{1} & -\widetilde{\varepsilon}_{1} \sqrt{2 
\varepsilon b} \\
\frac{\widetilde{\varepsilon}_{1} \sqrt{2}}{\sqrt{\varepsilon b}} & \sqrt{3}+
\widetilde{\varepsilon}_{1}
\end{pmatrix} \! + \! \mathcal{O} \! \left(\tau^{-(\epsilon_{0}-\delta_{
\diamondsuit})} \right) \right) \! (b(\tau))^{\frac{1}{2} \sigma_{3}}.
\end{equation*}
\end{proposition:s}

\emph{Proof.} The result follows {}from the definition of $T(\widetilde{\mu})$ 
given in Equation~(70) and conditions~(61) and~(63). \hfill $\square$
\begin{proposition:s}
For the conditions stated in Proposition~{\rm 4.1.4},
\begin{equation}
\int_{\widetilde{\mu}^{\mathrm{TP}}}^{\widetilde{\mu}} \operatorname{diag} \! 
\left((T(\xi))^{-1} \partial_{\xi}T(\xi) \right) \md \xi \underset{\tau \to 
+\infty}{=} o(1) \sigma_{3}.
\end{equation}
\end{proposition:s}

\emph{Proof.} Using the definition of $T(\widetilde{\mu})$ given in 
Equation~(70), one shows that the integrand is proportional to the 
$\xi$-independent factor $\mathscr{D}(\tau) \! := \! 4 \sqrt{\smash[b]{-a(\tau)
b(\tau)}}+ \! \tfrac{r(\tau)(\varepsilon b)^{4/3}}{2 \sqrt{\smash[b]{-a(\tau)
b(\tau)}}}$: this factor, as $\tau \! \to \! +\infty$, is $o(1)$ due to 
conditions~(61) and~(63). \hfill $\square$
\subsection{Asymptotics Near the Turning Points}
For the calculation of the monodromy data, one needs a more accurate 
approximation for the solution of Equation~(57) in the neighbourhoods of the 
turning points than that given by the WKB formula (cf. Equation~(68)). For 
the pair of coalescing turning points (resp., single turning points), this 
approximation is known to be given in terms of parabolic cylinder (resp., 
Airy) functions (see, for example, \cite{W,F}): more accurate statements 
are given below.
\begin{proposition:s}
Set $\widetilde{\mu} \! = \! \tfrac{(\varepsilon b)^{1/6}}{\sqrt{2}}+\tau^{-
1/3} \widetilde{\Lambda}$, where $\arg \widetilde{\Lambda} \! \in \! \{0,\pi
\}$, and $\vert \widetilde{\Lambda} \vert \! <_{\tau \to +\infty} \! \mathcal{
O}(\tau^{\epsilon_{0}})$, $0 \! < \! \epsilon_{0} \! < \! \tfrac{1}{9}$. 
Then, under conditions~{\rm (60)--(64)}, for any fundamental solution of 
Equation~{\rm (57)} the following asymptotic representation is valid,
\begin{equation}
\widetilde{\Psi}(\widetilde{\mu}) \underset{\tau \to +\infty}{=} \mathcal{N}
(\tau) \! \left(\mathrm{I} \! + \! \mathcal{O} \! \left(\tau^{-\epsilon} 
\right) \right) \! \widetilde{\Psi}_{0}(\widetilde{\Lambda}),
\end{equation}
where
\begin{equation}
\mathcal{N}(\tau) \! := \! \dfrac{\mi \sqrt{\smash[b]{b(\tau)}}}{(6 
\varepsilon b)^{1/4}} \! 
\begin{pmatrix}
-\frac{(\sqrt{\smash[b]{3}} \, -1) \sqrt{\smash[b]{\varepsilon b}}}{\sqrt{
\smash[b]{2}} \, b(\tau)} & \frac{\sqrt{\smash[b]{2}} \, \sqrt{\smash[b]{
\varepsilon b}}}{(\sqrt{\smash[b]{3}} \, -1)b(\tau)} \\
1 & 1
\end{pmatrix},
\end{equation}
$\epsilon$ is some positive number, and $\widetilde{\Psi}_{0}(\widetilde{
\Lambda})$ is a fundamental solution of
\begin{equation}
\dfrac{\partial \widetilde{\Psi}_{0}(\widetilde{\Lambda})}{\partial \widetilde{
\Lambda}} \! = \! \left(4 \sqrt{3} \, \mi \widetilde{\Lambda} \sigma_{3} \! + 
\! 
\begin{pmatrix}
0 & \widetilde{p} \\
\widetilde{q} & 0
\end{pmatrix} \! \right) \! \widetilde{\Psi}_{0}(\widetilde{\Lambda}),
\end{equation}
with
\begin{equation}
\begin{gathered}
\widetilde{p} \! := \! \dfrac{\mi (\varepsilon b)^{1/6}(\sqrt{\smash[b]{3}}+ 
\! 1)(4u_{0}(\tau) \! + \! (\sqrt{\smash[b]{3}}+ \! 1)r_{0}(\tau))}{2 \sqrt{
\smash[b]{2}}}, \\
\widetilde{q} \! := \! -\dfrac{\mi (\varepsilon b)^{1/6}(\sqrt{\smash[b]{3}}- 
\! 1)(4u_{0}(\tau) \! - \! (\sqrt{\smash[b]{3}}- \! 1)r_{0}(\tau))}{2 \sqrt{
\smash[b]{2}}}.
\end{gathered}
\end{equation}
\end{proposition:s}

\emph{Proof.} Making the change of variables $(\widetilde{\Psi}$, $\widetilde{
\mu}) \! \to \! (\widetilde{\Psi}^{\natural}$, $\widetilde{\Lambda})$ in 
Equation~(57), one shows that
\begin{equation*}
\dfrac{\partial \widetilde{\Psi}^{\natural}(\widetilde{\Lambda})}{\partial 
\widetilde{\Lambda}} \! \underset{\tau \to +\infty}{=} \! \left(P_{1}^{
\natural} \widetilde{\Lambda} \! + \! P_{0}^{\natural} \! + \! (b(\tau))^{-
\frac{1}{2} \sigma_{3}} \mathcal{O} \! \left(\tau^{-1/3} \widetilde{\Lambda}^{
2} \right) \! (b(\tau))^{\frac{1}{2} \sigma_{3}} \right) \! \widetilde{\Psi}^{
\natural}(\widetilde{\Lambda}),
\end{equation*}
where
\begin{gather*}
P_{1}^{\natural} \! := \! 
\begin{pmatrix}
\mi (r(\tau)-2) & -\frac{4 \sqrt{\smash[b]{2}} \, \mi (\varepsilon b)^{1/2}}{
b(\tau)} \\
-\frac{4 \sqrt{\smash[b]{2}} \, \mi b(\tau)}{(\varepsilon b)^{1/2}} & -\mi 
(r(\tau)-2)
\end{pmatrix}, \\
P_{0}^{\natural} \! := \! 
\begin{pmatrix}
-\frac{\mi (\varepsilon b)^{1/6}(r(\tau)+2) \tau^{1/3}}{\sqrt{\smash[b]{2}}} & 
-\frac{4 \mi \sqrt{\smash[b]{-a(\tau)b(\tau)}} \, \tau^{1/3}}{b(\tau)}+\frac{2 
\mi (\varepsilon b)^{2/3} \tau^{1/3}}{b(\tau)} \\
\frac{2 \mi b(\tau) \tau^{1/3}}{(\varepsilon b)^{1/3}}+\frac{\mi b(\tau)r(\tau)
(\varepsilon b)^{1/3} \tau^{1/3}}{2 \sqrt{\smash[b]{-a(\tau)b(\tau)}}} & \frac{
\mi (\varepsilon b)^{1/6}(r(\tau)+2) \tau^{1/3}}{\sqrt{\smash[b]{2}}}
\end{pmatrix}.
\end{gather*}
Now, {}from conditions~(61) and~(63), one finds that
\begin{gather*}
P_{1}^{\natural} \! = \! -4 \mi \! 
\begin{pmatrix}
1 & \frac{\sqrt{\smash[b]{2}} \, (\varepsilon b)^{1/2}}{b(\tau)} \\
\frac{\sqrt{\smash[b]{2}} \, b(\tau)}{(\varepsilon b)^{1/2}} & -1
\end{pmatrix} \! + \! \mi r_{0}(\tau) \tau^{-1/3} \sigma_{3}, \\
P_{0}^{\natural} \! = \! 
\begin{pmatrix}
-\frac{\mi (\varepsilon b)^{1/6}r_{0}(\tau)}{\sqrt{\smash[b]{2}}} & -\frac{2 
\mi (\varepsilon b)^{2/3}u_{0}(\tau)}{b(\tau)} \\
\frac{2 \mi b(\tau)}{(\varepsilon b)^{1/3}} \! \left( u_{0}(\tau)+\frac{r_{0}
(\tau)}{2} \right) & \frac{\mi (\varepsilon b)^{1/6}r_{0}(\tau)}{\sqrt{
\smash[b]{2}}}
\end{pmatrix} \! + \! b(\tau) \mathscr{E} \sigma_{-},
\end{gather*}
where
\begin{equation*}
\mathscr{E} \! := \! \mathcal{O} \! \left(r_{0}(\tau)u_{0}^{2}(\tau) \tau^{-
\frac{2}{3}} \right) \! + \! \mathcal{O} \! \left(r_{0}(\tau)u_{0}(\tau) \tau^{
-\frac{1}{3}} \right) \! + \! \mathcal{O} \! \left(u_{0}^{2}(\tau) \tau^{-
\frac{1}{3}} \right) \! + \! \mathcal{O} \! \left(u_{0}(\tau) \tau^{-\frac{2}{
3}} \right) \! + \! \mathcal{O} \! \left(\tau^{-\frac{1}{3}} \right).
\end{equation*}
Making the gauge transformation
\begin{equation*}
\widetilde{\Psi}^{\natural}(\widetilde{\Lambda}) \! = \! \mathcal{N}(\tau) 
\widehat{\Psi}(\widetilde{\Lambda}),
\end{equation*}
where $\mathcal{N}(\tau)$ is given by Equation~(84), one arrives at
\begin{equation*}
\dfrac{\partial \widehat{\Psi}(\widetilde{\Lambda})}{\partial \widetilde{
\Lambda}} \! \underset{\tau \to +\infty}{=} \! \left(4 \sqrt{3} \, \mi 
\widetilde{\Lambda} \sigma_{3} \! + \! 
\begin{pmatrix}
0 & \widetilde{p} \\
\widetilde{q} & 0
\end{pmatrix} \! + \! \mi r_{0}(\tau) \tau^{-1/3} \widetilde{\Lambda} \sigma_{
3} \! + \! \mathscr{E} \sigma_{-} \! + \! \mathcal{O} \! \left(\tau^{-1/3} 
\widetilde{\Lambda}^{2} \right) \right) \! \widehat{\Psi}(\widetilde{\Lambda}).
\end{equation*}
{}From Equation~(67), it follows that $\widetilde{p} \, \widetilde{q} \! =_{
\tau \to +\infty} \! \mathcal{O}(1)$. Since one of either $\widetilde{p}$ or 
$\widetilde{q}$ exhibits the power-like growth $\mathcal{O}(\tau^{\delta_{
\diamondsuit}})$, the other decays like $\mathcal{O}(\tau^{-\delta_{
\diamondsuit}})$. Using the fact that (cf. Equation~(64)) $\delta_{
\diamondsuit} \! \in \! (0,\tfrac{1}{9})$, $\widetilde{\Lambda}$ is real, and 
$\vert \widetilde{\Lambda} \vert \! \leqslant_{\tau \to +\infty} \! \mathcal{O}
(\tau^{\epsilon_{0}})$, one proceeds analogously as in \cite{a9,a18} to prove 
that
\begin{equation*}
\widehat{\Psi}(\widetilde{\Lambda}) \! \underset{\tau \to +\infty}{=} \! \left(
\mathrm{I} \! + \! \mathcal{O} \! \left(\tau^{-\epsilon} \right) \right) \! 
\widetilde{\Psi}_{0}(\widetilde{\Lambda}),
\end{equation*}
where $\widetilde{\Psi}_{0}(\widetilde{\Lambda})$ is a fundamental solution of 
Equation~(85). \hfill $\square$
\begin{proposition:s}
A fundamental solution of Equation~{\rm (85)} is
\begin{equation*}
\widetilde{\Psi}_{0}(\widetilde{\Lambda}) \! = \! 
\begin{pmatrix}
D_{-1-\nu} \! \left(2 \sqrt{2 \sqrt{3}} \, \mi \me^{\frac{\mi \pi}{4}} 
\widetilde{\Lambda} \right) & D_{\nu} \! \left(2 \sqrt{2 \sqrt{3}} \, \me^{
\frac{\mi \pi}{4}} \widetilde{\Lambda} \right) \\
\eth_{\mathscr{D}}D_{-1-\nu} \! \left(2 \sqrt{2 \sqrt{3}} \, \mi \me^{\frac{
\mi \pi}{4}} \widetilde{\Lambda} \right) & \eth_{\mathscr{D}}D_{\nu} \! \left(
2 \sqrt{2 \sqrt{3}} \, \me^{\frac{\mi \pi}{4}} \widetilde{\Lambda} \right)
\end{pmatrix},
\end{equation*}
where $\eth_{\mathscr{D}} \! := \! (\widetilde{p})^{-1}(\partial_{\widetilde{
\Lambda}} \! - \! 4 \sqrt{3} \, \mi \widetilde{\Lambda})$, and $D_{\star}
(\boldsymbol{\cdot})$ is the parabolic cylinder function \rm{\cite{a24}}.
\end{proposition:s}

\emph{Proof.} Changing the independent variable according to the rule 
$\widetilde{\Lambda} \! = \! \alpha x$, where $\alpha \! := \! \me^{-\frac{\mi 
\pi}{4}} \! \left(2 \sqrt{2 \sqrt{3}} \right)^{-1}$, one shows that 
$\widetilde{\Psi}_{0}(\widetilde{\Lambda}) \! := \! \mathcal{D}(x)$ satisfies 
$\partial_{x} \mathcal{D}(x) \! = \! \left(\tfrac{x}{2} \sigma_{3} \! + \! 
\left(
\begin{smallmatrix}
0 & p^{\ast} \\
q^{\ast} & 0
\end{smallmatrix}
\right) \right) \! \mathcal{D}(x)$, with $(p^{\ast},q^{\ast}) \! := \! (\alpha 
\widetilde{p},\alpha \widetilde{q})$, whose fundamental solution is 
\cite{a9,a18} $\mathcal{D}(x) \! = \! 
\left(
\begin{smallmatrix}
D_{-1-\nu}(\mi x) & D_{\nu}(x) \\
\dot{D}_{-1-\nu}(\mi x) & \dot{D}_{\nu}(x)
\end{smallmatrix}
\right)$, where $\dot{D}_{\star}(z) \! := \! (p^{\ast})^{-1}(\partial_{z}D_{
\star}(z) \! - \! \tfrac{z}{2}D_{\star}(z))$, $D_{\star}(\cdot)$ is the 
parabolic cylinder function, and
\begin{equation}
\nu \! + \! 1 \! := \! -p^{\ast} q^{\ast} \! \underset{\tau \to +\infty}{=} \! 
-\dfrac{\mi}{2 \sqrt{3}} \! \left(\! \left(a \! - \! \dfrac{\mi}{2} \right) \! 
+ \! \dfrac{2h_{0}(\tau)}{(\varepsilon b)^{1/3}} \right) \! \left(1 \! + \! 
\mathcal{O} \! \left(\tau^{-\frac{1}{3}+3\delta_{\diamondsuit}} \right) 
\right) \! \underset{\tau \to +\infty}{=} \! \mathcal{O}(1). \qquad \quad 
\square
\end{equation}
\begin{eeee}
In the proof of Proposition~4.2.2, the matrix-valued function $\mathcal{D}(x)$ 
is introduced. In Subsection~4.3, the following asymptotics of $\mathcal{D}
(x)$ as $\vert x \vert \! \to \! \infty$ are required (see, for example, 
\cite{a19}):
\begin{equation*}
\mathcal{D}(x) \! \underset{\underset{\arg x=\frac{\pi \mi}{4}+k \frac{\pi 
\mi}{2}}{x \to \infty}}{=} \! \left(\mathrm{I} \! + \! \mathcal{O} \! \left(
x^{-1} \right) \right) \! \exp \! \left(\! \left(\dfrac{x^{2}}{4} \! - \! 
(\nu \! + \! 1) \ln x \right) \! \sigma_{3} \right) \! \mathcal{R}_{k}, 
\quad k \! = \! -1,0,1,2,
\end{equation*}
where
\begin{gather*}
\mathcal{R}_{-1} \! := \! 
\begin{pmatrix}
\me^{-\frac{\pi \mi}{2}(\nu + 1)} & 0 \\
0 & -(p^{\ast})^{-1}
\end{pmatrix}, \qquad \quad \mathcal{R}_{0} \! := \! 
\begin{pmatrix}
\me^{-\frac{\pi \mi}{2}(\nu +1)} & 0 \\
-\frac{\mi}{p^{\ast}} \frac{\sqrt{2 \pi}}{\Gamma (\nu +1)} \me^{-\frac{\pi 
\mi}{2}(\nu +1)} & -(p^{\ast})^{-1}
\end{pmatrix}, \\
\mathcal{R}_{1} \! := \! 
\begin{pmatrix}
\me^{\frac{3 \pi \mi}{2}(\nu +1)} & \frac{\sqrt{2 \pi}}{\Gamma (-\nu)} \me^{
\pi \mi (\nu +1)} \\
-\frac{\mi}{p^{\ast}} \frac{\sqrt{2 \pi}}{\Gamma (\nu +1)} \me^{-\frac{\pi 
\mi}{2}(\nu +1)} & -(p^{\ast})^{-1}
\end{pmatrix}, \qquad \mathcal{R}_{2} \! := \! 
\begin{pmatrix}
\me^{\frac{3 \pi \mi}{2}(\nu +1)} & \frac{\sqrt{2 \pi}}{\Gamma (-\nu)} \me^{
\pi \mi (\nu +1)} \\
0 &  -\frac{1}{p^{\ast}} \me^{-2 \pi \mi (\nu +1)}
\end{pmatrix},
\end{gather*}
and $\Gamma (\boldsymbol{\cdot})$ is the gamma function \cite{a24}. \hfill 
$\blacksquare$
\end{eeee}
\subsection{Asymptotic Matching as $\tau \! \to \! +\infty$}
Let $\widetilde{\Psi}(\widetilde{\mu})$ be a fundamental solution of 
Equation~(57) corresponding to (see Equation~(83)) the function $\widetilde{
\Psi}_{0}(\widetilde{\Lambda})$ given in Proposition~4.2.2. Define $\widetilde{
Y}_{0}^{\infty}(\widetilde{\mu}) \! := \! \tau^{-(1/12) \sigma_{3}}Y_{0}^{
\infty}(\widetilde{\mu} \tau^{-1/6})$, where $Y_{0}^{\infty}(\cdot)$ is the 
canonical solution of Equation~(55) (cf. Proposition~2.2), and
\begin{equation*}
L_{\infty} \! := \! \left(\widetilde{\Psi}(\widetilde{\mu}) \right)^{-1} \! 
\widetilde{Y}_{0}^{\infty}(\widetilde{\mu}).
\end{equation*}
\begin{cccc}
Let conditions~{\rm (60)--(64)} be valid, and
\begin{equation}
b(\tau) \tau^{-\frac{1}{3} \Im (a)} \! \underset{\tau \to +\infty}{=} \! 
\mathcal{O}(1).
\end{equation}
Then there exists $\epsilon_{\nu}^{\infty} \! > \! 0$ and $\delta_{\infty} \! 
> \! 0$ such that for
\begin{equation}
\vert \Re (\nu \! + \! 1) \vert \! \underset{\tau \to +\infty}{<} \! 
\epsilon_{\nu}^{\infty},
\end{equation}
\begin{align}
L_{\infty} \underset{\tau \to +\infty}{=}& \, -\mi (\mathcal{R}_{0})^{-1} \exp 
\! \left(\! \left(\! \dfrac{(3 \sqrt{\smash[b]{3}}- \! 2) \mi (\varepsilon b)^{
1/3} \tau^{2/3}}{2} \! + \! \left(-\dfrac{\mi a}{6} \! + \! \dfrac{(\nu \! + 
\! 1)}{3} \right) \! \ln \tau \! + \! \dfrac{\mi \pi (\nu \! + \! 1)}{4} 
\right. \right. \nonumber \\
+&\left. \left. \, \widetilde{C}_{\infty}^{\mathrm{WKB}} \right) \! \sigma_{3} 
\right) \! \left(\! \dfrac{(\varepsilon b)^{1/4} \sqrt{\smash[b]{\sqrt{
\smash[b]{3}}+ \! 1}}}{2^{1/4} \sqrt{\smash[b]{b(\tau)}}} \right)^{\sigma_{3}} 
\! \sigma_{1} \! \left({\rm I} \! + \! \mathcal{O} \! \left(\tau^{-\delta_{
\infty}} \right) \right),
\end{align}
where $\mathcal{R}_{0}$ is defined in Remark~{\rm 4.2.1}, and
\begin{align}
\widetilde{C}_{\infty}^{\mathrm{WKB}} :=& \, \left(-\dfrac{3 \mi}{2} \! \left(
a \! - \! \dfrac{\mi}{2} \right) \! + \! 2(\nu \! + \! 1) \right) \! \ln 2 \! 
+ \! \dfrac{5}{4}(\nu \! + \! 1) \ln 3 \! + \! \dfrac{1}{6} \! \left(\mi \! 
\left(a \! - \! \dfrac{\mi}{2} \right) \! + \! (\nu \! + \! 1) \right) \! \ln 
(\varepsilon b) \nonumber \\
+& \, \left(\mi \! \left(a \! - \! \dfrac{\mi}{2} \right) \! - \! 2(\nu \! + 
\! 1) \right) \! \ln \! \left(\sqrt{\smash[b]{3}} + \! 1 \right).
\end{align}
\end{cccc}

\emph{Proof.} Denote by $\widetilde{\Psi}_{\mathrm{WKB}}(\widetilde{\mu})$ 
the solution of Equation~(57) which has WKB asymptotics (68) in the 
canonical domain containing the Stokes curve approaching the positive 
$\widetilde{\mu}$-real axis {}from above as $\widetilde{\mu} \! \to \! 
\infty$. Rewriting $L_{\infty}$ as $L_{\infty} \! = \! \left((\widetilde{\Psi}
(\widetilde{\mu}))^{-1} \widetilde{\Psi}_{\mathrm{WKB}}(\widetilde{\mu}) 
\right) \! \left((\widetilde{\Psi}_{\mathrm{WKB}}(\widetilde{\mu}))^{-1} 
\widetilde{Y}_{0}^{\infty}(\widetilde{\mu}) \right)$, and noting that the 
members of the pairs in parentheses are solutions of the same equation 
(Equation~(57)), they differ by right-hand, $\widetilde{\mu}$-independent 
matrices: the latter matrices are calculated by taking different $\widetilde{
\mu}$ limits $(\widetilde{\mu} \! \to \! +\infty$ and $\widetilde{\mu} \! \to 
\! (\varepsilon b)^{1/6}/\sqrt{2})$. More precisely, modulo factors that are 
$(\mathrm{I} \! + \! \mathcal{O}(\tau^{-\widetilde{\delta}}))$, for some 
$\widetilde{\delta} \! > \! 0$,
\begin{equation*}
L_{\infty} \! \underset{\tau \to +\infty}{=} \! (\, \underbrace{(\mathcal{N}
(\tau) \widetilde{\Psi}_{0}(\widetilde{\Lambda}))^{-1}T(\widetilde{\mu})}_{
\underset{\widetilde{\mu}=\frac{(\varepsilon b)^{1/6}}{\sqrt{2}}+\tau^{-1/3} 
\widetilde{\Lambda}}{\widetilde{\Lambda} \sim \tau^{\epsilon_{o}}, \quad \arg 
\widetilde{\Lambda}=0}} \,)(\, \underbrace{(\widetilde{\Psi}_{\mathrm{WKB}}
(\widetilde{\mu}))^{-1} \widetilde{Y}_{0}^{\infty}(\widetilde{\mu})}_{
\widetilde{\mu} \to \infty, \quad \arg \widetilde{\mu}=0} \,),
\end{equation*}
where $\mathcal{N}(\tau)$ and $\widetilde{\Psi}_{0}(\widetilde{\Lambda})$, 
respectively, are given in Propositions~4.2.1 and~4.2.2 (see, also, 
Remark~4.2.1), $T(\widetilde{\mu})$ is given in Proposition~4.1.4, $\widetilde{
\Psi}_{\mathrm{WKB}}(\widetilde{\mu})$ is defined by Equation~(68) (with $T
(\widetilde{\mu})$ given in Equation~(80) and the phase-part of the WKB 
formula estimated in Equations~(74), (75), and~(82)), and $\widetilde{Y}_{0}^{
\infty}(\cdot)$ is defined in the beginning of this subsection. The additional 
conditions~(88) and~(89) are necessary in order to guarantee the power-like 
decay, $\mathcal{O}(\tau^{-\delta_{\infty}})$, $\delta_{\infty} \! > \! 0$, 
of the off-diagonal entries of $L_{\infty}$ given in Equation~(90). \hfill 
$\square$

Define $\widetilde{X}_{0}^{0}(\widetilde{\mu}) \! := \! \tau^{-(1/12) \sigma_{
3}}X_{0}^{0}(\widetilde{\mu} \tau^{-1/6})$, where $X_{0}^{0}(\cdot)$ is the 
canonical solution of Equation~(55) (cf. Proposition~2.2), and
\begin{equation*}
L_{0} \! := \! \left(\widetilde{\Psi}(\widetilde{\mu}) \right)^{-1} \! 
\widetilde{X}_{0}^{0}(\widetilde{\mu}).
\end{equation*}
\begin{cccc}
Under conditions~{\rm (60)--(64)} and~{\rm (89)}, there exists $\epsilon_{
\nu}^{0} \! > \! 0$ and $\delta_{0} \! > \! 0$ such that for
\begin{equation}
\vert \Re (\nu \! + \! 1) \vert \! \underset{\tau \to +\infty}{<} \!
\epsilon_{\nu}^{0},
\end{equation}
\begin{align}
L_{0} \underset{\tau \to +\infty}{=}& \, -\mi (\mathcal{R}_{2})^{-1} \exp \! 
\left(\! \left(\! \dfrac{3 \sqrt{\smash[b]{3}} \, \mi (\varepsilon b)^{1/3} 
\tau^{2/3}}{2} \! + \! \dfrac{(\nu \! + \! 1)}{3} \ln \tau \! + \! \dfrac{\mi 
\pi (\nu \! + \! 1)}{4} \! + \! \widetilde{C}_{0}^{\mathrm{WKB}} \right) \! 
\sigma_{3} \right) \nonumber \\
\times& \, \left(\! \dfrac{\sqrt{\smash[b]{\sqrt{\smash[b]{3}}+ \! 1}}}{2^{
1/4}} \right)^{\sigma_{3}} \! \sigma_{2} \! \left(\mathrm{I} \! + \! \mathcal{
O} \! \left(\tau^{-\delta_{0}} \right) \right),
\end{align}
where $\mathcal{R}_{2}$ is defined in Remark~{\rm 4.2.1}, and
\begin{equation}
\widetilde{C}_{0}^{\mathrm{WKB}} \! := \! \left(\! -\dfrac{\mi}{2} \! \left(a 
\! - \! \dfrac{\mi}{2} \right) \! + \! (\nu \! + \! 1) \right) \! \ln 2 \! + 
\! \frac{5}{4}(\nu \! + \! 1) \ln 3 \! + \! \dfrac{1}{6}(\nu \! + \! 1) \ln 
(\varepsilon b) \! + \! \mi \! \left(a \! - \! \frac{\mi}{2} \right) \! \ln 
\! \left(\sqrt{\smash[b]{3}} \! + \! 1 \right).
\end{equation}
\end{cccc}

\emph{Proof.} Denote by $\widetilde{\Psi}_{\mathrm{WKB}}(\widetilde{\mu})$ 
the solution of Equation~(57) which has WKB asymptotics (68) in the 
canonical domain containing the Stokes curve approaching the positive 
$\widetilde{\mu}$-real axis {}from above as $\widetilde{\mu} \! \to \! 0$. 
Rewriting $L_{0}$ as $L_{0} \! = \! \left((\widetilde{\Psi}(\widetilde{\mu})
)^{-1} \widetilde{\Psi}_{\mathrm{WKB}}(\widetilde{\mu}) \right) \! \left(
(\widetilde{\Psi}_{\mathrm{WKB}}(\widetilde{\mu}))^{-1} \widetilde{X}_{0}^{0}
(\widetilde{\mu}) \right)$, and arguing analogously as in the proof of 
Lemma~4.3.1, one can estimate $L_{0}$ by taking different $\widetilde{\mu}$ 
limits $(\widetilde{\mu} \! \to \! +0$ and $\widetilde{\mu} \! \to \! 
(\varepsilon b)^{1/6}/\sqrt{2})$. More precisely, modulo factors that are 
$(\mathrm{I} \! + \! \mathcal{O}(\tau^{-\widehat{\delta}}))$, for some 
$\widehat{\delta} \! > \! 0$,
\begin{equation*}
L_{0} \! \underset{\tau \to +\infty}{=} \! (\, \underbrace{(\mathcal{N}(\tau) 
\widetilde{\Psi}_{0}(\widetilde{\Lambda}))^{-1}T(\widetilde{\mu})}_{\underset{
\widetilde{\mu}=\frac{(\varepsilon b)^{1/6}}{\sqrt{2}}+\tau^{-1/3} \widetilde{
\Lambda}}{\widetilde{\Lambda} \sim \tau^{\epsilon_{o}}, \quad \arg \widetilde{
\Lambda}=\pi}} \,)(\, \underbrace{(\widetilde{\Psi}_{\mathrm{WKB}}(\widetilde{
\mu}))^{-1} \widetilde{X}_{0}^{0}(\widetilde{\mu})}_{\widetilde{\mu} \to 0, 
\quad \arg \widetilde{\mu}=0} \,),
\end{equation*}
where $\mathcal{N}(\tau)$ and $\widetilde{\Psi}_{0}(\widetilde{\Lambda})$, 
respectively, are given in Propositions~4.2.1 and~4.2.2 (see, also, 
Remark~4.2.1), $T(\widetilde{\mu})$ is given in Proposition~4.1.4, 
$\widetilde{\Psi}_{\mathrm{WKB}}(\widetilde{\mu})$ is defined by Equation~(68) 
(with $T(\widetilde{\mu})$ given in Equation~(81) and the phase-part of the 
WKB formula estimated in Equations~(76), (77), and~(82)), and $\widetilde{X}_{
0}^{0}(\cdot)$ is defined in the paragraph preceding Lemma~4.3.2. \hfill 
$\square$
\begin{dddd}
Under conditions~{\rm (60)--(64)}, {\rm (88)}, {\rm (89)}, and~{\rm (92)}, the 
connection matrix $G$ (cf. Equation~{\rm (28))} has the following asymptotics,
\begin{align}
G \underset{\tau \to +\infty}{=}& 
\begin{pmatrix}
-\frac{\mi \me^{-2 \pi \mi (\nu +1)} \sqrt{\smash[b]{b(\tau)}} \, \me^{(z_{2}
-z_{1})}}{(\varepsilon b)^{1/4}} & \frac{\sqrt{\smash[b]{2 \pi}} \, 
(\varepsilon b)^{1/4} \sqrt{\smash[b]{2 + \sqrt{\smash[b]{3}}}} \, \me^{-2 \pi 
\mi (\nu +1)} \me^{(z_{2}+z_{1})}}{p^{\ast} \sqrt{\smash[b]{b(\tau)}} \, 
\Gamma (\nu +1)} \\
-\frac{\mi p^{\ast} \sqrt{\smash[b]{2 \pi}} \, \me^{\pi \mi (\nu+1)} \sqrt{
\smash[b]{b (\tau)}} \, \me^{-(z_{2}+z_{1})}}{(\varepsilon b)^{1/4} \sqrt{
\smash[b]{2 +\sqrt{\smash[b]{3}}}} \, \Gamma (-\nu)} & \frac{\mi (\varepsilon 
b)^{1/4} \me^{-(z_{2}-z_{1})}}{\sqrt{\smash[b]{b(\tau)}}}
\end{pmatrix} \nonumber \\
\times& \, \left(\mathrm{I} \! + \! \mathcal{O} \! \left(\tau^{-\delta} 
\right) \right),
\end{align}
where
\begin{align}
z_{2} \! - \! z_{1} :=& \, \mi (\varepsilon b)^{1/3} \tau^{2/3} \! + \! \dfrac{
\mi a}{6} \ln \tau \! + \! \mi \! \left(a \! - \! \dfrac{\mi}{2} \right) \! 
\ln 2 \! - \! \dfrac{\mi}{6} \! \left(a \! - \! \dfrac{\mi}{2} \right) \! \ln 
(\varepsilon b) \! + \! (\nu \! + \! 1) \ln \! \left(\! \dfrac{\sqrt{\smash[b]{
3}}+ \! 1}{\sqrt{\smash[b]{3}}- \! 1} \right), \\
z_{2} \! + \! z_{1} :=& \, \! \left(3 \sqrt{\smash[b]{3}} \! - \! 1 \right) \! 
\mi (\varepsilon b)^{1/3} \tau^{2/3} \! + \! \left(\dfrac{2}{3}(\nu \! + \! 1) 
\! - \! \dfrac{\mi a}{6} \right) \! \ln \tau \! + \! \dfrac{\mi \pi}{2}(\nu \! 
+ \! 1) \! + \! \dfrac{5}{2}(\nu \! + \! 1) \ln 3 \nonumber \\
+& \, \left(-\mi \! \left(a \! - \! \dfrac{\mi}{2} \right) \! + \! 2(\nu \! + 
\! 1) \right) \! \ln 2 \! + \! \dfrac{1}{6} \! \left(\mi \! \left(a \! - \! 
\dfrac{\mi}{2} \right) \! + \! 2(\nu \! + \! 1) \right) \! \ln (\varepsilon b) 
\! + \! \left(\mi \! \left(a \! - \! \dfrac{\mi}{2} \right) \right. \nonumber 
\\
-&\left. \, (\nu \! + \! 1) \right) \ln \! \left(\! \dfrac{\sqrt{\smash[b]{3}}
+ \! 1}{\sqrt{\smash[b]{3}}- \! 1} \right),
\end{align}
and $\delta \! > \! 0;$ in particular, $G \! =_{\tau \to +\infty} \! \mathcal{
O}(1)$.
\end{dddd}

\emph{Proof.} {}From the definition of $G$ (Equation~(28)), $L_{\infty}$, 
and $L_{0}$, one arrives at $G \! = \! (L_{0})^{-1}L_{\infty}$: now, via 
straightforward calculations, the result stated in the Theorem is a 
consequence of Lemmae~4.3.1 and~4.3.2. \hfill $\square$
\begin{eeee}
It follows {}from Equation~(28) that $\det (G) \! = \! 1$: one shows, via the 
well-known identity $\Gamma (z) \Gamma (1 \! - \! z) \! = \! \tfrac{\pi}{\sin 
\pi z}$, that asymptotics~(95) is consistent with this fact. \hfill 
$\blacksquare$
\end{eeee}
\begin{ffff}
Let $g_{ij}$, $i,j \! = \! 1,2$, be the matrix elements of $G$. There exists 
$\epsilon_{\nu} \! > \! 0$ such that for
\begin{equation}
\vert \arg (g_{11}g_{22}) \vert \! < \! \epsilon_{\nu}, \qquad \quad g_{11}
g_{12}g_{21}g_{22} \! \not= \! 0,
\end{equation}
the functions $b(\tau)$, $r_{0}(\tau)$, $u_{0}(\tau)$, and $h_{0}(\tau)$ have 
the following asymptotic representation,
\begin{align}
b(\tau) \underset{\tau \to + \infty}{=}& \, -g_{11}^{2} \sqrt{\smash[b]{
\varepsilon b}} \, \me^{4 \pi \mi (\nu +1)} \exp \! \left(\! -2 \! \left(\mi 
(\varepsilon b)^{1/3} \tau^{2/3} \! + \! \dfrac{\mi a}{6} \ln \tau \! + \! 
\mi \! \left(\! a \! - \! \dfrac{\mi}{2} \right) \ln 2 \right. \right. 
\nonumber \\
-&\left. \left. \, \dfrac{\mi}{6} \! \left(\! a \! - \! \dfrac{\mi}{2} \right) 
\! \ln (\varepsilon b) \! + \! (\nu \! + \! 1) \ln \! \left(\! \dfrac{\sqrt{
\smash[b]{3}}+ \! 1}{\sqrt{\smash[b]{3}}- \! 1} \right) \! + \! o \! \left(
\tau^{-\delta} \right) \right) \right), \\
r_{0}(\tau) \underset{\tau \to +\infty}{=}& \, \dfrac{4 \sqrt{\smash[b]{2}} 
\, \sqrt{\smash[b]{\nu \! + \! 1}} \, \me^{\frac{3 \pi \mi}{4}}}{3^{1/4}
(\varepsilon b)^{1/6}} \sinh \! \left(\! \ln \widehat{x} \! - 2 \pi \mi (\nu 
\! + \! 1) \! + \! \dfrac{1}{2} \ln \! \left(2 \! + \! \sqrt{\smash[b]{3}} 
\right) \! + \! \mi \widehat{y} \! + \! z_{n} \! + \! o \! \left(\tau^{-
\delta} \right) \right), \\
u_{0}(\tau) \underset{\tau \to +\infty}{=}& \, \dfrac{2 \sqrt{\smash[b]{\nu \! 
+ \! 1}} \, \me^{\frac{3 \pi \mi}{4}}}{3^{1/4}(\varepsilon b)^{1/6}} \cosh \! 
\left(\ln \widehat{x} \! - 2 \pi \mi (\nu \! + \! 1) \! + \! \mi \widehat{y} 
\! + \! z_{n} \! + \! o \! \left(\tau^{-\delta} \right) \right), \\
h_{0}(\tau) \underset{\tau \to +\infty}{=}& \, -\dfrac{(\varepsilon b)^{1/3}}{
2} \! \left(\! \dfrac{\sqrt{\smash[b]{3}}}{\pi} \ln (g_{11}g_{22}) \! + \! 
\left(\! a \! - \! \dfrac{\mi}{2} \right) \! + \! o \! \left(\tau^{-\delta} 
\right) \right),
\end{align}
where
\begin{equation*}
\widehat{x} \! := \! \left\vert \dfrac{g_{22}g_{21} \Gamma(-(\nu \! + \! 1))}{
g_{11}g_{12} \Gamma(\nu \! + \! 1)} \right\vert^{1/2}, \qquad \quad \widehat{
y} \! := \! \dfrac{1}{2} \arg \! \left(\! \dfrac{g_{22}g_{21} \Gamma(-(\nu \! 
+ \! 1))}{g_{11}g_{12} \Gamma(\nu \! + \! 1)} \right),
\end{equation*}
\begin{align*}
z_{n} :=& \, 3 \sqrt{\smash[b]{3}} \, \mi (\varepsilon b)^{1/3} \tau^{2/3} 
\! + \! \dfrac{2}{3}(\nu \! + \! 1) \ln \tau \! + \! \dfrac{\pi \mi}{4} \! 
+ \! 2 (\nu \! + \! 1) \ln 2  \! + \! \dfrac{5}{2}(\nu \! + \! 1) \ln 3 \\
+& \, \dfrac{1}{3}(\nu \! + \! 1) \ln (\varepsilon b) \! + \! \mi a \ln \! 
\left(\! \dfrac{\sqrt{\smash[b]{3}}+ \! 1}{\sqrt{\smash[b]{3}}- \! 1} \right),
\end{align*}
and $\delta \! > \! 0$.
\end{ffff}

\emph{Proof.} Multiplying the diagonal elements of $G$ (Equation~(95)), one 
arrives at
\begin{equation}
(\nu \! + \! 1) \! = \! \dfrac{\mi}{2 \pi} \ln (g_{11}g_{22}) \! \left(1 \! + 
\! o \! \left(\tau^{-\delta} \right) \right):
\end{equation}
taking conditions~(89) and~(92) into account, and setting $\epsilon_{\nu} \! 
:= \! 2 \pi \min \{\epsilon_{\nu}^{\infty},\epsilon_{\nu}^{0}\}$, one obtains 
the first of conditions~(98). Equation~(99) is obtained {}from the $(1 \, 
1)$-element of $G$. {}From the definition of $p^{\ast}$ and $q^{\ast}$ given 
in the proof of Proposition~4.2.2 and Equations~(86), one obtains
\begin{equation*}
r_{0}(\tau) \! = \! \dfrac{4 \me^{-\frac{\mi \pi}{4}}}{3^{1/4}(\varepsilon 
b)^{1/6}} \! \left(\! \dfrac{p^{\ast}}{\sqrt{\smash[b]{3}}+ \! 1} \! + \! 
\dfrac{q^{\ast}}{\sqrt{\smash[b]{3}}- \! 1} \right), \qquad u_{0}(\tau) \! = 
\!\dfrac{\me^{-\frac{\mi \pi}{4}}}{3^{1/4}(\varepsilon b)^{1/6}} \! \left(\! 
\dfrac{\sqrt{\smash[b]{3}}- \! 1}{\sqrt{\smash[b]{3}}+ \! 1}p^{\ast} \! - \! 
\dfrac{\sqrt{\smash[b]{3}} + \! 1}{\sqrt{\smash[b]{3}}- \! 1}q^{\ast} \right).
\end{equation*}
Using the second of conditions~(98), $p^{\ast}q^{\ast} \! = \! -(\nu \! + \! 
1)$ (cf. Equation~(87)), and the identity $\Gamma (-\nu) \! = \! -(\nu \! + \! 
1) \Gamma (-(\nu \! + \! 1))$, one deduces, {}from the second, respectively, 
first, column of $G$, the following representation for $p^{\ast}$, 
respectively, $q^{\ast}$:
\begin{equation*}
p^{\ast} \! = \! -\mi \sqrt{\smash[b]{2 \pi (2 \! + \! \sqrt{\smash[b]{3}})}} 
\, \dfrac{g_{22} \me^{-2 \pi \mi(\nu +1)} \me^{2z_{2}}}{g_{12} \Gamma (\nu 
+1)}, \qquad \quad q^{\ast} \! = \! \dfrac{\sqrt{\smash[b]{2 \pi}} \, g_{11} 
\me^{3 \pi \mi (\nu +1)} \me^{-2z_{2}}}{\sqrt{\smash[b]{2 \! + \! \sqrt{
\smash[b]{3}}}} \, g_{21} \Gamma(-(\nu +1))}.
\end{equation*}
Substituting the latter formulae into the relations for $r_{0}(\tau)$ and $u_{
0}(\tau)$ given above, one obtains Equations~(100) and~(101). Equation~(102) 
is a direct consequence of Equations~(87) and (103). \hfill $\square$
\begin{proposition:s}
Let $G$ be the connection matrix of Equation~{\rm (55)} with 
$\tau$-independent elements satisfying conditions~{\rm (98)}. Then the 
corresponding isomonodromy deformations have the following asymptotic 
representation,
\begin{align}
\sqrt{\smash[b]{-a(\tau)b(\tau)}} \underset{\tau \to +\infty}{=}& \, \dfrac{
(\varepsilon b)^{2/3}}{2} \! + \! \dfrac{\sqrt{\smash[b]{\varepsilon b(\nu \! 
+ \! 1)}} \, \me^{\frac{3 \pi \mi}{4}}}{3^{1/4} \tau^{1/3}} \cosh \! \left(
\ln \widehat{x} \! - \! 2 \pi \mi (\nu \! + \! 1) \! + \! \mi \widehat{y} \! 
+ \! z_{n} \! + \! o \! \left(\tau^{-\delta} \right) \right), \\
a(\tau)d(\tau) \underset{\tau \to +\infty}{=}& -\dfrac{\mi \varepsilon b}{4} 
\! + \! \dfrac{\mi (\varepsilon b)^{5/6} \sqrt{\smash[b]{\nu \! + \! 1}} \, 
\me^{\frac{3 \pi \mi}{4}}}{2 \cdot 3^{1/4} \tau^{1/3}} \! \left(\! \sqrt{2} 
\sinh \! \left(\ln \widehat{x} \! - \! 2 \pi \mi (\nu \! + \! 1) \! + \! 
\dfrac{1}{2} \ln \! \left(2 \! + \! \sqrt{\smash[b]{3}} \right) \right. 
\right. \nonumber \\
+&\left. \left. \, \mi \widehat{y} \! + \! z_{n} \! + \! o \! \left(\tau^{-
\delta} \right) \right) \! - \! \cosh \! \left(\ln \widehat{x} \! - \! 2 \pi 
\mi (\nu \! + \! 1) \! + \! \mi \widehat{y} \! + \! z_{n} \! + \! o \! \left(
\tau^{-\delta} \right) \right) \right), \\
b(\tau)c(\tau) \underset{\tau \to +\infty}{=}& \, -\dfrac{\mi \varepsilon b}{
4} \! - \! \dfrac{\mi (\varepsilon b)^{5/6} \sqrt{\smash[b]{\nu \! + \! 1}} 
\, \me^{\frac{3 \pi \mi}{4}}}{2 \cdot 3^{1/4} \tau^{1/3}} \! \left(\! \sqrt{
\smash[b]{2}} \sinh \! \left(\ln \widehat{x} \! - \! 2 \pi \mi (\nu \! + \! 
1) \! + \! \dfrac{1}{2} \ln \! \left(2 \! + \! \sqrt{\smash[b]{3}} \right) 
\right. \right. \nonumber \\
+&\left. \left. \, \mi \widehat{y} \! + \! z_{n} \! + \! o \! \left(\tau^{-
\delta} \right) \right) \! - \! \cosh \! \left(\ln \widehat{x} \! - \! 2 \pi 
\mi (\nu \! + \! 1) \! + \! \mi \widehat{y} \! + \! z_{n} \! + \! o \! \left(
\tau^{-\delta} \right) \right) \right), \\
c(\tau)d(\tau) \underset{\tau \to +\infty}{=}& \, \dfrac{(\varepsilon b)^{2/3}
}{4} \! - \! \dfrac{\sqrt{\smash[b]{\varepsilon b(\nu \! + \! 1)}} \, \me^{
\frac{3 \pi \mi}{4}}}{3^{1/4} \tau^{1/3}} \cosh \! \left(\ln \widehat{x} \! - 
\! 2 \pi \mi (\nu \! + \! 1) \! + \! \mi \widehat{y} \! + \! z_{n} \! + \! o 
\! \left(\tau^{-\delta} \right) \right),
\end{align}
and $b(\tau)$ is given in Equation~{\rm (99)}.
\end{proposition:s}

\emph{Proof.} If the elements of $G$ are $\tau$-independent, then any 
functions whose asymptotics (as $\tau \! \to \! +\infty)$ are given by 
Equations~(99)--(102) satisfy conditions~(60)--(64), (88), (89), and~(92); 
therefore, one can now use the justification scheme suggested in \cite{a20} 
(see, also, \cite{a22}). Equations~(104)--(107) are obtained {}from 
Equations~(62), (63), (65), and~(66) by the direct substitution of 
Equations~(100)--(102). \hfill $\square$

According to Proposition~1.2, $u(\tau)$, the solution of the degenerate third 
Painlev\'{e} equation~(\ref{eq:dp3}), can be written in terms of the 
functions $a(\tau)$ and $b(\tau)$ as
\begin{equation*}
u(\tau) \! = \! \varepsilon \tau^{1/3} \sqrt{\smash[b]{-a(\tau)b(\tau)}}, 
\quad \varepsilon \! = \! \pm 1
\end{equation*}
(cf. Equations~(56)); therefore, the result (as $\tau \! \to \! +\infty)$ 
stated in Theorem~3.1 follows {}from Equation~(104).

Comparing Equation~(37) with Equation~(62), and taking into account 
Equation~(56), one shows that
\begin{equation*}
\mathcal{H}(\tau) \! \underset{\tau \to +\infty}{=} \! 3(\varepsilon b)^{2/3} 
\tau^{1/3} \! - \! 4h_{0}(\tau) \tau^{-1/3} \! + \! \dfrac{(a \! - \! \mi/2)^{
2}}{2 \tau}:
\end{equation*}
now, the asymptotics of $\mathcal{H}(\tau)$ given in Theorem~3.1, 
Equation~(40) follows {}from Equations~(87) and~(103).
\section{Calculation of the Monodromy Data as $\tau \! \to \! +0$}
In this section, the first equation of System~(\ref{eq:UV}) is studied 
asymptotically as $\tau \! \to \! +0$ under specific conditions on the 
elements of $\mathscr{U}(\lambda,\tau)$. Fundamental solutions of this 
equation in the neighbourhoods of the essential singularities $(0$ and 
$\infty)$ are approximated in terms of Hankel and Whittaker functions: this 
allows one to calculate, asymptotically, the corresponding monodromy data.

Denoting by $\mathscr{U}_{0}$ and $\mathscr{V}_{0}$ the following matrices,
\begin{equation*}
\mathscr{U}_{0} \! = \! \tau \! \left(-\mi \sigma_{3} \! - \! \dfrac{a \mi}{2 
\tau \lambda} \sigma_{3} \! - \! \frac{1}{\lambda} \! 
\begin{pmatrix}
0 & C \\
D & 0
\end{pmatrix} \right), \qquad \quad \mathscr{V}_{0} \! = \! \dfrac{\mi \tau}{2 
\lambda^{2}} \! 
\begin{pmatrix}
\sqrt{\smash[b]{-AB}} & A \\
B & -\sqrt{\smash[b]{-AB}}
\end{pmatrix},
\end{equation*}
the first equation of System~(\ref{eq:UV}) can be rewritten as
\begin{equation}
\partial_{\lambda} \Phi (\lambda) \! = \! (\mathscr{U}_{0} \! + \! \mathscr{
V}_{0}) \Phi (\lambda).
\end{equation}
\begin{proposition}
Consider
\begin{equation}
\partial_{\lambda} \mathbf{W}(\lambda) \! = \! \mathscr{U}_{0} \mathbf{W}
(\lambda).
\end{equation}
A fundamental solution of Equation~{\rm (109)} is given by
\begin{equation}
\mathbf{W}(\lambda) \! = \! \dfrac{\me^{-\frac{\pi a}{4}}}{\sqrt{\smash[b]{2 
\mi \lambda \tau}}} \! 
\begin{pmatrix}
W_{\varkappa_{1},\widehat{\rho}}(2 \mi \lambda \tau) & \mi \widehat{\gamma}
W_{-\varkappa_{1},\widehat{\rho}}(-2 \mi \lambda \tau) \\
\widehat{\delta} \, W_{\varkappa_{1}-1,\widehat{\rho}}(2 \mi \lambda \tau) & 
\mi W_{-(\varkappa_{1}-1),\widehat{\rho}}(-2 \mi \lambda \tau)
\end{pmatrix} \! \me^{\mi a \ln (\sqrt{2 \tau}) \sigma_{3}},
\end{equation}
where
\begin{equation}
\varkappa_{1} \! := \! \dfrac{1}{2}(1 \! - \! a \mi), \qquad \quad \widehat{
\rho}^{2} \! := \! \widehat{\gamma} \, \widehat{\delta} \! - \! \dfrac{a^{2}}{
4}, \qquad \quad \widehat{\gamma} \! := \! \tau C, \qquad \quad \widehat{
\delta} \! := \! \tau D,
\end{equation}
and $W_{\star,\widehat{\rho}}(\boldsymbol{\cdot})$ is the Whittaker function 
{\rm \cite{a24};} moreover, $\det (\mathbf{W}(\lambda)) \! = \! 1$, and
\begin{equation}
\mathbf{W}(\lambda) \! \underset{\underset{\arg \lambda =0}{\lambda \to 
\infty}}{=} \! \left(\mathrm{I} \! + \! \dfrac{1}{2 \mi \tau \lambda} \! 
\begin{pmatrix}
\widehat{\gamma} \, \widehat{\delta} & -\widehat{\gamma} \\
\widehat{\delta} & -\widehat{\gamma} \, \widehat{\delta}
\end{pmatrix} \! + \! \mathcal{O} \! \left(\dfrac{1}{\lambda^{2}} \right) 
\right) \! \me^{-\mi (\tau \lambda +\frac{a}{2} \ln \lambda) \sigma_{3}}.
\end{equation}
\end{proposition}

\emph{Proof.} Rewriting Equation~(109) as a system for its components, one 
deduces {}from it, for $(\mathbf{W}(\lambda))_{11}$ and $(\mathbf{W}(\lambda)
)_{22}$, the Whittaker ODE: then, $(\mathbf{W}(\lambda))_{12}$ and $(\mathbf{
W}(\lambda))_{21}$ are obtained {}from the latter system by applying certain 
identities for the Whittaker function. \hfill $\square$
\begin{proposition}
Let $\epsilon_{i} \! > \! 0$, $i \! = \! 1,2$, and the parameters of 
Equation~{\rm (108)} satisfy the following restrictions:
\begin{equation}
\begin{gathered}
\vert \Im (a) \vert \! < \! 1, \qquad \qquad \widehat{\rho} \! \underset{\tau 
\to +0}{=} \! \mathcal{O}(1), \qquad \qquad A \tau^{1+\mi a} \! \underset{\tau 
\to +0}{\sim} \! \tau^{\epsilon_{1}}, \\
B \tau^{1-\mi a} \! \underset{\tau \to +0}{\sim} \! \tau^{\epsilon_{1}}, 
\qquad \qquad \vert C \vert \! \underset{\tau \to +0}{>} \! \vert A \vert 
\tau^{-\epsilon_{1}}, \qquad \qquad \vert D \vert \! \underset{\tau \to +0}{>} 
\! \vert B \vert \tau^{-\epsilon_{1}}.
\end{gathered}
\end{equation}
Then there exists a fundamental solution of Equation~{\rm (108)} which has the 
asymptotic representation
\begin{equation}
\Phi (\lambda) \! \underset{\underset{\lambda > \epsilon_{2}, \, \, \arg 
\lambda =0}{\tau \to +0}}{=} \! \mathbf{W}(\lambda) \! \left(\mathrm{I} \! + 
\! o \! \left(\! \left(\dfrac{\tau}{\lambda} \right)^{\delta_{1}} \right) 
\right),
\end{equation}
with $\delta_{1} \! > \! 0$.
\end{proposition}

\emph{Proof.} Follows {}from a successive approximations argument applied to 
Equation~(108). \hfill $\square$

For the purpose of approximating the solution of the first equation of 
System~(\ref{eq:UV}) as $\lambda \! \to \! 0$, it is convenient to rewrite 
it as
\begin{equation}
\partial_{\lambda} \Phi (\lambda) \! = \! (\widetilde{\mathscr{U}}_{0}
\! + \! \widetilde{\mathscr{V}}_{0}) \Phi (\lambda),
\end{equation}
where
\begin{equation*}
\widetilde{\mathscr{U}}_{0}= \tau \! \left(\! -\dfrac{a \mi}{2 \tau \lambda} 
\sigma_{3} \! - \! \dfrac{1}{\lambda} \! 
\begin{pmatrix}
0 & C \\
D & 0
\end{pmatrix} \! + \! \dfrac{\mi}{2 \lambda^{2}} \! 
\begin{pmatrix}
\sqrt{\smash[b]{-AB}} & A \\
B & -\sqrt{\smash[b]{-AB}}
\end{pmatrix} \right), \qquad \quad \widetilde{\mathscr{V}}_{0} \! = \! -\mi 
\tau \sigma_{3}:
\end{equation*}
this equation is compared to the following model system,
\begin{equation}
\partial_{\lambda} \mathfrak{B}(\lambda) \! = \! \widetilde{\mathscr{U}}_{0} 
\mathfrak{B}(\lambda).
\end{equation}
\begin{proposition}
The fundamental solution of Equation~{\rm (116)} can be written as follows:
\begin{equation}
\mathfrak{B}(\lambda) \! = \! 
\begin{pmatrix}
\mathfrak{B}_{11}(\lambda) & \mathfrak{B}_{12}(\lambda) \\
\mathfrak{B}_{21}(\lambda) & \mathfrak{B}_{22}(\lambda)
\end{pmatrix},
\end{equation}
where
\begin{align*}
\mathfrak{B}_{11}(\lambda) \! &= \! -\dfrac{\sqrt{\smash[b]{-AB}} \, \sqrt{
\smash[b]{\pi}} \, \me^{-\mi (\frac{\pi \nu_{o}}{2}+\frac{\pi}{4})}}{2 \sqrt{
\smash[b]{\varepsilon b}} \, \sqrt{\smash[b]{B}}} \! \left(\mathfrak{x}^{
\uparrow}(\tau)H_{\nu_{o}}^{(2)} \! \left(\! \sqrt{\frac{\tau \varepsilon b}
{\lambda}} \right) \! + \! \sqrt{\frac{\tau \varepsilon b}{\lambda}}H_{\nu_{o}
-1}^{(2)} \! \left(\! \sqrt{\frac{\tau \varepsilon b}{\lambda}} \right) 
\right), \\
\mathfrak{B}_{12}(\lambda) \! &= \! -\dfrac{\sqrt{\smash[b]{-AB}} \, \sqrt{
\smash[b]{\pi}} \, \me^{\mi (\frac{\pi \nu_{o}}{2}+\frac{\pi}{4})}}{2 \sqrt{
\smash[b]{\varepsilon b}} \, \sqrt{\smash[b]{B}}} \! \left(\mathfrak{x}^{
\uparrow}(\tau)H_{\nu_{o}}^{(1)} \! \left(\! \sqrt{\frac{\tau \varepsilon b}
{\lambda}} \right) \! + \! \sqrt{\frac{\tau \varepsilon b}{\lambda}}H_{\nu_{o}
-1}^{(1)} \! \left(\! \sqrt{\frac{\tau \varepsilon b}{\lambda}} \right) 
\right), \\
\mathfrak{B}_{21}(\lambda) \! &= \! -\dfrac{\sqrt{\smash[b]{B}} \, \sqrt{
\smash[b]{\pi}} \, \me^{-\mi (\frac{\pi \nu_{o}}{2}+\frac{\pi}{4})}}{2 \sqrt{
\smash[b]{\varepsilon b}}} \! \left(\mathfrak{x}^{\downarrow}(\tau)H_{\nu_{o}
}^{(2)} \! \left(\! \sqrt{\frac{\tau \varepsilon b}{\lambda}} \right) \! + \! 
\sqrt{\frac{\tau \varepsilon b}{\lambda}}H_{\nu_{o}-1}^{(2)} \! \left(\! 
\sqrt{\frac{\tau \varepsilon b}{\lambda}} \right) \right), \\
\mathfrak{B}_{22}(\lambda) \! &= \! -\dfrac{\sqrt{\smash[b]{B}} \, \sqrt{
\smash[b]{\pi}} \, \me^{\mi (\frac{\pi \nu_{o}}{2}+\frac{\pi}{4})}}{2 \sqrt{
\smash[b]{\varepsilon b}}} \! \left(\mathfrak{x}^{\downarrow}(\tau)H_{\nu_{o}
}^{(1)} \! \left(\! \sqrt{\frac{\tau \varepsilon b}{\lambda}} \right) \! + \! 
\sqrt{\frac{\tau \varepsilon b}{\lambda}}H_{\nu_{o}-1}^{(1)} \! \left(\! 
\sqrt{\frac{\tau \varepsilon b}{\lambda}} \right) \right),
\end{align*}
with
\begin{equation}
\mathfrak{x}^{\uparrow}(\tau) \! := \! -\nu_{o} \! + \! a \mi \! + \! \dfrac{
2B \widehat{\gamma}}{\sqrt{\smash[b]{-AB}}}, \quad \quad \mathfrak{x}^{
\downarrow}(\tau) \! := \! -\nu_{o} \! - \! a \mi \! + \! \dfrac{2 \widehat{
\delta} \sqrt{\smash[b]{-AB}}}{B}, \quad \quad \nu_{o}^{2} \! := \! 4 
\widehat{\rho}^{2},
\end{equation}
and $H^{(j)}_{\star}(\boldsymbol{\cdot})$, $j \! = \! 1,2$, the Hankel 
functions of the first $(j \! = \! 1)$ and second $(j \! = \! 2)$ kind 
{\rm \cite{a24}}.

Furthermore, $\det (\mathfrak{B}(\lambda)) \! = \! 1$, and $\mathfrak{B}
(\lambda)$ has the asymptotic expansion
\begin{align}
\mathfrak{B}(\lambda) \underset{\underset{\arg \lambda =0}{\lambda \to 
+0}}{=}& \, \left(\! \dfrac{\mi \tau^{1/4} \lambda^{-1/4}}{\sqrt{\smash[b]{2}} 
\, (\varepsilon b)^{1/4}} \! 
\begin{pmatrix}
\frac{\sqrt{\smash[b]{-AB}}}{\sqrt{\smash[b]{B}}} & -\frac{\sqrt{\smash[b]{-
AB}}}{\sqrt{\smash[b]{B}}} \\
\sqrt{\smash[b]{B}} & -\sqrt{\smash[b]{B}}
\end{pmatrix} \! + \! \dfrac{\lambda^{1/4} \tau^{-1/4}}{\sqrt{\smash[b]{2}} 
\, (\varepsilon b)^{3/4}} \right. \nonumber \\
\times&\left. 
\begin{pmatrix}
-\frac{\sqrt{\smash[b]{-AB}}}{\sqrt{\smash[b]{B}}} \! \left(\mathfrak{x}^{
\uparrow}(\tau)-\frac{(4(\nu_{o}-1)^{2}-1)}{8} \right) & -\frac{\sqrt{
\smash[b]{-AB}}}{\sqrt{\smash[b]{B}}} \! \left(\mathfrak{x}^{\uparrow}(\tau)
-\frac{(4(\nu_{o}-1)^{2}-1)}{8} \right) \\
-\sqrt{\smash[b]{B}} \left(\mathfrak{x}^{\uparrow}(\tau)-\frac{(4(\nu_{o}-1)^{
2}-1)}{8} \right)-\frac{\mi \varepsilon b \sqrt{\smash[b]{B}}}{\sqrt{\smash[b]{
-AB}}} & -\sqrt{\smash[b]{B}} \left(\mathfrak{x}^{\uparrow}(\tau)-\frac{(4(
\nu_{o}-1)^{2}-1)}{8} \right)-\frac{\mi \varepsilon b \sqrt{\smash[b]{B}}}{
\sqrt{\smash[b]{-AB}}}
\end{pmatrix} \right. \nonumber \\
+&\left. \, \mathcal{O} \! \left(\dfrac{\lambda^{3/4}}{\tau^{3/4}} \right) 
\right) \exp \! \left(\! -\mi \sqrt{\frac{\tau \varepsilon b}{\lambda}} \, 
\sigma_{3} \right).
\end{align}
\end{proposition}

\emph{Proof.} Defining the functions $\widetilde{\Phi}^{(0)}_{\uparrow 
\downarrow}(x)$ by $\mathfrak{B}(x) \! := \! \widehat{\phi}^{0}_{\uparrow 
\downarrow} \widetilde{\Phi}^{(0)}_{\uparrow \downarrow}(x)$, where $x \! := 
\! \mi \tau/\lambda$, $\widehat{\phi}^{0}_{\uparrow} \! = \! 
\left(
\begin{smallmatrix}
1 & 0 \\
\frac{B}{\sqrt{\smash[b]{-AB}}} & 1
\end{smallmatrix} \right)$, and $\widehat{\phi}^{0}_{\downarrow} \! = \! 
\left(
\begin{smallmatrix}
1 & \frac{\sqrt{\smash[b]{-AB}}}{B} \\
0 & 1
\end{smallmatrix}
\right)$, one arrives at
\begin{align}
\partial_{x} \widetilde{\Phi}^{(0)}_{\uparrow}(x) \! &= \!
\left(\dfrac{a \mi}{2x} \! 
\begin{pmatrix}
1 & 0 \\
-\frac{2B}{\sqrt{\smash[b]{-AB}}} & -1
\end{pmatrix} \! + \! \dfrac{\tau}{x} \! 
\begin{pmatrix}
\frac{BC}{\sqrt{\smash[b]{-AB}}} & C \\
\frac{(AD+BC)}{A} & -\frac{BC}{\sqrt{\smash[b]{-AB}}}
\end{pmatrix} \! - \! \dfrac{1}{2} \! 
\begin{pmatrix}
0 & A \\
0 & 0
\end{pmatrix} \right) \! \widetilde{\Phi}^{(0)}_{\uparrow}(x), \\
\partial_{x} \widetilde{\Phi}^{(0)}_{\downarrow}(x) \! &= \! \left(
\dfrac{a \mi}{2x} \! 
\begin{pmatrix}
1 & \frac{2 \sqrt{\smash[b]{-AB}}}{B} \\
0 & -1
\end{pmatrix} \! + \! \dfrac{\tau}{x} \! 
\begin{pmatrix}
-\frac{D \sqrt{\smash[b]{-AB}}}{B} & \frac{(AD+BC)}{B} \\
D & \frac{D \sqrt{\smash[b]{-AB}}}{B}
\end{pmatrix} \! - \! \dfrac{1}{2} \! 
\begin{pmatrix}
0 & 0 \\
B & 0
\end{pmatrix} \right) \! \widetilde{\Phi}^{(0)}_{\downarrow}(x).
\end{align}
Consider, for example, Equation~(120). Rewriting it in component form, one 
obtains (modulo the change $x \! \to \! (-\mi \varepsilon b x)^{1/2})$ for 
them the Bessel equation. Choosing
\begin{equation*}
\left(\widetilde{\Phi}^{(0)}_{\uparrow}(\lambda) \right)_{21} \! = \! c_{21}
H_{\nu_{o}}^{(2)} \! \left(\! \sqrt{\smash[b]{(\tau \varepsilon b)/\lambda}} 
\right) \qquad \text{and} \qquad \left(\widetilde{\Phi}^{(0)}_{\uparrow}
(\lambda) \right)_{22} \! = \! c_{22}H_{\nu_{o}}^{(1)} \! \left(\! \sqrt{
\smash[b]{(\tau \varepsilon b)/\lambda}} \right),
\end{equation*}
where $c_{21} \! = \! -\tfrac{\mi \sqrt{\smash[b]{\pi \varepsilon 
b}} \, \sqrt{\smash[b]{B}} \, \me^{-\mi (\frac{\pi \nu_{o}}{2}+\frac{\pi}{4})}
}{2 \sqrt{\smash[b]{-AB}}}$, $c_{22} \! = \! -\tfrac{\mi \sqrt{\smash[b]{\pi 
\varepsilon b}} \, \sqrt{\smash[b]{B}} \, \me^{\mi (\frac{\pi \nu_{o}}{2}+
\frac{\pi}{4})}}{2 \sqrt{\smash[b]{-AB}}}$, and $H^{(j)}_{\star}(\cdot)$ are 
the Hankel functions of the first $(j \! = \! 1)$ and second $(j \! = \! 2)$ 
kind, one obtains
\begin{equation*}
\left(\widetilde{\Phi}^{(0)}_{\uparrow}(\lambda) \right)_{11} \! = \! \dfrac{
\mi c_{21}A}{\varepsilon b} \! \left(\mathfrak{x}^{\uparrow}(\tau)H_{\nu_{o}
}^{(2)} \! \left(\! \sqrt{\smash[b]{(\tau \varepsilon b)/\lambda}} \right) \! 
+ \! \sqrt{\smash[b]{(\tau \varepsilon b)/\lambda}} \, H_{\nu_{o}-1}^{(2)} \! 
\left(\! \sqrt{\smash[b]{(\tau \varepsilon b)/\lambda}} \right) \right)
\end{equation*}
and
\begin{equation*}
\left(\widetilde{\Phi}^{(0)}_{\uparrow}(\lambda) \right)_{12} \! = \! \dfrac{
\mi c_{22}A}{\varepsilon b} \! \left(\mathfrak{x}^{\uparrow}(\tau)H_{\nu_{o}
}^{(1)} \! \left(\! \sqrt{\smash[b]{(\tau \varepsilon b)/\lambda}} \right) \! 
+ \! \sqrt{\smash[b]{(\tau \varepsilon b)/\lambda}} \, H_{\nu_{o}-1}^{(1)} \! 
\left(\! \sqrt{\smash[b]{(\tau \varepsilon b)/\lambda}} \right) \right)
\end{equation*}
{}from the corresponding component equations of Equation~(120) and the 
identity $z \partial_{z}H_{\star}^{(j)}(z) \! = \! -\star H_{\star}^{(j)}(z) 
\! + \! zH_{\star -1}^{(j)}(z)$, $j \! = \! 1,2$. The components $(\widetilde{
\Phi}^{(0)}_{\uparrow}(\lambda))_{21}$ and $(\widetilde{\Phi}^{(0)}_{\uparrow}
(\lambda))_{22}$ are chosen in such a way that the leading term of asymptotics 
of $\mathfrak{B}(\lambda)$ as $\lambda \! \to \! +0$ (Equation~(119)) matches 
with the canonical asymptotics given in Equation~(16): to verify this, one 
uses the following asymptotic expansions for the Hankel functions \cite{a27},
\begin{align*}
H_{\nu_{o}}^{(1)}(z) &\underset{\underset{\vert \arg z \vert < \pi}{z \to 
\infty}}{=} \sqrt{\dfrac{2}{\pi z}} \, \me^{\mi (z-\frac{\pi \nu_{o}}{2}-
\frac{\pi}{4})} \! \left(1 \! + \! \dfrac{\mi (4 \nu_{o}^{2} \! - \! 1)}{8z} 
\! + \! \mathcal{O} \! \left(\dfrac{1}{z^{2}} \right) \right), \\
H_{\nu_{o}}^{(2)}(z) &\underset{\underset{\vert \arg z \vert < \pi}{z \to 
\infty}}{=} \sqrt{\dfrac{2}{\pi z}} \, \me^{-\mi (z-\frac{\pi\nu_{o}}{2}-
\frac{\pi}{4})} \! \left(1 \! - \! \dfrac{\mi (4 \nu_{o}^{2} \! - \! 1)}{8z} 
\! + \! \mathcal{O} \! \left(\dfrac{1}{z^{2}} \right) \right).
\end{align*}
Noting that $\det (\mathfrak{B}(\lambda)) \! = \! \det (\widetilde{\Phi}^{
(0)}_{\uparrow}(\lambda))$ and the trace of the coefficient matrix of 
Equation~(120) is zero, one shows that $\det (\mathfrak{B}(\lambda)) \! = \! 
\mathrm{const.}$: one proves that $\mathrm{const.} \! = \! 1$ by calculating 
$\det (\mathfrak{B}(\lambda))$ as $\lambda \! \to \! +0$ by means of 
asymptotics~(119). \hfill $\square$
\begin{proposition}
Let $0 \! < \! \epsilon_{i} \! < \! 1$, $i \! = \! 3,4$. If the parameters of 
Equation~{\rm (115)} satisfy the restrictions
\begin{equation}
\nu_{o} \! \underset{\tau \to +0}{=} \! \mathcal{O}(1), \qquad \qquad 
\left\vert AD \! + \! BC \right\vert \! \underset{\tau \to +0}{>} \! 
\left\vert \sqrt{-AB} \right\vert \! \tau^{-\epsilon_{3}},
\end{equation}
and either
\begin{equation}
\tau^{-\epsilon_{3}} \! \underset{\tau \to +0}{<} \! \dfrac{\vert BC \vert}{
\left\vert \sqrt{\smash[b]{-AB}} \right\vert} \! \underset{\tau \to +0}{
\leqslant} \! \mathcal{O} \! \left(\dfrac{1}{\tau} \right),
\end{equation}
or
\begin{equation}
\tau^{-\epsilon_{3}} \! \underset{\tau \to +0}{<} \! \dfrac{D \sqrt{\smash[b]{
-AB}}}{B} \! \underset{\tau \to +0}{\leqslant} \! \mathcal{O} \! \left(\dfrac{
1}{\tau} \right),
\end{equation}
then there exists a fundamental solution of Equation~{\rm (115)} which has the 
asymptotic representation
\begin{equation}
\Phi (\lambda) \! \underset{\underset{0< \lambda < \epsilon_{4}, \, \, \, \arg 
\lambda =0}{\tau \to +0}}{=} \! \mathfrak{B}(\lambda) \! \left(\mathrm{I} \! + 
\! o \! \left(\! \left(\tau \lambda \right)^{\delta} \right) \right),
\end{equation}
with some $\delta \! > \! 0$.
\end{proposition}

\emph{Proof.} Consider, for example, the restrictions~(122) and~(123) imposed 
on the coefficients of Equation~(115). Consider the transformation $\Phi 
(\lambda) \! = \! \widehat{\phi}^{0}_{\uparrow} \widetilde{\Phi}_{\uparrow}
(\lambda)$, where $\widehat{\phi}^{0}_{\uparrow}$ is defined in the proof of 
Proposition~5.3. Note that, under the assumed conditions, $(\widehat{\phi}^{
0}_{\uparrow})^{-1} \widetilde{\mathscr{V}}_{0} \widehat{\phi}^{0}_{\uparrow} 
\! =_{\underset{0< \lambda < \epsilon_{4}}{\tau \to +0}} \! \mathcal{O} \! 
\left(\tau^{-\epsilon}(\widehat{\phi}^{0}_{\uparrow})^{-1} \widetilde{\mathscr{
U}}_{0} \widehat{\phi}^{0}_{\uparrow} \right)$, with some $\epsilon_{4},
\epsilon \! > \! 0$: now, using a successive approximations argument, one 
shows that $\widetilde{\Phi}_{\uparrow}(\lambda) \! =_{\underset{0< \lambda < 
\epsilon_{4}, \, \, \, \arg \lambda =0}{\tau \to +0}} \! \widetilde{\Phi}^{
(0)}_{\uparrow}(\lambda) \! \left(1 \! + \! o((\tau \lambda)^{\delta}) 
\right)$. For the restrictions~(122) and~(124), one proceeds as above, but, 
instead of for $\widetilde{\Phi}_{\uparrow}(\lambda)$, for the function 
$\widetilde{\Phi}_{\downarrow}(\lambda) \! = \! (\widehat{\phi}^{0}_{
\downarrow})^{-1} \Phi (\lambda)$, where $\widehat{\phi}^{0}_{\downarrow}$ 
is defined in the proof of Proposition~5.3. \hfill $\square$
\begin{ccc}
Let $0 \! < \! \epsilon \! < \! 1$,
\begin{gather*}
\widehat{\rho} \! \underset{\tau \to +0}{\not=} \! 0, \qquad \qquad \left\vert 
\Re (\widehat{\rho}) \right\vert \! \underset{\tau \to +0}{<} \! \dfrac{1}{2}, 
\\
\left\vert \sqrt{\smash[b]{-AB}} \right\vert \! + \! \left\vert A/\widehat{
\gamma} \right\vert \! \underset{\tau \to +0}{>} \! \left(\left\vert B \, 
\widehat{\gamma} \right\vert \! + \! \mathcal{O}(1) \right) \! \tau^{2-4 \vert 
\Re (\widehat{\rho}) \vert}, \\
\left\vert \sqrt{\smash[b]{-AB}} \right\vert \! + \! \left\vert A \, \widehat{
\delta} \right\vert \! \underset{\tau \to +0}{>} \! \left(\left\vert B \, 
\widehat{\gamma} \right\vert \! + \! \mathcal{O}(1) \right) \! \tau^{2-4 \vert 
\Re (\widehat{\rho}) \vert}, \\
\left\vert \sqrt{\smash[b]{-AB}} \right\vert \! + \! \left\vert B/\widehat{
\delta} \right\vert \! \underset{\tau \to +0}{>} \! \left(\left\vert A 
\widehat{\delta} \right\vert \! + \! \mathcal{O}(1) \right) \! \tau^{2-4 \vert 
\Re (\widehat{\rho}) \vert}, \\
\left\vert \sqrt{\smash[b]{-AB}} \right\vert \! + \! \left\vert B \, \widehat{
\gamma} \right\vert \! \underset{\tau \to +0}{>} \! \left(\left\vert A 
\widehat{\delta} \right\vert \! + \! \mathcal{O}(1) \right) \! \tau^{2-4 \vert 
\Re (\widehat{\rho}) \vert}, \\
\left\vert \widehat{\gamma} \sqrt{B} \, \tau^{-\frac{a \mi}{2}} \right\vert \! 
\underset{\tau \to +0}{\leqslant} \! \mathcal{O} \! \left(\tau^{-2+2 \vert \Re 
(\widehat{\rho}) \vert +\epsilon} \right), \qquad \qquad \left\vert \sqrt{B} 
\, \tau^{\frac{a \mi}{2}} \right\vert \! \underset{\tau \to +0}{\leqslant} \! 
\mathcal{O} \! \left(\tau^{-2+2 \vert \Re (\widehat{\rho}) \vert+\epsilon} 
\right).
\end{gather*}
Then, under conditions~{\rm (113)}, {\rm (122)} and {\rm (123)}, or 
{\rm (113)}, {\rm (122)} and~{\rm (124)}, the connection matrix, $G$ 
(Equation~{\rm (28))}, has the following asymptotics $(g_{ij} \! := \! (G)_{i
j}$, $i,j \! = \! 1,2)$,
\begin{align}
g_{11} \underset{\tau \to +0}{=}& \, \left(\, \sum\limits_{l = \pm 1} \! 
\dfrac{2^{l \widehat{\rho}+1} \me^{\frac{\mi \pi l \widehat{\rho}}{2}}(\Gamma 
(2l \widehat{\rho}))^{2} \tau^{-2l \widehat{\rho}}(2 \tau)^{\frac{a \mi}{2}} 
\me^{-\frac{\mi \pi}{4}}l \widehat{\rho}}{\sqrt{\smash[b]{\pi \varepsilon b}} 
\, (\varepsilon b)^{l \widehat{\rho}} \! \left(l \widehat{\rho} \! + \! \frac{
a \mi}{2} \right) \! \me^{\frac{\pi a}{4}} \Gamma \! \left(l \widehat{\rho} \! 
+ \! \frac{a \mi}{2} \right)} \! \left(\! \sqrt{\smash[b]{B}} \left(l \widehat{
\rho} \! + \! \dfrac{a \mi}{2} \right) \! - \! \dfrac{\widehat{\delta} \sqrt{
\smash[b]{-AB}}}{\sqrt{\smash[b]{B}}} \, \right) \right. \nonumber \\
+&\left. \, o \! \left(\tau^{\delta_{1}} \right) \right) \! \left(1 \! + \! 
o \! \left(\tau^{\delta_{2}} \right) \right), \\
g_{12} \underset{\tau \to +0}{=}& \, \left(\, \sum\limits_{l = \pm 1} \! 
\dfrac{2^{l \widehat{\rho}+1} \me^{\frac{3 \pi \mi l \widehat{\rho}}{2}}
(\Gamma (2l \widehat{\rho}))^{2} \tau^{-2l \widehat{\rho}}(2 \tau)^{-\frac{a 
\mi}{2}} \me^{-\frac{\mi \pi}{4}}l \widehat{\rho}}{\sqrt{\smash[b]{\pi 
\varepsilon b}} \, (\varepsilon b)^{l \widehat{\rho}} \! \left(l \widehat{
\rho} \! - \! \frac{a \mi}{2} \right) \! \me^{\frac{\pi a}{4}} \Gamma \! 
\left(l \widehat{\rho} \! - \! \frac{a \mi}{2} \right)} \! \left(\! \sqrt{
\smash[b]{B}} \, \widehat{\gamma} \! - \! \dfrac{\sqrt{\smash[b]{-AB}}}{\sqrt{
\smash[b]{B}}} \! \left(l \widehat{\rho} \! - \! \dfrac{a \mi}{2} \right) 
\right) \right. \nonumber \\
+&\left. \, o \! \left(\tau^{\delta_{1}} \right) \right) \! \left(1 \! + \! 
o \! \left(\tau^{\delta_{2}} \right) \right), \\
g_{21} \underset{\tau \to +0}{=}& \, \left(-\sum\limits_{l = \pm 1} \! 
\dfrac{2^{l \widehat{\rho}+1} \me^{-\frac{3 \pi \mi l \widehat{\rho}}{2}}
(\Gamma (2l \widehat{\rho}))^{2} \tau^{-2l \widehat{\rho}}(2 \tau)^{\frac{a 
\mi}{2}} \me^{\frac{\mi \pi}{4}}l \widehat{\rho}}{\sqrt{\smash[b]{\pi 
\varepsilon b}} \, (\varepsilon b)^{l \widehat{\rho}} \! \left(l \widehat{
\rho} \! + \! \frac{a \mi}{2} \right) \! \me^{\frac{\pi a}{4}} \Gamma \! 
\left(l \widehat{\rho} \! + \! \frac{a \mi}{2} \right)} \! \left(\! \sqrt{
\smash[b]{B}} \left(l \widehat{\rho} \! + \! \dfrac{a \mi}{2} \right) \! - \! 
\dfrac{\widehat{\delta} \sqrt{\smash[b]{-AB}}}{\sqrt{\smash[b]{B}}} \right) 
\right. \nonumber \\
+&\left. \, o \! \left(\tau^{\delta_{1}} \right) \right) \! \left(1 \! + \! 
o \! \left(\tau^{\delta_{2}} \right) \right), \\
g_{22} \underset{\tau \to +0}{=}& \, \left(-\sum\limits_{l = \pm 1} \! 
\dfrac{2^{l \widehat{\rho}+1} \me^{-\frac{\mi \pi l \widehat{\rho}}{2}}(\Gamma 
(2l \widehat{\rho}))^{2} \tau^{-2l \widehat{\rho}}(2 \tau)^{-\frac{a \mi}{2}} 
\me^{\frac{\mi \pi}{4}}l \widehat{\rho}}{\sqrt{\smash[b]{\pi \varepsilon b}} 
\, (\varepsilon b)^{l \widehat{\rho}} \! \left(l \widehat{\rho} \! - \! \frac{
a \mi}{2} \right) \! \me^{\frac{\pi a}{4}} \Gamma \! \left(l \widehat{\rho} \! 
- \! \frac{a \mi}{2} \right)} \! \left(\! \sqrt{\smash[b]{B}} \, \widehat{
\gamma} \! - \! \dfrac{\sqrt{\smash[b]{-AB}}}{\sqrt{\smash[b]{B}}} \left(l 
\widehat{\rho} \! - \! \dfrac{a \mi}{2} \right) \right) \right. \nonumber \\
+&\left. \, o \! \left(\tau^{\delta_{1}} \right) \right) \! \left(1 \! + \! 
o \! \left(\tau^{\delta_{2}} \right) \right),
\end{align}
with $\delta_{j} \! > \! 0$, $j \! = \! 1,2$.
\end{ccc}

\emph{Proof.} Consider, say, the derivation of the results for $G$ 
corresponding to conditions~(113), (122) and~(123): the derivation for the 
case corresponding to conditions~(113), (122) and~(124) is analogous. For this 
purpose, the results stated in Propositions~5.1 and~5.3 are applicable; 
therefore,
\begin{equation}
G \! \underset{\tau \to +0}{=} \! \left(1 \! + \! o \! \left((\tau \lambda)^{
\delta_{2}} \right) \right) \! \left(\mathfrak{B}(\lambda) \right)^{-1} 
\mathbf{W}(\lambda) \! \left(1 \!+ \! o \! \left(\! \left(\tau / \lambda 
\right)^{\delta_{1}} \right) \right),
\end{equation}
with $\delta_{j} \! > \! 0$, $j \! = \! 1,2$. This formula allows one to 
calculate $G$ to the desired order by taking (real) $\lambda$ in the finite 
domain $\epsilon_{2} \! < \! \lambda \! < \! \epsilon_{1}$ and applying 
small-argument asymptotic expansions of the Whittaker and Hankel functions. 
Towards this end, one notes that \cite{a24} $H_{\star}^{(1)}(z) \! = \! (2/\pi 
z)^{1/2} \, \exp \! \left(-\mi \! \left(\tfrac{\pi \star}{2} \! + \! \tfrac{
\pi}{4} \right) \right) \! W_{0,\star}(-2 \mi z)$ and $H_{\star}^{(2)}(z) \! 
= \! (2/\pi z)^{1/2} \, \exp \! \left(\mi \! \left(\tfrac{\pi \star}{2} \! 
+ \! \tfrac{\pi}{4} \right) \right) \! W_{0,\star}(2 \mi z)$; he\-n\-c\-e, 
as a consequence, one needs the small-argument asymptotic expansion of the 
Whittaker function:
\begin{equation}
W_{z_{1},z_{2}}(z) \! \underset{\underset{\vert \arg z \vert < \pi}{z \to 0}
}{=} \! z^{1/2} \sum_{k=0}^{2} \sum_{l = \pm 1}w_{k}(z_{1},lz_{2})z^{k} \! + 
\! \mathcal{O} \! \left(z^{3}z^{-\vert \Re (z_{2}) \vert} \right),
\end{equation}
where
\begin{gather*}
w_{0}(z_{1},lz_{2}) \! = \! \dfrac{\Gamma (-2lz_{2})z^{lz_{2}}}{\Gamma (\frac{
1}{2} \! - \! lz_{2} \! - \! z_{1})}, \qquad \qquad \dfrac{w_{1}(z_{1},lz_{2})
}{w_{0}(z_{1},lz_{2})} \! = \! -\dfrac{1}{2} \! + \! \dfrac{(lz_{2} \! - \! 
z_{1} \! + \! \frac{1}{2})}{(2lz_{2} \! + \! 1)}, \\
\dfrac{w_{2}(z_{1},lz_{2})}{w_{0}(z_{1},lz_{2})} \! = \! \dfrac{1}{8} \! - \! 
\dfrac{(lz_{2} \! - \! z_{1} \! + \! \frac{1}{2})}{2(2lz_{2} \! + \! 1)} \! + 
\! \dfrac{(lz_{2} \! - \! z_{1} \! + \! \frac{1}{2})(lz_{2} \! - \! z_{1} \! 
+ \! \frac{3}{2})}{2(2lz_{2} \! + \! 1)(2lz_{2} \! + \! 2)}.
\end{gather*}
The latter expansion is obtained by using a representation of $W_{z_{1},z_{2}}
(\cdot)$ in terms of the confluent hypergeometric function \cite{a27}, 
$\boldsymbol{\Phi}(\alpha,\beta;z) \! := \! \sum_{n=0}^{\infty} \tfrac{
(\alpha)_{n}z^{n}}{(\beta)_{n}n!}$:
\begin{align}
W_{z_{1},z_{2}}(z) \! &= \me^{-\frac{z}{2}}z^{\frac{1}{2}+z_{2}} \! \left(
\dfrac{\Gamma (-2z_{2})}{\Gamma (\frac{1}{2} \! - \! z_{2} \! - \! z_{1})} 
\boldsymbol{\Phi}(z_{2} \! - \! z_{1} \! + \! \tfrac{1}{2},2z_{2} \! + \! 1;z) 
\! + \! \dfrac{\Gamma (2z_{2})}{\Gamma (\frac{1}{2} \! + \! z_{2} \! - \! z_{
1})} \right. \nonumber \\
&\left. \times \, z^{-2z_{2}} \boldsymbol{\Phi}(-z_{2} \! - \! z_{1} \! + \! 
\tfrac{1}{2},-2z_{2} \! + \! 1;z) \right).
\end{align}
Consider, for example, $g_{11} \! := \! (G)_{11}$; the remaining matrix 
elements are estimated analogously. {}From Equations~(130) and~(131), and the 
duplication formula for the gamma function \cite{a24}, $\Gamma (2x) \! = \! 
\tfrac{2^{2x-1}}{\sqrt{\pi}} \Gamma (x) \Gamma (x \! + \! \tfrac{1}{2})$, one 
shows that $g_{11} \! =_{\underset{\epsilon_{2}< \lambda < \epsilon_{1}}{\tau 
\to +0}} \! \left((g_{11})_{0} \! + \! (g_{11})_{\Delta} \! + \! \Delta 
\mathscr{E}_{11} \right) \! \left(1 \! + \! o(\tau^{\delta_{2}}) \right)$, 
$\delta_{2} \! > \! 0$, where
\begin{align*}
(g_{11})_{0} :=& \, \sum\limits_{l = \pm 1} \! \dfrac{2^{l \widehat{\rho}} 
\me^{\frac{\mi \pi l \widehat{\rho}}{2}}(\varepsilon b)^{-l \widehat{\rho}}
(\Gamma (2l \widehat{\rho}))^{2} \tau^{-2l \widehat{\rho}} \me^{-\frac{\pi 
a}{4}}(2 \tau)^{\frac{a \mi}{2}} \me^{-\frac{\mi \pi}{4}}}{4 \sqrt{\smash[b]{
\pi \varepsilon b}} \, \Gamma \! \left(l \widehat{\rho} \! + \! \frac{a \mi}{
2} \right)} \\
\times& \, \left(\! \left(-4l \widehat{\rho} \! + \! 2a \mi \! + \! \dfrac{4B 
\, \widehat{\gamma}}{\sqrt{\smash[b]{-AB}}} \right) \! \left(\! -\sqrt{
\smash[b]{B}} \! + \! \dfrac{\widehat{\delta} \sqrt{\smash[b]{-AB}}}{\sqrt{
\smash[b]{B}} \left(l \widehat{\rho} \! + \! \frac{a \mi}{2} \right)} \right) 
\! - \! \dfrac{2 \mi \varepsilon b \sqrt{\smash[b]{B}}}{\sqrt{\smash[b]{-AB}}} 
\right), \\
(g_{11})_{\Delta} :=& \, \sum\limits_{l = \pm 1} \! \dfrac{2^{3l \widehat{
\rho}} \me^{\frac{3 \pi \mi l \widehat{\rho}}{2}}(\varepsilon b)^{-l \widehat{
\rho}} \lambda^{2l \widehat{\rho}} \Gamma (l \widehat{\rho}) \Gamma (-l 
\widehat{\rho}) \me^{-\frac{\pi a}{4}}(2 \tau)^{\frac{a \mi}{2}} \me^{-\frac{
\mi \pi}{4}}}{4 \sqrt{\smash[b]{\pi \varepsilon b}} \, \Gamma \! \left(-l 
\widehat{\rho} \! + \! \frac{a \mi}{2} \right)} \\
\times& \, \left(\! \left(-4l \widehat{\rho} \! + \! 2a \mi \! + \! \dfrac{4B 
\, \widehat{\gamma}}{\sqrt{\smash[b]{-AB}}} \right) \! \left(\! -\sqrt{
\smash[b]{B}} \! + \! \dfrac{\widehat{\delta} \sqrt{\smash[b]{-AB}}}{\sqrt{
\smash[b]{B}} \left(-l \widehat{\rho} \! + \! \frac{a \mi}{2} \right)} \right) 
\! - \! \dfrac{2 \mi \varepsilon b \sqrt{\smash[b]{B}}}{\sqrt{\smash[b]{-AB}}} 
\right), \\
\Delta \mathscr{E}_{11} :=& \, \, \tau^{2+ \frac{a \mi}{2}} \sum\limits_{l_{1} 
= \tau,\lambda} \sum\limits_{l_{2} = \pm 1}l_{1}^{2l_{2} \widehat{\rho}} \! 
\left(a_{0}^{l_{1}}(l_{2}) \sqrt{\smash[b]{B}}+ \! \dfrac{a_{1}^{l_{1}}(l_{2}) 
\widehat{\delta} \sqrt{\smash[b]{-AB}}}{\sqrt{\smash[b]{B}}} \! + \! \dfrac{
a_{2}^{l_{1}}(l_{2}) B \sqrt{\smash[b]{B}} \, \widehat{\gamma}}{\sqrt{
\smash[b]{-AB}}} \right. \\
+&\left. \dfrac{a_{3}^{l_{1}}(l_{2}) \sqrt{\smash[b]{B}}}{\sqrt{\smash[b]{-A
B}}} \right),
\end{align*}
with $a_{k}^{l_{1}}(l_{2})$, $k \! = \! 0,1,2,3$, $l_{1} \! = \! \tau,
\lambda$, $l_{2} \! = \! \pm 1$, some $\mathcal{O}(1)$ coefficients (only the 
structure of $\Delta \mathscr{E}_{11}$ is essential). Recalling that $\widehat{
\rho}^{2} \! := \! \widehat{\gamma} \widehat{\delta} \! - \! a^{2}/4$ 
(Equation~(111)), it follows {}from Equations~(13) and~(14) that $(g_{11})_{
\Delta} \! = \! 0$. Using the restrictions stated in the Lemma, one estimates 
$\Delta \mathscr{E}_{11}$ as $o(\tau^{\delta_{1}})$, $\delta_{1} \! > \! 0$. 
Finally, Equation~(126) is obtained after simplification of the expression for 
$(g_{11})_{0}$ upon using Equations~(13) and~(14), and the definition of 
$\widehat{\rho}^{2}$. \hfill $\square$
\begin{ccc}
Let $0 \! < \! \epsilon \! < \! 1$,
\begin{equation*}
a \! \not= \! 0, \qquad \text{and} \qquad \widehat{\rho} \! = \! 0 \, \, \, 
(\equiv \! \nu_{o} \! = \! 0).
\end{equation*}
Then, under conditions~{\rm (113)}, {\rm (122)}, {\rm (123)}, and
\begin{align*}
\max \! \left\{\! \left\vert \sqrt{-A} \right\vert,\left\vert \widehat{\gamma} 
\sqrt{\smash[b]{B}} \right\vert \right\} \! &\underset{\tau \to +0}{>} \! 
\tau^{1-\epsilon} \max \! \left\{\! \left\vert \sqrt{\smash[b]{B}} \right\vert,
\left\vert \dfrac{\widehat{\gamma}}{\sqrt{\smash[b]{-A}}} \right\vert,
\left\vert \dfrac{\widehat{\gamma}^{2}B}{\sqrt{\smash[b]{-A}}} \right\vert 
\right\}, \\
\max \! \left\{\! \left\vert \sqrt{\smash[b]{B}} \right\vert,\left\vert 
\widehat{\delta} \sqrt{\smash[b]{-A}} \right\vert \right\} \! &\underset{\tau 
\to +0}{>} \! \tau^{1-\epsilon} \max \! \left\{\! \left\vert \sqrt{-A} 
\right\vert,\dfrac{1}{\left\vert \sqrt{\smash[b]{-A}} \right\vert},\left\vert 
\dfrac{B \widehat{\gamma}}{\sqrt{\smash[b]{-A}}} \right\vert \right\},
\end{align*}
or conditions~{\rm (113)}, {\rm (122)}, {\rm (124)}, and
\begin{align*}
\max \! \left\{\! \left\vert \sqrt{\smash[b]{B}} \right\vert,\left\vert 
\widehat{\delta} \sqrt{\smash[b]{-A}} \right\vert \right\} \! &\underset{\tau 
\to +0}{>} \! \tau^{1-\epsilon} \max \! \left\{\! \left\vert \dfrac{\widehat{
\delta}}{\sqrt{\smash[b]{B}}} \right\vert,\left\vert \dfrac{\widehat{\delta}^{
2}A}{\sqrt{\smash[b]{B}}} \right\vert \right\}, \\
\max \! \left\{\! \left\vert \sqrt{\smash[b]{-A}} \right\vert,\left\vert 
\widehat{\gamma} \sqrt{\smash[b]{B}} \right\vert \right\} \! &\underset{\tau 
\to +0}{>} \! \tau^{1-\epsilon} \max \! \left\{\dfrac{1}{\left\vert \sqrt{
\smash[b]{B}} \right\vert},\left\vert \dfrac{\widehat{\delta}A}{\sqrt{
\smash[b]{B}}} \right\vert \right\},
\end{align*}
the connection matrix, $G$ (Equation~{\rm (28))}, has the following 
asymptotics $(g_{ij} \! := \! (G)_{ij}$, $i,j \! = \! 1,2)$,
\begin{align}
g_{11} \underset{\tau \to +0}{=}& \, \dfrac{(2 \tau)^{\frac{a \mi}{2}} \me^{-
\frac{\pi a}{4}} \me^{-\frac{\mi \pi}{4}}}{\sqrt{\smash[b]{\pi \varepsilon b}} 
\, \Gamma \! \left(\frac{a \mi}{2} \right)} \! \left(\! \left(\! \sqrt{
\smash[b]{B}} + \! \dfrac{2 \mi \widehat{\delta} \sqrt{\smash[b]{-AB}}}{a 
\sqrt{\smash[b]{B}}} \right) \! \left(4 \psi (1) \! - \! \psi (\tfrac{a \mi}{
2}) \! + \! \ln 2 \! + \! \dfrac{\pi \mi}{2} \! - \! \ln (\varepsilon b) 
\right. \right. \nonumber \\
-&\left. \left. \, 2 \ln \tau \right) \! - \! \dfrac{4 \widehat{\delta} \sqrt{
\smash[b]{-AB}}}{a^{2} \sqrt{\smash[b]{B}}} \right) \! \left(1 \! + \! o \! 
\left(\tau^{\delta} \right) \right), \\
g_{12} \underset{\tau \to +0}{=}& \, -\dfrac{(2 \tau)^{-\frac{a \mi}{2}} \me^{
-\frac{\pi a}{4}} \me^{-\frac{\mi \pi}{4}}}{\sqrt{\pi \varepsilon b} \, \Gamma 
\! \left(-\frac{a \mi}{2} \right)} \! \left(\! \left(\! \dfrac{\sqrt{\smash[b]{
-AB}}}{\sqrt{\smash[b]{B}}} \! - \! \dfrac{2 \mi \widehat{\gamma} \sqrt{
\smash[b]{B}}}{a} \right) \! \left(4 \psi (1) \! - \! \psi (-\tfrac{a \mi}{2}) 
\! + \! \ln 2 \! + \! \dfrac{3 \pi \mi}{2} \! - \! \ln (\varepsilon b) \right. 
\right. \nonumber \\
-&\left. \left. \, 2 \ln \tau \right) \! - \! \dfrac{4 \widehat{\gamma} \sqrt{
\smash[b]{B}}}{a^{2}} \right) \! \left(1 \! + \! o \! \left(\tau^{\delta} 
\right) \right), \\
g_{21} \underset{\tau \to +0}{=}& \, -\dfrac{(2 \tau)^{\frac{a \mi}{2}} \me^{
-\frac{\pi a}{4}} \me^{\frac{\mi \pi}{4}}}{\sqrt{\smash[b]{\pi \varepsilon b}} 
\, \Gamma \! \left(\frac{a \mi}{2} \right)} \! \left(\! \left(\! \sqrt{
\smash[b]{B}} + \! \dfrac{2 \mi \widehat{\delta} \sqrt{\smash[b]{-AB}}}{a 
\sqrt{\smash[b]{B}}} \right) \! \left(4 \psi (1) \! - \! \psi (\tfrac{a \mi}{
2}) \! + \! \ln 2 \! - \! \dfrac{3 \pi \mi}{2} \! - \! \ln (\varepsilon b) 
\right. \right. \nonumber \\
-&\left. \left. \, 2 \ln \tau \right) \! - \! \dfrac{4 \widehat{\delta} \sqrt{
\smash[b]{-AB}}}{a^{2} \sqrt{\smash[b]{B}}} \right) \! \left(1 \! + \! o \! 
\left(\tau^{\delta} \right) \right), \\
g_{22} \underset{\tau \to +0}{=}& \, \dfrac{(2 \tau)^{-\frac{a \mi}{2}} \me^{
-\frac{\pi a}{4}} \me^{\frac{\mi \pi}{4}}}{\sqrt{\smash[b]{\pi \varepsilon b}} 
\, \Gamma \! \left(-\frac{a \mi}{2} \right)} \! \left(\! \left(\! \dfrac{
\sqrt{\smash[b]{-AB}}}{\sqrt{\smash[b]{B}}} \! - \! \dfrac{2 \mi \widehat{
\gamma} \sqrt{\smash[b]{B}}}{a} \right) \! \left(4 \psi (1) \! - \! \psi 
(-\tfrac{a \mi}{2}) \! + \! \ln 2 \! - \! \dfrac{\pi \mi}{2} \! - \! \ln 
(\varepsilon b) \right. \right. \nonumber \\
-&\left. \left. \, 2 \ln \tau \right) \! - \! \dfrac{4 \widehat{\gamma} \sqrt{
\smash[b]{B}}}{a^{2}} \right) \! \left(1 \! + \! o \! \left(\tau^{\delta} 
\right) \right),
\end{align}
where $\psi (z) \! := \! \frac{\md}{\md z} \ln \Gamma (z)$ is the psi 
function, $\psi (1) \! = \! -0.57721566490 \ldots$ {\rm \cite{a24}}, and 
$\delta \! > \! 0$.
\end{ccc}

\emph{Proof.} One proceeds as in the proof of Lemma~5.1; but, since $\widehat{
\rho} \! = \! 0$, the following representation \cite{a27} for the Whittaker 
function (instead of Equation~(132)) is used,
\begin{align*}
W_{z_{1},z_{2}}&(z) \! = \! \dfrac{(-1)^{2z_{2}}z^{z_{2}+\frac{1}{2}} \me^{-
\frac{z}{2}}}{\Gamma \! \left(\frac{1}{2} \! - \! z_{2} \! - \! z_{1} \right) 
\Gamma \! \left(\frac{1}{2} \! + \! z_{2} \! - \! z_{1} \right)} \! \left(
\sum\limits_{k=0}^{\infty} \dfrac{\Gamma (z_{2} \! + \! k \! - \! z_{1} \! + 
\! \frac{1}{2})}{k!(2z_{2} \! + \! k)!} \! \left(\psi (k \! + \! 1) \! + \! 
\psi (2z_{2} \! + \! k \! + \! 1) \right. \right. \\
-&\left. \left. \, \psi (z_{2} \! + \! k \! - \! z_{1} \! + \! \tfrac{1}{2}) 
\! - \! \ln z \right) \! z^{k} \! + \! (-z)^{-2z_{2}} \sum\limits_{k=0}^{2z_{
2}-1} \dfrac{\Gamma (2z_{2} \! - \! k) \Gamma (k \! - \! z_{2} \! - \! z_{1} 
\! + \! \frac{1}{2})}{k!}(-z)^{k} \right),
\end{align*}
$\vert \arg z \vert \! < \! 3 \pi/2$, $2z_{2} \! + \! 1 \! \in \! \mathbb{N}$, 
and, when $z_{2} \! = \! 0$, the second sum in the above expansion is equal to 
zero. The following identities for the Whittaker and psi functions were used 
\cite{a24}: $W_{0,-1}(\cdot) \! = \! W_{0,1}(\cdot)$, $\psi (z \! + \! 1) \! = 
\! \psi (z) \! + \! \tfrac{1}{z}$, and $\psi (\tfrac{1}{2}) \! = \! \psi (1) 
\! - \! 2 \ln 2$. \hfill $\square$
\begin{proposition}
Let $G$ be the connection matrix of Equation~{\rm (55)} with 
$\tau$-independent elements satisfying conditions~{\rm (43)}, and $\mathfrak{p}
(z_{1},z_{2})$, $\chi_{1}(\vec{g};z_{3}) \! := \! \chi_{1}(\vec{g}(0,0);
z_{3})$ and $\chi_{2}(\vec{g};z_{4}) \! := \! \chi_{2}(\vec{g}(0,0);z_{4})$, 
where $g_{ij}(0,0) \! = \! g_{ij}$, $i,j \! = \! 1,2$, be defined as in 
Equations~{\rm (46)} and~{\rm (47)}. Then the corresponding isomonodromy 
deformations have the following asymptotic representation,
\begin{align}
\dfrac{\sqrt{\smash[b]{-AB}}}{\sqrt{\smash[b]{B}}} \underset{\tau \to +0}{=}& 
\, \dfrac{\mi \me^{\frac{\pi a}{4}}(2 \tau)^{\frac{a \mi}{2}}}{4 \sqrt{
\smash[b]{\pi}} \, (\varepsilon b)^{-1/2}} \! \left(\mathfrak{p}(-a,\widehat{
\rho}) \me^{-\mi \pi \widehat{\rho}} \chi_{2}(\vec{g};\widehat{\rho}) \tau^{2 
\widehat{\rho}} \! + \! \mathfrak{p}(-a,-\widehat{\rho}) \me^{\mi \pi \widehat{
\rho}} \chi_{2}(\vec{g};-\widehat{\rho}) \tau^{-2 \widehat{\rho}} \right)& 
\nonumber \\
\times& \, \left(1 \! + \! o \! \left(\tau^{\delta} \right) \right), \\
\sqrt{\smash[b]{B}} \, \widehat{\gamma} \underset{\tau \to +0}{=}& \, \dfrac{
\mi \me^{\frac{\pi a}{4}}(2 \tau)^{\frac{a \mi}{2}}}{4 \sqrt{\smash[b]{\pi}} 
\, (\varepsilon b)^{-1/2}} \! \left(\! \left(-\widehat{\rho} \! - \! \dfrac{a 
\mi}{2} \right) \! \mathfrak{p}(-a,\widehat{\rho}) \me^{-\mi \pi \widehat{
\rho}} \chi_{2}(\vec{g};\widehat{\rho}) \tau^{2 \widehat{\rho}} \! + \! \left(
\widehat{\rho} \! - \! \frac{a \mi}{2} \right) \! \mathfrak{p}(-a,-\widehat{
\rho}) \me^{\mi \pi \widehat{\rho}} \right.& \nonumber \\
\times&\left. \, \chi_{2}(\vec{g};-\widehat{\rho}) \tau^{-2 \widehat{\rho}} 
\right) \! \left(1 \! + \! o \! \left(\tau^{\delta} \right) \right), \\
\sqrt{\smash[b]{B}} \underset{\tau \to +0}{=}& \, -\dfrac{\mi \me^{\frac{\pi 
a}{4}}(2 \tau)^{-\frac{a \mi}{2}}}{4 \sqrt{\smash[b]{\pi}} \, (\varepsilon 
b)^{-1/2}} \! \left(\mathfrak{p}(a,\widehat{\rho}) \chi_{1}(\vec{g};\widehat{
\rho}) \tau^{2 \widehat{\rho}} \! + \! \mathfrak{p}(a,-\widehat{\rho}) \chi_{1}
(\vec{g};-\widehat{\rho}) \tau^{-2 \widehat{\rho}} \right) \nonumber \\
\times& \, \left(1 \! + \! o \! \left(\tau^{\delta} \right) \right), \\
\dfrac{\widehat{\delta} \sqrt{\smash[b]{-AB}}}{\sqrt{\smash[b]{B}}} \underset{
\tau \to +0}{=}& \, -\frac{\mi \me^{\frac{\pi a}{4}}(2 \tau)^{-\frac{a \mi}{2}
}}{4 \sqrt{\smash[b]{\pi}} \, (\varepsilon b)^{-1/2}} \! \left(\! \left(-
\widehat{\rho} \! + \! \dfrac{a \mi}{2} \right) \! \mathfrak{p}(a,\widehat{
\rho}) \chi_{1}(\vec{g};\widehat{\rho}) \tau^{2 \widehat{\rho}} \! + \! \left(
\widehat{\rho} \! + \! \dfrac{a \mi}{2} \right) \! \mathfrak{p}(a,-\widehat{
\rho}) \chi_{1}(\vec{g};-\widehat{\rho}) \right. \nonumber \\
\times&\left. \, \tau^{-2 \widehat{\rho}} \right) \! \left(1 \! + \! o \! 
\left(\tau^{\delta} \right) \right), \\
\sqrt{\smash[b]{-AB}} \underset{\tau \to +0}{=}& \, \dfrac{\varepsilon b \me^{
\frac{\pi a}{2}}}{16 \pi} \! \left(\mathfrak{p}(a,\widehat{\rho}) \chi_{1}
(\vec{g};\widehat{\rho}) \tau^{2 \widehat{\rho}} \! + \! \mathfrak{p}(a,-
\widehat{\rho}) \chi_{1}(\vec{g};-\widehat{\rho}) \tau^{-2 \widehat{\rho}} 
\right) \nonumber \\
\times& \, \left(\mathfrak{p}(-a,\widehat{\rho}) \me^{-\mi \pi \widehat{\rho}} 
\chi_{2}(\vec{g};\widehat{\rho}) \tau^{2 \widehat{\rho}} \! + \! \mathfrak{p}
(-a,-\widehat{\rho}) \me^{\mi \pi \widehat{\rho}} \chi_{2}(\vec{g};-\widehat{
\rho}) \tau^{-2 \widehat{\rho}} \right) \! \left(1 \! + \! o \! \left(\tau^{
\delta} \right) \right), \\
\widehat{\gamma} \! := \! \tau C \underset{\tau \to +0}{=}& \, -\dfrac{\left((
-\widehat{\rho} \! - \! \frac{a \mi}{2}) \mathfrak{p}(-a,\widehat{\rho}) \me^{
-\mi \pi \widehat{\rho}} \chi_{2}(\vec{g};\widehat{\rho}) \tau^{2 \widehat{
\rho}} \! + \! (\widehat{\rho} \! - \! \frac{a \mi}{2}) \mathfrak{p}(-a,-
\widehat{\rho}) \me^{\mi \pi \widehat{\rho}} \chi_{2}(\vec{g};-\widehat{\rho}) 
\tau^{-2 \widehat{\rho}} \right)}{(2 \tau)^{-a \mi} \left(\mathfrak{p}(a,
\widehat{\rho}) \chi_{1}(\vec{g};\widehat{\rho}) \tau^{2 \widehat{\rho}} \! + 
\! \mathfrak{p}(a,-\widehat{\rho}) \chi_{1}(\vec{g};-\widehat{\rho}) \tau^{-2 
\widehat{\rho}} \right)} \nonumber \\
\times& \, \left(1 \! + \! o \! \left(\tau^{\delta} \right) \right), \\
\widehat{\delta} \! := \! \tau D \underset{\tau \to +0}{=}& \, -\dfrac{\left((-
\widehat{\rho} \! + \! \frac{a \mi}{2}) \mathfrak{p}(a,\widehat{\rho}) \chi_{1}
(\vec{g};\widehat{\rho}) \tau^{2 \widehat{\rho}} \! + \! (\widehat{\rho} \! + 
\! \frac{a \mi}{2}) \mathfrak{p}(a,-\widehat{\rho}) \chi_{1}(\vec{g};-\widehat{
\rho}) \tau^{-2 \widehat{\rho}} \right)}{(2 \tau)^{a \mi} \left(\mathfrak{p}(-
a,\widehat{\rho}) \me^{-\mi \pi \widehat{\rho}} \chi_{2}(\vec{g};\widehat{\rho}
) \tau^{2 \widehat{\rho}} \! + \! \mathfrak{p}(-a,-\widehat{\rho}) \me^{\mi 
\pi \widehat{\rho}} \chi_{2}(\vec{g};-\widehat{\rho}) \tau^{-2 \widehat{\rho}} 
\right)} \nonumber \\
\times& \, \left(1 \! + \! o \! \left(\tau^{\delta} \right) \right),
\end{align}
with $\delta \! > \! 0$.
\end{proposition}

\emph{Proof.} Algebraically inverting Equations~(126)--(129) and using the 
well-known identities for the gamma function \cite{a24} $\Gamma (2z) \! = \! 
\tfrac{2^{2z-1}}{\sqrt{\pi}} \Gamma (z) \Gamma (z \! + \! \tfrac{1}{2})$ and 
$\Gamma (\tfrac{1}{2} \! + \! z) \Gamma (\tfrac{1}{2} \! - \! z) \! = \! 
\tfrac{\pi}{\cos (\pi z)}$, one obtains Equations~(137)--(143). It is easy to 
verify that the system of Equations~(137)--(143) is compatible. Assuming that 
$g_{ij}$, $i,j \! = \! 1,2$, are constant, $\widehat{\rho} \! \not= \! 0$, and 
$\vert \Re (\widehat{\rho}) \vert \! < \! 1/2$, one proves that any functions 
$A$, $B$, $C$, and $D$ with asymptotics~(137)--(143) satisfy 
restrictions~(113), (122), (123), (124), and those given in Lemma~5.1; 
therefore, one is now in a position to use the justification scheme suggested 
in \cite{a20} to complete the proof. \hfill $\square$
\begin{proposition}
Under the conditions of Proposition~{\rm 5.5},
\begin{equation}
\cos (2 \pi \widehat{\rho}) \underset{\tau \to +0}{=} \cos (2 \pi \rho) \! 
\left(1 \! + \! o \! \left(\tau^{\delta} \right) \right),
\end{equation}
where
$\cos (2 \pi \rho)$ is defined in Equation~{\rm (44)}, and $\delta \! > \! 0$.
\end{proposition}

\emph{Proof.} One recalls {}from Equations~(111) that $\widehat{\rho}^{2} \! 
+ \! a^{2}/4 \! = \! \widehat{\gamma} \, \widehat{\delta}$. Substituting the 
asymptotic expressions for $\widehat{\gamma}$ and $\widehat{\delta}$ given in 
Equations~(142) and~(143) into this relation and simplifying, one shows that
\begin{align*}
&\left(\widehat{\rho} \! - \! \dfrac{a \mi}{2} \right) \! \mathfrak{p}(a,
\widehat{\rho}) \mathfrak{p}(-a,-\widehat{\rho}) \me^{\mi \pi \widehat{\rho}} 
\chi_{1}(\vec{g};\widehat{\rho}) \chi_{2}(\vec{g};-\widehat{\rho}) \\
+& \left(\widehat{\rho} \! + \! \dfrac{a \mi}{2} \right) \! \mathfrak{p}(a,-
\widehat{\rho}) \mathfrak{p}(-a,\widehat{\rho}) \me^{-\mi \pi \widehat{\rho}} 
\chi_{1}(\vec{g};-\widehat{\rho}) \chi_{2}(\vec{g};\widehat{\rho}) \! 
\underset{\tau \to +0}{=} \! o \! \left(\tau^{\delta} \right), \quad \delta 
\! > \! 0,
\end{align*}
with $\mathfrak{p}(z_{1},z_{2})$, $\chi_{1}(\vec{g};z_{3}) \! := \! \chi_{1}
(\vec{g}(0,0);z_{3})$ and $\chi_{2}(\vec{g};z_{4}) \! := \! \chi_{2}(\vec{g}
(0,0);z_{4})$, where $g_{ij}(0,0) \! = \! g_{ij}$, $i,j \! = \! 1,2$, defined 
in Equations~(46) and~(47). Applying to the last formula the gamma function 
identity $\Gamma (z) \Gamma (1 \! - \! z) \! = \! \tfrac{\pi}{\sin \pi z}$, 
one arrives at
\begin{align*}
&\dfrac{\me^{\mi \pi \widehat{\rho}}(\mi g_{11}g_{12} \! - \! \mi g_{21}g_{22} 
\! + \! g_{11}g_{22} \me^{2 \pi \mi \widehat{\rho}} \! + \! g_{12}g_{21} \me^{
-2 \pi \mi \widehat{\rho}})}{\sin \! \left(\pi \! \left(\widehat{\rho} \! + \! 
\frac{a \mi}{2} \right) \right)} \\
+& \, \dfrac{\me^{-\mi \pi \widehat{\rho}}(\mi g_{11}g_{12} \! - \! \mi g_{21}
g_{22} \! + \! g_{11}g_{22} \me^{-2 \pi \mi \widehat{\rho}} \! + \! g_{12}g_{
21} \me^{2\pi \mi \widehat{\rho}})}{\sin \! \left(\pi \! \left(\widehat{\rho} 
\! - \! \frac{a \mi}{2} \right) \right)} \! \underset{\tau \to +0}{=} \! o \! 
\left(\tau^{\delta} \right).
\end{align*}
Now, using trigonometric identities and the fact that $g_{11}g_{22} \! - \! 
g_{12}g_{21} \! = \! 1$, one obtains, in the case $g_{11}g_{22} \! \not= \! 
0$, Equation~(144), where $\cos (2 \pi \rho)$ is defined in Equation~(44): the 
second relation of Equation~(44) is a direct consequence of Equation~(33). 
\hfill $\square$
\begin{proposition}
Let $G$ be the connection matrix of Equation~{\rm (55)} with 
$\tau$-independent elements satisfying conditions~{\rm (49)}, and $Q(z)$, 
$\chi_{1}(\vec{g};z_{3}) \! := \! \chi_{1}(\vec{g}(0,0);z_{3})$ and $\chi_{2}
(\vec{g};z_{4}) \! := \! \chi_{2}(\vec{g}(0,0);z_{4})$, where $g_{ij}(0,0) 
\! = \! g_{ij}$, $i,j \! = \! 1,2$, be as defined in Equations~{\rm (51)} 
and~{\rm (47)}. Then the corresponding isomonodromy deformations have the 
following asymptotic representation,
\begin{align}
\dfrac{\sqrt{\smash[b]{-AB}}}{\sqrt{\smash[b]{B}}} \underset{\tau \to +0}{=}& 
\, \dfrac{\mi \sqrt{\varepsilon b} \, \Gamma \! \left(-\frac{a \mi}{2} \right) 
\! \me^{\frac{\pi a}{4}}}{2 \sqrt{\smash[b]{\pi}} \, (2 \tau)^{-\frac{a \mi}{
2}}} \! \left(\chi_{2}(\vec{g};0) \! \left(1 \! + \! \dfrac{a \mi}{2}Q(-a) 
\right) \! + \! \dfrac{\pi a}{4} \! \left(g_{12} \me^{\frac{\mi \pi}{4}} \! - 
\! 3g_{22} \me^{-\frac{\mi \pi}{4}} \right) \right. \nonumber \\
-&\left. \, a \mi \chi_{2}(\vec{g};0) \ln \tau \right) \! \left(1 \! + \! o \! 
\left(\tau^{\delta} \right) \right), \\
\sqrt{\smash[b]{B}} \, \widehat{\gamma} \underset{\tau \to +0}{=}& \, \dfrac{
\mi \sqrt{\smash[b]{\varepsilon b}} \, \Gamma \! \left(-\frac{a \mi}{2} 
\right) \! a^{2} \me^{\frac{\pi a}{4}}}{8 \sqrt{\smash[b]{\pi}} \, (2 \tau)^{
-\frac{a \mi}{2}}} \! \left( \chi_{2}(\vec{g};0)Q(-a) \! + \! \dfrac{\mi \pi}{
2} \! \left(3g_{22} \me^{-\frac{\mi \pi}{4}} \! - \! g_{12} \me^{\frac{\mi 
\pi}{4}} \right) \right. \nonumber \\
-&\left. \, 2 \chi_{2}(\vec{g};0) \ln \tau \right) \! \left(1 \! + \! o \! 
\left(\tau^{\delta} \right) \right), \\
\sqrt{\smash[b]{B}} \underset{\tau \to +0}{=}& \, -\dfrac{\mi \sqrt{\smash[b]{
\varepsilon b}} \, \Gamma \! \left(\frac{a \mi}{2} \right) \! \me^{\frac{\pi 
a}{4}}}{2 \sqrt{\smash[b]{\pi}} \, (2 \tau)^{\frac{a \mi}{2}}} \! \left(\chi_{
1}(\vec{g};0) \! \left(1 \! - \! \dfrac{a \mi}{2}Q(a) \right) \! + \! \dfrac{
\pi a}{4} \! \left(g_{21} \me^{-\frac{\mi \pi}{4}} \! - \! 3g_{11} \me^{\frac{
\mi \pi}{4}} \right) \right. \nonumber \\
+&\left. \, a \mi \chi_{1}(\vec{g};0) \ln \tau \right) \! \left(1 \! + \! o \! 
\left(\tau^{\delta} \right) \right), \\
\dfrac{\widehat{\delta} \sqrt{\smash[b]{-AB}}}{\sqrt{\smash[b]{B}}} \underset{
\tau \to +0}{=}& \, -\dfrac{\mi \sqrt{\smash[b]{\varepsilon b}} \, \Gamma 
\! \left(\frac{a \mi}{2} \right) \! a^{2} \me^{\frac{\pi a}{4}}}{8 \sqrt{
\smash[b]{\pi}} \, (2 \tau)^{\frac{a \mi}{2}}} \! \left(\chi_{1}(\vec{g};0)Q
(a) \! + \! \dfrac{\mi \pi}{2} \! \left(g_{21} \me^{-\frac{\mi \pi}{4}} \! - 
\! 3g_{11} \me^{\frac{\mi \pi}{4}} \right) \right. \nonumber \\
-&\left. \, 2 \chi_{1}(\vec{g};0) \ln \tau \right) \! \left(1 \! + \! o \! 
\left(\tau^{\delta} \right) \right), \\
\sqrt{\smash[b]{-AB}} \underset{\tau \to +0}{=}& \, \dfrac{\varepsilon b 
\me^{\frac{\pi a}{2}}}{2 a \sinh \! \left(\frac{\pi a}{2} \right)} \! \left(
\chi_{1}(\vec{g};0) \! \left(1 \! - \! \dfrac{a \mi}{2}Q(a) \right) \! + \! 
\dfrac{\pi a}{4} \! \left(g_{21} \me^{-\frac{\mi \pi}{4}} \! - \! 3g_{11} 
\me^{\frac{\mi \pi}{4}} \right) \! + \! a \mi \chi_{1}(\vec{g};0) \right. 
\nonumber \\
\times&\left. \, \ln \tau \right) \! \left(\chi_{2}(\vec{g};0) \! \left(1 \! + 
\! \dfrac{a \mi}{2}Q(-a) \right) \! + \! \dfrac{\pi a}{4} \! \left(g_{12} \me^{
\frac{\mi \pi}{4}} \! - \! 3g_{22} \me^{-\frac{\mi \pi}{4}} \right) \! - \! a 
\mi \chi_{2}(\vec{g};0) \ln \tau \right) \nonumber \\
\times& \, \left(1 \! + \! o \! \left(\tau^{\delta} \right) \right), \\
\widehat{\gamma} \! := \! \tau C \underset{\tau \to +0}{=}& \, -\tfrac{\pi 
a (2 \tau)^{a \mi}}{2 \left(\Gamma \left(\frac{a \mi}{2} \right) \right)^{
2} \sinh \left(\frac{\pi a}{2} \right)} \tfrac{\left(\chi_{2}(\vec{g};0)Q(-a)
+\frac{\mi \pi}{2}(3g_{22} \me^{-\frac{\mi \pi}{4}}-g_{12} \me^{\frac{\mi 
\pi}{4}})-2 \chi_{2}(\vec{g};0) \ln \tau \right)}{\left(\chi_{1}(\vec{g};0)
(1-\frac{a \mi}{2}Q(a))+\frac{\pi a}{4}(g_{21} \me^{-\frac{\mi \pi}{4}}-3
g_{11} \me^{\frac{\mi \pi}{4}})+a \mi \chi_{1}(\vec{g};0) \ln \tau \right)} 
\nonumber \\
\times& \, \left(1 \! + \! o \! \left(\tau^{\delta} \right) \right), \\
\widehat{\delta} \! := \! \tau D \underset{\tau \to +0}{=}& \, -\tfrac{a^{3} 
\left(\Gamma \left(\frac{a \mi}{2} \right) \right)^{2} \sinh \left(\frac{\pi 
a}{2} \right)}{8 \pi (2 \tau)^{a \mi}} \tfrac{\left(\chi_{1}(\vec{g};0)Q(a)+
\frac{\mi \pi}{2}(g_{21} \me^{-\frac{\mi \pi}{4}}-3g_{11} \me^{\frac{\mi \pi}{
4}})-2 \chi_{1}(\vec{g};0) \ln \tau \right)}{\left(\chi_{2}(\vec{g};0)(1+
\frac{a \mi}{2}Q(-a))+\frac{\pi a}{4}(g_{12} \me^{\frac{\mi \pi}{4}}-3g_{22} 
\me^{-\frac{\mi \pi}{4}})-a \mi \chi_{2}(\vec{g};0) \ln \tau \right)} 
\nonumber \\
\times& \, \left(1 \! + \! o \! \left(\tau^{\delta} \right) \right),
\end{align}
with $\delta \! > \! 0$.
\end{proposition}

\emph{Proof.} Follows {}from Lemma~5.2 in a manner analogous to the proof of 
Proposition~5.5. \hfill $\square$

Now, one can complete the proof (as $\tau \! \to \! +\infty)$ of Theorems~3.4 
and~3.5. The asymptotics for $u(\tau)$ given in Theorem~3.4 (resp., 
Theorem~3.5), Equation~(45) (resp., Equation~(50)) follows {}from 
Proposition~1.2 $(u(\tau) \! = \! \varepsilon \tau \sqrt{\smash[b]{-AB}}$, 
$\varepsilon \! = \! \pm 1)$ and Proposition~5.5, Equation~(141) (resp., 
Proposition~5.7, Equation~(149)). To get the asymptotics for $\mathcal{H}
(\tau)$ given in these Theorems, one uses the second relation of 
Equation~(37), where $u^{\prime}(\tau)$ is calculated via
\begin{equation*}
u^{\prime}(\tau) \! = \! \dfrac{u(\tau)}{\tau} \! + \! 2 \varepsilon \tau (AD 
\! - \! BC),
\end{equation*}
and: (1) in the case $\widehat{\rho} \! \not= \! 0$, Equation~(44) and the 
following identities,
\begin{gather*}
\mathfrak{p}(a,\rho) \mathfrak{p}(-a,-\rho) \me^{\mi \pi \rho} \! + \! \dfrac{
\pi \! \left(\rho \! + \! \frac{a \mi}{2} \right) \! \me^{\mi \pi \rho}}{\rho^{
2} \sin \! \left(\pi \! \left(\rho \! + \! \frac{a \mi}{2} \right) \right)} \! 
= \! 0, \\
\dfrac{\me^{\mi \pi \rho} \chi_{1}(\vec{g};\rho) \chi_{2}(\vec{g};-
\rho)}{\sin \! \left(\pi \! \left(\rho \! + \! \frac{a \mi}{2} \right) \right)
} \! + \! \dfrac{\me^{-\mi \pi \rho} \chi_{1}(\vec{g};-\rho) \chi_{2}(\vec{g};
\rho)}{\sin \! \left(\pi \! \left(\rho \! - \! \frac{a \mi}{2} \right) \right)
} \! = \! 0, \\
\dfrac{\me^{\mi \pi \rho} \chi_{1}(\vec{g};\rho) \chi_{2}(\vec{g};-\rho)}{\sin 
\! \left(\pi \! \left(\rho \! + \! \frac{a \mi}{2} \right) \right)} \! - \! 
\dfrac{\me^{-\mi \pi \rho} \chi_{1}(\vec{g};-\rho) \chi_{2}(\vec{g};\rho)}{
\sin \! \left(\pi \! \left(\rho \! - \! \frac{a \mi}{2} \right) \right)} \! 
= \! 4 \mi \me^{-\frac{\pi a}{2}};
\end{gather*}
and (2) in the case $\widehat{\rho} \! = \! 0$, Equation~(49) and
\begin{equation*}
\psi (z \! + \! 1) \! = \! \psi (z) \! + \! \dfrac{1}{z}, \qquad \psi (1 \! 
- \! z) \! = \! \psi (z) \! + \! \pi \cot (\pi z).
\end{equation*}
\section{Transformations}
\label{sec:6}
In this section we collect B\"{a}cklund transformations and Lie-point 
symmetries for Equation~(\ref{eq:dp3}) and System~(\ref{sys:ABCD}). To apply 
these transformations for 
connection results, we consider their actions on solutions of 
Systems~(\ref{eq:UV}) and~(12) $(\Phi$ and $\Psi)$, and the manifold of 
monodromy data $(\mathscr{M})$.
\subsection{B\"{a}cklund Transformations}
 \label{subsec:BT}
The B\"{a}cklund transformations for System~(\ref{sys:ABCD}), as well as for 
Equation~(\ref{eq:dp3}), 
can be constructed in a standard way via the Schlesinger transformations for 
the first equation of System~(\ref{eq:UV}) \cite{a10,a11,a12}. For the 
degenerate Painlev\'{e} V equation (equivalent to the ``complete'' third 
Painlev\'{e} equation), these transformations (with a parameter 
$\theta \! = \! \pm 1)$ were constructed in \cite{a16}; however, as noted 
in \cite{a16}, these formulae are also applicable, without modification 
(simply set $\theta \! = \! 0)$, to the degenerate Painlev\'{e} III equation. 
For the reader's convenience, the latter formulae are given below. The set of 
Schlesinger transformations forms a group which acts covariantly on the set 
of solutions of System~(\ref{eq:UV}):
\begin{equation*}
\Phi_{1}(\lambda,\tau) \! = \! \mathscr{R} \Phi (\lambda,\tau),
\end{equation*}
where $\Phi_{1}(\lambda,\tau)$ is the solution of System~(\ref{eq:UV}) with 
some functions $A_{1}(\tau)$, $B_{1}(\tau)$, $C_{1}(\tau)$, and $D_{1}(\tau)$, 
respectively, in place of $A(\tau)$, $B(\tau)$, $C(\tau)$, and $D(\tau)$ 
for $\Phi (\lambda,\tau)$. For System~(\ref{eq:UV}), the group of Schlesinger 
transformations is a free cyclic group with generator $\mathscr{R}_{1,
2}$\footnote{Note that there are two misprints in \cite{a16}: (i) the factor 
$DA$ which appears in the $(2 \, 2)$-element of the right-most matrix of 
Equation~(6.2) should be changed to $-DA$; and (ii) the function $\mathscr{W}$ 
which appears in Equation~(6.4) should be changed to $\mathscr{W} \! = \! 
\tfrac{1+a \mi}{\tau^{2}} \! - \! \tfrac{B}{\tau D}$.}:
\begin{gather*}
\mathscr{R}_{1,2} \! = \! 
\begin{pmatrix}
0 & 0 \\
0 & 1
\end{pmatrix} \! \left(\dfrac{\lambda}{\tau} \right)^{1/2} \! + \! 
\begin{pmatrix}
1 & \frac{A}{\sqrt{\smash[b]{-AB}}} \\
-\frac{D}{2 \mi \tau} & -\frac{AD}{2 \mi \tau \sqrt{\smash[b]{-AB}}}
\end{pmatrix} \! \left(\dfrac{\tau}{\lambda} \right)^{1/2}, \\
A_{1} \! = \! \dfrac{2 \mi \tau}{D} \sqrt{\smash[b]{-A_{1}B_{1}}}, \qquad 
B_{1} \! = \! -\dfrac{D}{2 \mi \tau} \sqrt{\smash[b]{-A_{1}B_{1}}}, \qquad 
\sqrt{\smash[b]{-A_{1}B_{1}}} \! = \! -\frac{\mi \varepsilon b D}{2 \tau B}, 
\qquad C_{1} \! = \! -\frac{2 \mi \tau A}{\sqrt{\smash[b]{-AB}}}, \\
D_{1} \! = \! \frac{B}{2 \mi \tau} \! - \! \frac{D}{2 \mi \tau^{2}} \! \left(
1 \! + \! a \mi \! + \! \frac{\tau AD}{\sqrt{\smash[b]{-AB}}} \right), \qquad 
\sqrt{\smash[b]{-A_{1}B_{1}}}= \! \sqrt{\smash[b]{-AB}}+ \! \frac{1}{2 \tau} 
\dfrac{\md}{\md \tau} \! \left(\frac{\tau AD}{\sqrt{\smash[b]{-AB}}} \right), 
\\
\mi a_{1} \! = \! \mi a \! + \! 1, \qquad u_{1}(\tau) \! = \! -\frac{\mi 
\varepsilon b}{8u^{2}(\tau)} \! \left(\tau (-u^{\prime}(\tau) \! + \! \mi b) 
\! + \! (2a \mi\!+\!1)u(\tau) \right).
\end{gather*}
The inverse transformation of $\mathscr{R}_{1,2}$ is
\begin{gather*}
\mathscr{R}_{3,4} \! = \! 
\begin{pmatrix}
1 & 0 \\
0 & 0
\end{pmatrix} \! \left(\frac{\lambda}{\tau} \right)^{1/2} \! + \! 
\begin{pmatrix}
-\frac{BC}{2 \mi \tau \sqrt{\smash[b]{-AB}}} & \frac{C}{2 \mi \tau} \\
-\frac{B}{\sqrt{\smash[b]{-AB}}} & 1
\end{pmatrix} \! \left(\dfrac{\tau}{\lambda} \right)^{1/2}, \\
A_{1} \! = \! -\frac{C}{2 \mi \tau} \sqrt{\smash[b]{-A_{1}B_{1}}}, \qquad B_{
1} \! = \! \frac{2 \mi \tau}{C} \sqrt{\smash[b]{-A_{1}B_{1}}}, \qquad C_{1} \! 
= \! \frac{A}{2 \mi \tau} \! + \! \frac{C}{2 \mi \tau^{2}} \! \left(1 \! - \! 
a \mi \! - \! \frac{\tau BC}{\sqrt{\smash[b]{-AB}}} \right), \\
D_{1} \! = \! -\dfrac{2 \mi \tau B}{\sqrt{\smash[b]{-AB}}}, \qquad \sqrt{
\smash[b]{-A_{1}B_{1}}}= \! \sqrt{\smash[b]{-AB}}- \! \dfrac{1}{2 \tau} \dfrac{
\md}{\md \tau} \! \left(\frac{\tau BC}{\sqrt{\smash[b]{-AB}}} \right), \\
\mi a_{1} \! = \! \mi a \! - \! 1, \qquad u_{1}(\tau) \! = \! -\frac{\mi 
\varepsilon b}{8u^{2}(\tau)} \! \left(\tau (u^{\prime}(\tau) \! + \! \mi b) \! 
+ \! (2a \mi\!-\!1)u(\tau) \right).
\end{gather*}
Denoting $v_{n}(\tau) \! = \! u_{n}(\tau)/\tau$, $n \! \in \! \mathbb{Z}$, 
where $u_{n}(\tau)$ is the general solution of Equation~(\ref{eq:dp3}) 
corresponding to the coefficient $a \! := \! a_{0} \! - \! \mi n$, one arrives 
at the following differential-difference and difference equations:
\begin{equation*}
\dfrac{\mi \varepsilon b}4v^{\prime}_{n} \! = \! v_{n}(v_{n+1} \! - \! v_{n-
1}), \quad \qquad v_{n}^2(v_{n+1} \! + \! v_{n-1}) \! = \! \dfrac{\varepsilon 
b}{4 \tau^{2}}(b \! + \! 2(a_{0} \! - \! \mi n)v_{n}).
\end{equation*}
The first one is the Kac-Moerbeke \cite{KM} equation, whilst the second should 
be equivalent to one of the so-called difference Painlev\'{e} equations. 
Consider the function $f_{n}(\tau) \! = \! p_{n}(\tau)q_{n}(\tau)/2$, where 
the functions $p_{n}(\tau)$ and $q_{n}(\tau)$ solve System~(11) for $a \! = \! 
a_{0} \! - \! \mi n$ and $\varepsilon_{1} \! = \! -1$. One finds that $f_{n}
(\tau) \! = \! \tfrac{2 \tau^{2}}{\mi \varepsilon b}v_{n+1}(\tau)v_{n}(\tau) 
\! := \! \tfrac{2 \tau^{2}}{\mi \varepsilon b}g_{n}(\tau)$. Using, now, the 
above equations for the function $v_{n}(\tau)$, one proves that
\begin{equation*}
2f_{n}(f_{n+1} \! + \! f_{n} \! + \! (\mi a_{0} \! + \! n \! + \! 1))(f_{n} \! 
+ \! f_{n-1} \! + \! (\mi a_{0} \! + \! n)) \! = \! \mi \varepsilon b \tau^{
2},
\end{equation*}
and
\begin{equation*}
\left(\dfrac{\mi \varepsilon b}{4} \right)^{2} \dfrac{\md^{2} \ln g_{n}}{\md 
\tau^{2}} \! = \! \Delta^{2}g_{n} \! := \! g_{n+1} \! + \! g_{n-1} \! - \! 2
g_{n}.
\end{equation*}
The last equation for $g_{n} \! = \! g_{n}(\tau)$ is a form of the Toda chain 
equation.

The double iteration of $\mathscr{R}_{1,2}$ (with abuse of notation, 
$\mathscr{R}_{5} \! = \! \mathscr{R}_{1,2} \circ \mathscr{R}_{1,2})$ reads
\begin{gather*}
\mathscr{R}_{5} \! = \! 
\begin{pmatrix}
0 & 0 \\
0 & 1
\end{pmatrix} \! \dfrac{\lambda}{\tau} \! + \! 
\begin{pmatrix}
0 & \frac{2 \mi \tau}{D} \\
-\frac{D}{2 \mi \tau} & \frac{\mathcal{X}}{2 \mi}
\end{pmatrix}, \qquad \quad \mathcal{X} \! := \! \frac{1 \! + \! a \mi}{\tau^{
2}} \! - \! \frac{B}{\tau D}, \\
\sqrt{\smash[b]{-A_{1}B_{1}}}= \! -\sqrt{\smash[b]{-AB}}+ \! \dfrac{\tau B 
\mathcal{X}}{D}, \qquad \sqrt{\smash[b]{-A_{1}B_{1}}}= \! \sqrt{\smash[b]{-AB}}
+ \! \dfrac{1}{2 \tau} \dfrac{\md}{\md \tau} \! \left(\dfrac{\tau B}{D} 
\right), \\
A_{1} \! = \! \dfrac{4 \tau^{2}B}{D^{2}}, \qquad B_{1} \! = \! (2 \tau)^{-2} 
\! \left(AD^{2} \! - \! \tau \mathcal{X}(\tau B \mathcal{X} \! - \! 2D \sqrt{
\smash[b]{-AB}}) \right), \\
C_{1} \! = \! \dfrac{4 \tau^{2}}{D}, \qquad D_{1} \! = \! (2 \tau)^{-2} \! 
\left(CD^{2} \! + \! 2D \sqrt{\smash[b]{-AB}}- \! (\tau B \! + \! D) \mathcal{
X} \right), \\
\mi a_{1} \! = \! \mi a \! + \! 2.
\end{gather*}
The inverse of $\mathscr{R}_{5}$ is
\begin{gather*}
\mathscr{R}_{5} \! = \! 
\begin{pmatrix}
1 & 0 \\
0 & 0
\end{pmatrix} \! \dfrac{\lambda}{\tau} \! + \! 
\begin{pmatrix}
\frac{1}{4 \mi C} \frac{\md}{\md \tau} \! \left( \frac{C}{\tau} 
\right) & \frac{C}{2 \mi \tau} \\
-\frac{2 \mi \tau}{C} & 0
\end{pmatrix}, \\
\sqrt{\smash[b]{-A_{1}B_{1}}}= \! \sqrt{\smash[b]{-AB}}- \! \dfrac{1}{2 \tau} 
\dfrac{\md}{\md \tau} \! \left(\dfrac{\tau A}{C} \right), \qquad \sqrt{
\smash[b]{-A_{1}B_{1}}}= \! -\sqrt{\smash[b]{-AB}}+ \! \dfrac{\tau A}{2C^{2}} 
\dfrac{\md}{\md \tau} \! \left(\dfrac{C}{\tau} \right), \\
A_{1} \! = \! (2 \tau)^{-2} \! \left(BC^{2} \! + \! \tau \mathcal{Y}(2C \sqrt{
\smash[b]{-AB}}+ \! \tau A \mathcal{Y}) \right), \qquad B_{1} \! = \! \dfrac{4 
\tau^{2}A}{C^{2}}, \\
C_{1} \! = \! (2 \tau)^{-2} \! \left(C^{2}D \! + \! 2C \sqrt{\smash[b]{-AB}}- 
\! (C \! - \! \tau A) \mathcal{Y} \right), \qquad D_{1} \! = \! \dfrac{4 
\tau^{2}}{C}, \\
\mathcal{Y} \! := \! \dfrac{1 \! - \! a \mi}{\tau^{2}} \! + \! \frac{A}{\tau 
C}, \qquad \quad \mi a_{1} \! = \! \mi a \! - \! 2.
\end{gather*}

For the purpose of the application of these transformations to asymptotic 
results, one has to find the action of the group of B\"{a}cklund 
transformations on the manifold of monodromy data. In terms of $\Psi (\mu,
\tau)$ introduced in Proposition~2.1, the transformation $\mathscr{R}_{1,2}$ 
reads
\begin{equation*}
\Psi_{1}(\mu,\tau) \! = \! \dfrac{\mu}{\sqrt{\tau}} \! \left(\! 
\begin{pmatrix}
0 & 0 \\
0 & 1
\end{pmatrix} \! + \! \frac{1}{\mu} \! 
\begin{pmatrix}
0 & \frac{2 \mi \tau}{D} \\
-\frac{D}{2 \mi} & 0
\end{pmatrix} \! \right) \! \Psi (\mu,\tau):
\end{equation*}
applying the latter relation to the canonical solutions, one shows that
\begin{equation*}
\mathscr{R}_{1,2}Y^{\infty}_{k}(\mu) \! = \! (Y^{\infty}_{k}(\mu))_{1}, \qquad 
\qquad \mathscr{R}_{1,2}X^{0}_{k}(\mu) \! = \! \mi (X^{0}_{k}(\mu))_{1} 
\sigma_{3};
\end{equation*} 
hence,
\begin{equation*}
(a,s_{0}^{0},s_{0}^{\infty},s_{1}^{\infty},g_{11},g_{12},g_{21},g_{22}) \! 
\to \! (a \! - \! \mi, -s_{0}^{0},s_{0}^{\infty},s_{1}^{\infty},\mi g_{11},
\mi g_{12},-\mi g_{21},-\mi g_{22}).
\end{equation*}
\subsection{Lie-Point Symmetries}
\label{subsec:LPS}
In this subsection, the subscripts ``$n$'' and ``$o$'' are used to denote 
``new'' and ``old'' variables, respectively.
\subsubsection{$\tau \! \to \! -\tau$}
Let $A_{o}(\tau_{o})$, $B_{o}(\tau_{o})$, $C_{o}(\tau_{o})$, and $D_{o}
(\tau_{o})$ be a solution of System~(\ref{sys:ABCD}), and the function 
$u_{o}(\tau_{o})$ solve Equation~(\ref{eq:dp3}). Set $\tau_{n} \! = \! 
\tau_{o} \me^{-\mi \pi p}$, $p \! = \! \pm 1$. Then
\begin{equation*}
A_{n}(\tau_{n}) \! = \! -A_{o}(\tau_{o}), \quad B_{n}(\tau_{n}) \! = \! -B_{o}
(\tau_{o}), \quad C_{n}(\tau_{n}) \! = \! C_{o}(\tau_{o}), \quad \text{and} 
\quad D_{n}(\tau_{n}) \! = \! D_{o}(\tau_{o})
\end{equation*}
is a solution of System~(\ref{sys:ABCD}) for $\sqrt{\smash[b]{-A_{n}(\tau_{n})
B_{n}(\tau_{n})}}=\sqrt{\smash[b]{-A_{o}(\tau_{o})B_{o}(\tau_{o})}}$, and the 
function $u_{n}(\tau_{n}) \! = \! -u_{o}(\tau_{o})$ solves 
Equation~(\ref{eq:dp3}). Note that $\widetilde{\alpha}_{o}(\tau_{o}) \! = \! 
\widetilde{\alpha}_{n}(\tau_{n})$ and $-\mi \widetilde{\alpha}_{n}(\tau_{n})
B_{n}(\tau_{n}) \! = \! \varepsilon b$, $\varepsilon \! = \! \pm 1$. On the 
corresponding fundamental solutions of Systems~(\ref{eq:UV}) and~(12), the 
above transformation acts as follows:
\begin{equation*}
\lambda_{o} \! = \! -\lambda_{n}, \qquad \Phi_{o}(\lambda_{o},\tau_{o}) \! = 
\! \sigma_{3} \Phi_{n}(\lambda_{n},\tau_{n})
\end{equation*}
and
\begin{equation*}
\mu_{o} \! = \! \mu_{n} \me^{\frac{\mi \pi l}{2}}, \quad l \! = \! \pm 1, 
\qquad \Psi_{o}(\mu_{o},\tau_{o}) \! = \! \me^{\frac{\mi \pi l}{4} \sigma_{3}} 
\Psi_{n}(\mu_{n},\tau_{n}).
\end{equation*}
In terms of the canonical solutions of System~(12), this action reads:
\begin{gather*}
Y^{\infty}_{o,k}(\mu_{o}) \! = \! \me^{\frac{\mi \pi l}{4} \sigma_{3}}
Y^{\infty}_{n,k-p-l}(\mu_{n}) \me^{-\frac{\mi \pi l}{4} \sigma_{3}} \me^{
\frac{\pi l}{2}(a-\mi /2) \sigma_{3}}, \\
X^{0}_{o,k}(\mu_{o}) \! = \! 
\begin{cases}
\me^{\frac{\mi \pi l}{4} \sigma_{3}}X^{0}_{n,k}(\mu_{n}), &\text{$p \! = \! 
l$,} \\
-\mi l \me^{\frac{\mi \pi l}{4} \sigma_{3}}X^{0}_{n,k-l}(\mu_{n}) \sigma_{1}, 
&\text{$p \! = \! -l$.}
\end{cases}
\end{gather*}
These formulae for the canonical solutions imply the following action on 
$\mathscr{M}$:
\begin{gather*}
S^{\infty}_{n,k-p-l} \! = \! \me^{-\frac{\mi \pi l}{4} \sigma_{3}} \me^{\frac{
\pi l}{2}(a-\mi /2) \sigma_{3}}S^{\infty}_{o,k} \me^{-\frac{\pi l}{2}(a-\mi 
/2) \sigma_{3}} \me^{\frac{\mi \pi l}{4} \sigma_{3}}, \\
S^{0}_{o,k} \! = \! 
\begin{cases}
S^{0}_{n,k}, &\text{$p \! = \! l$,} \\
\sigma_{1}S^{0}_{n,k-l} \sigma_{1}, &\text{$p \! = \! -l$,}
\end{cases} \\
G_{o} \! = \! 
\begin{cases}
-\mi \sigma_{1}(S^{0}_{n,0})^{-1}G_{n} \me^{\frac{\mi \pi}{4} \sigma_{3}} 
\me^{-\frac{\pi}{2}(a-\mi /2) \sigma_{3}}, &\text{$p \! = \! 1$,} \\
\mi S^{0}_{n,0} \sigma_{1}G_{n} \me^{-\frac{\mi \pi}{4} \sigma_{3}} 
\me^{\frac{\pi}{2}(a-\mi /2) \sigma_{3}}, &\text{$p \! = \! -1$.}
\end{cases}
\end{gather*}
Note that the action on $\mathscr{M}$ is independent of $l$.

One uses this transformation to prove the asymptotic results stated in 
Theorems~3.1--3.5 and formulate Conjectures~3.1 and~3.2 for $\varepsilon_{1} 
\! = \! \pm 1$ (negative $\tau)$ by using those for $\varepsilon_{1} \! = \! 
0$ (positive $\tau)$.
\subsubsection{$a \! \to \! -a$}
Let $A_{o}(\tau_{o})$, $B_{o}(\tau_{o})$, $C_{o}(\tau_{o})$, and $D_{o}(\tau_{
o})$ be a solution of System~(\ref{sys:ABCD}) for $a \! = \! a_{o}$, and the 
function $u_{o}(\tau_{o})$ solve Equation~(\ref{eq:dp3}) for $a \! = \! 
a_{o}$, $b \! = \! b_{o}$, and $\varepsilon \! = \! \varepsilon_{o}$ $(= \! 
\pm 1)$. Set $a_{n} \! = \! -a_{o}$, and either $b_{n} \! = \! b_{o}$ and 
$\varepsilon_{n} \! = \! \varepsilon_{o} \me^{-\mi \pi p}$, $p \! = \! \pm 1$, 
or $b_{n} \! = \! b_{o} \me^{-\mi \pi p}$ and $\varepsilon_{n} \! = \! 
\varepsilon_{o}$. Then
\begin{equation*}
A_{n}(\tau_{n}) \! = \! -B_{o}(\tau_{o}), \quad B_{n}(\tau_{n}) \! = \! -A_{o}
(\tau_{o}), \quad C_{n}(\tau_{n}) \! = \! D_{o}(\tau_{o}), \quad \text{and} 
\quad D_{n}(\tau_{n}) \! = \! C_{o}(\tau_{o})
\end{equation*}
is a solution of System~(\ref{sys:ABCD}) for $a \! = \! -a_{n}$ and $\sqrt{
\smash[b]{-A_{n}(\tau_{n})B_{n}(\tau_{n})}}=\sqrt{\smash[b]{-A_{o}(\tau_{o})
B_{o}(\tau_{o})}}$, and the function $u_{n}(\tau_{n}) \! = \! \varepsilon_{o} 
\varepsilon_{n}u_{o}(\tau_{o})$ solves Equation~(\ref{eq:dp3}) for $a \! = \! 
-a_{n}$, $b \! = \! b_{n}$, and $\varepsilon \! = \! \varepsilon_{n}$. Note 
that $\widetilde{\alpha}_{o}(\tau_{o}) \! = \! \tfrac{B_{n}(\tau_{n})}{A_{n}
(\tau_{n})} \widetilde{\alpha}_{n}(\tau_{n})$ and $-\mi \widetilde{\alpha}_{n}
(\tau_{n})B_{n}(\tau_{n}) \! = \! \varepsilon_{n}b_{n}$ $(\varepsilon_{n} \! = 
\! \pm 1)$. On the corresponding fundamental solutions of 
Systems~(\ref{eq:UV}) and~(12), the above transformation acts as follows:
\begin{equation*}
\lambda_{o} \! = \! -\lambda_{n}, \qquad \Phi_{o}(\lambda_{o},\tau_{o}) \! = 
\! \sigma_{1} \Phi_{n}(\lambda_{n},\tau_{n})
\end{equation*}
and
\begin{equation*}
\mu_{o} \! = \! \mu_{n} \me^{\frac{\mi \pi l}{2}}, \quad l \! = \! \pm 1, 
\qquad \Psi_{o}(\mu_{o},\tau_{o}) \! = \! \mathcal{Q}(\mu_{n},\tau_{n}) 
\Psi_{n}(\mu_{n},\tau_{n}),
\end{equation*} 
where
\begin{equation*}
\mathcal{Q}(\mu_{n},\tau_{n}) \! := \! \mu_{n} \me^{\frac{\mi \pi l}{4}} \! 
\begin{pmatrix}
0 & 0 \\
1 & 0
\end{pmatrix} \! - \! \left(\dfrac{B_{n}(\tau_{n}) \me^{-\frac{\mi \pi l}{4}}
}{\sqrt{\smash[b]{-A_{n}(\tau_{n})B_{n}(\tau_{n})}}} \right)^{\sigma_{3}} 
\sigma_{3}.
\end{equation*}
In terms of the canonical solutions of System~(12), this action reads:
\begin{gather*}
Y^{\infty}_{o,k}(\mu_{o}) \! = \! \mathcal{Q}(\mu_{n},\tau_{n})Y^{\infty}_{n,k
-l}(\mu_{n}) \me^{\frac{a_{n} \pi l}{2} \sigma_{3}} \sigma_{1}, \\
X^{0}_{o,k}(\mu_{o}) \! = \! 
\begin{cases}
-\mi p \mathcal{Q}(\mu_{n},\tau_{n})X^{0}_{n,k}(\mu_{n}), &\text{$p \! = \! 
l$,} \\
-\mathcal{Q}(\mu_{n},\tau_{n})X^{0}_{n,k+p}(\mu_{n}) \sigma_{1}, &\text{$p \! 
= \! -l$.}
\end{cases}
\end{gather*}
These formulae for the canonical solutions imply the following action on 
$\mathscr{M}$:
\begin{gather*}
S^{\infty}_{n,k-l} \! = \! \me^{\frac{a_{n} \pi l}{2} \sigma_{3}} \sigma_{1}
S^{\infty}_{o,k} \sigma_{1} \me^{-\frac{a_{n} \pi l}{2} \sigma_{3}}, \\
S^{0}_{o,k} \! = \! 
\begin{cases}
S^{0}_{n,k}, &\text{$p \! = \! l$,} \\
\sigma_{1}S^{0}_{n,k+p} \sigma_{1}, &\text{$p \! = \! -l$,}
\end{cases} \\
G_{o} \! = \! 
\begin{cases}
\mi G_{n} \me^{\pi (a_{n}-\mi /2) \sigma_{3}} \sigma_{3}(S^{\infty}_{n,1})^{-
1} \sigma_{3} \me^{-\pi (a_{n}-\mi /2) \sigma_{3}} \me^{\frac{a_{n} \pi}{2} 
\sigma_{3}} \sigma_{1}, &\text{$p \! = \! 1$,} \\
-S^{0}_{n,0} \sigma_{1}G_{n} \me^{\pi (a_{n}-\mi /2) \sigma_{3}} \sigma_{3}
(S^{\infty}_{n,1})^{-1} \sigma_{3} \me^{-\pi (a_{n}-\mi /2) \sigma_{3}} 
\me^{\frac{a_{n} \pi}{2} \sigma_{3}} \sigma_{1}, &\text{$p \! = \! -1$.}
\end{cases}
\end{gather*}
Note that the action on $\mathscr{M}$ is independent of $l$.

One uses this transformation to prove the asymptotic results stated in 
Theorems~3.1--3.5 and formulate Conjectures~3.1 and~3.2 for $\varepsilon_{2} 
\! = \! \pm 1$ (negative $\varepsilon b)$ by using those for $\varepsilon_{2} 
\! = \! 0$ (positive $\varepsilon b)$.
\subsubsection{$\tau \! \to \! \mi \tau$}
Let $A_{o}(\tau_{o})$, $B_{o}(\tau_{o})$, $C_{o}(\tau_{o})$, and $D_{o}(\tau_{
o})$ be a solution of System~(\ref{sys:ABCD}), and the function $u_{o}(\tau_{o}
)$ solve Equation~(\ref{eq:dp3}) for $b \! = \! b_{o}$ and $\varepsilon \! = 
\! \varepsilon_{o}$. Set $\tau_{n} \! = \! -\mi l \tau_{o}$, $l \! = \! \pm 
1$, and either $b_{n} \! = \! b_{o}$ and $\varepsilon_{n} \! = \! 
\varepsilon_{o} \me^{-\mi \pi p}$, $p \! = \! \pm 1$, or $b_{n} \! = \! b_{o} 
\me^{-\mi \pi p}$ and $\varepsilon_{n} \! = \! \varepsilon_{o}$. Then
\begin{equation*}
A_{n}(\tau_{n}) \! = \! -A_{o}(\tau_{o}), \quad B_{n}(\tau_{n}) \! = \! -B_{o}
(\tau_{o}), \quad C_{n}(\tau_{n}) \! = \! \mi l C_{o}(\tau_{o}), \quad 
\text{and} \quad D_{n}(\tau_{n}) \! = \! \mi l D_{o}(\tau_{o})
\end{equation*}
is a solution of System~(\ref{sys:ABCD}) for $\sqrt{\smash[b]{-A_{n}(\tau_{n})
B_{n}(\tau_{n})}}=-\sqrt{\smash[b]{-A_{o}(\tau_{o})B_{o}(\tau_{o})}}$, and the 
function $u_{n}(\tau_{n}) \! = \! \mi \varepsilon_{o} \varepsilon_{n}lu_{o}
(\tau_{o})$ solves Equation~(\ref{eq:dp3}) for $b \! = \! b_{n}$ and 
$\varepsilon \! = \! \varepsilon_{n}$. Note that $\widetilde{\alpha}_{o}
(\tau_{o}) \! = \! \widetilde{\alpha}_{n}(\tau_{n})$ and $-\mi \widetilde{
\alpha}_{n}(\tau_{n})B_{n}(\tau_{n}) \! = \! \varepsilon_{n}b_{n}$ 
$(\varepsilon_{n} \! = \! \pm 1)$. On the corresponding fundamental solutions 
of Systems~(\ref{eq:UV}) and~(12), the above transformation acts as follows:
\begin{equation*}
\lambda_{o} \! = \! -\mi l \lambda_{n}, \qquad \Phi_{o}(\lambda_{o},\tau_{o}) 
\! = \! \Phi_{n}(\lambda_{n},\tau_{n})
\end{equation*}
and
\begin{equation*}
\mu_{o} \! = \! \mu_{n} \me^{-\frac{\mi \pi l}{4}}, \qquad \Psi_{o}(\mu_{o},
\tau_{o}) \! = \! \me^{\frac{\mi \pi l}{8} \sigma_{3}} \Psi_{n}(\mu_{n},
\tau_{n}).
\end{equation*}
In terms of the canonical solutions of System~(12), this action reads:
\begin{gather*}
Y^{\infty}_{o,k}(\mu_{o}) \! = \! \me^{\frac{\mi \pi l}{8} \sigma_{3}}
Y^{\infty}_{n,k}(\mu_{n}) \me^{-\frac{\pi la}{4} \sigma_{3}}, \\
X^{0}_{o,k}(\mu_{o}) \! = \! 
\begin{cases}
\me^{\frac{\mi \pi l}{8} \sigma_{3}}X^{0}_{n,k}(\mu_{n}), &\text{$p \! = \! 
-l$,} \\
\mi l \me^{\frac{\mi \pi l}{8} \sigma_{3}}X^{0}_{n,k+l}(\mu_{n}) \sigma_{1}, 
&\text{$p \! = \! l$.}
\end{cases}
\end{gather*}
These formulae for the canonical solutions imply the following action on 
$\mathscr{M}$:
\begin{gather*}
S^{\infty}_{n,k} \! = \! \me^{-\frac{\pi la}{4} \sigma_{3}}S^{\infty}_{o,k} 
\me^{\frac{\pi la}{4} \sigma_{3}}, \\
S^{0}_{o,k} \! = \! 
\begin{cases}
S^{0}_{n,k}, &\text{$p \! = \! -l$,} \\
\sigma_{1}S^{0}_{n,k+l} \sigma_{1}, &\text{$p \! = \! l$,}
\end{cases} \\
G_{o} \! = \! 
\begin{cases}
\mi S^{0}_{n,0} \sigma_{1}G_{n} \me^{\frac{\pi a}{4} \sigma_{3}}, &\text{$p \! 
= \! l \! = \! -1$,} \\
G_{n} \me^{-\frac{\pi a}{4} \sigma_{3}}, &\text{$p \! = \! -l \! = \! -1$,} \\
G_{n} \me^{\frac{\pi a}{4} \sigma_{3}}, &\text{$p \! = \! -l \! = \! 1$,} \\
-\mi \sigma_{1}(S^{0}_{n,0})^{-1}G_{n} \me^{-\frac{\pi a}{4} \sigma_{3}}, 
&\text{$p \! = \! l \! = \! 1$.} \\
\end{cases}
\end{gather*}

One uses this transformation to prove the asymptotic results stated in (see 
the Appendix) Theorems~A.1--A.5 and formulate Conjectures~A.1 and~A.2 for 
(pure) imaginary $\tau$ by using the results for real $\tau$ in Section~3.

\vspace{1.5cm}

\textbf{Acknowledgements}

The authors are grateful to the Department of Pure Mathematics of the 
University of Adelaide, where part of this work was done, for hospitality. 
A.~V.~K. was partially supported by ARC grant \# A69803721 and RFBR grant \# 
01-01-01045.
\clearpage
\section*{Appendix: Asymptotics for Imaginary $\tau$}
\setcounter{section}{1}
\setcounter{equation}{0}
\renewcommand{\thesection}{\Alph{section}}
\renewcommand{\theequation}{\Alph{section}.\arabic{equation}}
Here we present the summary of results for asymptotics of $u(\tau)$, $\mathcal{
H}(\tau)$, and $\pmb{\pmb{\boldsymbol{\tau}}}(\tau)$ as $\tau \! \to \! \pm 
\mi 0$ and $\tau \! \to \! \pm \mi \infty$. These results are obtained by 
applying the transformations changing $\tau \! \to \! \mi \tau$ given in 
Subsection~6.2 to asymptotic results for real values of $\tau$ presented in 
Section~3.

In order to present these results, it is convenient to introduce the mapping 
$\widehat{\mathscr{F}}_{\varepsilon_{1},\varepsilon_{2}} \colon \mathscr{M} \! 
\to \! \mathscr{M}$, $(a,s^{0}_{0},s^{\infty}_{0},s^{\infty}_{1},g_{11},g_{12},
g_{21},g_{22}) \! \! \to \! \! (a,s^{0}_{0},\widehat{s}^{\infty}_{0}
(\varepsilon_{1},\varepsilon_{2}),\widehat{s}^{\infty}_{1}(\varepsilon_{1},
\varepsilon_{2}),\widehat{g}_{11}(\varepsilon_{1},\varepsilon_{2}),\widehat{
g}_{12}(\varepsilon_{1},\varepsilon_{2}),\widehat{g}_{21}(\varepsilon_{1},
\varepsilon_{2}),\linebreak[4]
\widehat{g}_{22}(\varepsilon_{1},\varepsilon_{2}))$, $\varepsilon_{1} \! = \! 
\pm 1$, $\varepsilon_{2} \! = \! 0,\pm 1$. Define:
\begin{enumerate}
\item[(1)] $\widehat{\mathscr{F}}_{-1,0}$ as: $\widehat{s}^{\infty}_{0}(-1,0) 
\! = \! s^{\infty}_{1} \me^{-\frac{3 \pi a}{2}}$, $\widehat{s}^{\infty}_{1}
(-1,0) \! = \! s^{\infty}_{0} \me^{-\frac{\pi a}{2}}$, $\widehat{g}_{11}(-1,
0) \! = \! -g_{22} \me^{-\frac{3 \pi a}{4}}$, $\widehat{g}_{12}(-1,0) 
\linebreak[4]
= \! -(g_{21} \! + \! s^{\infty}_{0}g_{22}) \me^{\frac{3 \pi a}{4}}$, 
$\widehat{g}_{21}(-1,0) \!= \! -(g_{12} \! - \! s^{0}_{0}g_{22}) \me^{-\frac{3 
\pi a}{4}}$, and $\widehat{g}_{22}(-1,0) \! = \! -(g_{11} \! + \! s^{\infty}_{
0}g_{12} \! - \! (g_{21} \! + \! s^{\infty}_{0}g_{22})s^{0}_{0}) \me^{\frac{3 
\pi a}{4}}$;
\item[(2)] $\widehat{\mathscr{F}}_{-1,-1}$ as: $\widehat{s}^{\infty}_{0}(-1,
-1) \! = \! s^{\infty}_{0} \me^{-\frac{\pi a}{2}}$, $\widehat{s}^{\infty}_{1}
(-1,-1) \! = \! s^{\infty}_{1} \me^{\frac{\pi a}{2}}$, $\widehat{g}_{11}(-1,
-1) \! = \! -\mi g_{21} \me^{-\frac{\pi a}{4}}$, $\widehat{g}_{12}(-1,-1) \! = 
\! -\mi g_{22} \me^{\frac{\pi a}{4}}$, $\widehat{g}_{21}(-1,-1) \! = \! -\mi 
(g_{11} \! - \! s^{0}_{0}g_{21}) \me^{-\frac{\pi a}{4}}$, and $\widehat{g}_{22}
(-1,-1) \! = \! -\mi (g_{12} \! - \! s^{0}_{0}g_{22}) \me^{\frac{\pi a}{4}}$;
\item[(3)] $\widehat{\mathscr{F}}_{-1,1}$ as: $\widehat{s}^{\infty}_{0}(-1,1) 
\! = \! s^{\infty}_{0} \me^{-\frac{\pi a}{2}}$, $\widehat{s}^{\infty}_{1}(-1,
1) \! = \! s^{\infty}_{1} \me^{\frac{\pi a}{2}}$, $\widehat{g}_{11}(-1,1) \! = 
\! g_{11} \me^{-\frac{\pi a}{4}}$, $\widehat{g}_{12}(-1,1) \! = \! g_{12} \me^{
\frac{\pi a}{4}}$, $\widehat{g}_{21}(-1,1) \! = \! g_{21} \me^{-\frac{\pi a}{
4}}$, and $\widehat{g}_{22}(-1,1) \! = \! g_{22} \me^{\frac{\pi a}{4}}$;
\item[(4)] $\widehat{\mathscr{F}}_{1,0}$ as: $\widehat{s}^{\infty}_{0}(1,0) \! 
= \! s^{\infty}_{1} \me^{-\frac{\pi a}{2}}$, $\widehat{s}^{\infty}_{1}(1,0) 
\! = \! s^{\infty}_{0} \me^{-\frac{3 \pi a}{2}}$, $\widehat{g}_{11}(1,0) \! = 
\! -\mi g_{12} \me^{-\frac{\pi a}{4}}$, $\widehat{g}_{12}(1,0) \! = \! -\mi 
(g_{11} \! + \! s^{\infty}_{0}g_{12}) \me^{\frac{\pi a}{4}}$, $\widehat{g}_{21}
(1,0) \! = \! -\mi g_{22} \me^{-\frac{\pi a}{4}}$, and $\widehat{g}_{22}(1,0) 
\! = \! -\mi (g_{21} \! + \! s^{\infty}_{0}g_{22}) \me^{\frac{\pi a}{4}}$;
\item[(5)] $\widehat{\mathscr{F}}_{1,-1}$ as: $\widehat{s}^{\infty}_{0}(1,-1) 
\! = \! s^{\infty}_{0} \me^{\frac{\pi a}{2}}$, $\widehat{s}^{\infty}_{1}(1,-1) 
\! = \! s^{\infty}_{1} \me^{-\frac{\pi a}{2}}$, $\widehat{g}_{11}(1,-1) \! = 
\! g_{11} \me^{\frac{\pi a}{4}}$, $\widehat{g}_{12}(1,-1) \! = \! g_{12} \me^{
-\frac{\pi a}{4}}$, $\widehat{g}_{21}(1,-1) \! = \! g_{21} \me^{\frac{\pi a}{
4}}$, and $\widehat{g}_{22}(1,-1) \! = \! g_{22} \me^{-\frac{\pi a}{4}}$; and
\item[(6)] $\widehat{\mathscr{F}}_{1,1}$ as: $\widehat{s}^{\infty}_{0}(1,1) \! 
= \! s^{\infty}_{0} \me^{\frac{\pi a}{2}}$, $\widehat{s}^{\infty}_{1}(1,1) \! 
= \! s^{\infty}_{1} \me^{-\frac{\pi a}{2}}$, $\widehat{g}_{11}(1,1) \! = \! 
\mi (g_{21} \! + \! s^{0}_{0}g_{11}) \me^{\frac{\pi a}{4}}$, $\widehat{g}_{12}
(1,1) \! = \! \mi (g_{22} \! + \! s^{0}_{0}g_{12}) \me^{-\frac{\pi a}{4}}$, 
$\widehat{g}_{21}(1,1) \! = \! \mi g_{11} \me^{\frac{\pi a}{4}}$, and 
$\widehat{g}_{22}(1,1) \! = \! \mi g_{12} \me^{-\frac{\pi a}{4}}$.
\end{enumerate}
\begin{ddddd}
Let $\varepsilon_{1} \! = \! \pm 1$, $\varepsilon_{2} \! = \! 0,\pm 1$, 
$\varepsilon b \! = \! \vert \varepsilon b \vert \me^{\mi \pi \varepsilon_{
2}}$, and $u(\tau)$ be a solution of Equation~{\rm (\ref{eq:dp3})} 
corresponding to the monodromy data 
$(a,s^{0}_{0},s^{\infty}_{0},s^{\infty}_{1},g_{11},g_{12},g_{21},g_{22})$. 
Suppose that
\begin{equation*}
\widehat{g}_{11}(\varepsilon_{1},\varepsilon_{2}) \widehat{g}_{12}
(\varepsilon_{1},\varepsilon_{2}) \widehat{g}_{21}(\varepsilon_{1},
\varepsilon_{2}) \widehat{g}_{22}(\varepsilon_{1},\varepsilon_{2}) \! \not= 
\! 0, \qquad \left \vert \Re \! \left(\tfrac{\mi}{2 \pi} \ln (\widehat{g}_{11}
(\varepsilon_{1},\varepsilon_{2}) \widehat{g}_{22}(\varepsilon_{1},
\varepsilon_{2})) \right) \right \vert \! < \! \dfrac{1}{6}.
\end{equation*}
Then $\exists \, \, \delta \! > \! 0$ such that $u(\tau)$ has the asymptotic 
expansion
\begin{align}
u(\tau) \underset{\tau \to \infty \me^{\frac{\mi \pi \varepsilon_{1}}{2}}}{=}& 
\, \dfrac{(-1)^{\frac{1+\varepsilon_{1}}{2}} \mi \varepsilon \sqrt{\smash[b]{
\vert \varepsilon b \vert}}}{3^{1/4}} \! \left(\! \sqrt{\dfrac{\vartheta 
(\tau)}{12}}+ \! \sqrt{\widehat{\nu}(\varepsilon_{1},\varepsilon_{2}) \! + \! 
1} \, \me^{\frac{3 \pi \mi}{4}} \cosh \! \left(\mi \vartheta (\tau) \! + \! 
(\widehat{\nu}(\varepsilon_{1},\varepsilon_{2}) \! + \! 1) \right. \right. 
\nonumber \\
\times&\left. \left. \ln \vartheta (\tau) \! + \! \widehat{z}(\varepsilon_{1},
\varepsilon_{2}) \! + \! o \! \left(\tau^{-\delta} \right) \right) \right),
\end{align}
where
\begin{equation*}
\vartheta (\tau) \! := \! 3 \sqrt{\smash[b]{3}} \, \vert \varepsilon b \vert^{
1/3} \vert \tau \vert^{2/3}, \qquad \quad \widehat{\nu}(\varepsilon_{1},
\varepsilon_{2}) \! + \! 1 \! := \! \dfrac{\mi}{2 \pi} \ln (\widehat{g}_{11}
(\varepsilon_{1},\varepsilon_{2}) \widehat{g}_{22}(\varepsilon_{1},
\varepsilon_{2})),
\end{equation*} 
\begin{align*}
\widehat{z}(\varepsilon_{1},\varepsilon_{2}) :=& \, \dfrac{1}{2} \ln (2 \pi) 
\! - \! \dfrac{\pi \mi}{2} \! - \! \dfrac{3 \pi \mi}{2} \! \left(\widehat{\nu}
(\varepsilon_{1},\varepsilon_{2}) \! + \! 1 \right) \! + \! (-1)^{\varepsilon_{
2}} \mi a \ln \! \left(2 \! + \! \sqrt{\smash[b]3} \right) \! + \! \left(
\widehat{\nu}(\varepsilon_{1},\varepsilon_{2}) \! + \! 1 \right) \! \ln 12 \\
-& \, \ln \! \left(\widehat{\omega}(\varepsilon_{1},\varepsilon_{2}) \sqrt{
\smash[b]{\widehat{\nu}(\varepsilon_{1},\varepsilon_{2}) \! + \! 1}} \, \Gamma 
(\widehat{\nu}(\varepsilon_{1},\varepsilon_{2}) \! + \! 1) \right),
\end{align*}
with
\begin{equation*}
\widehat{\omega}(\varepsilon_{1},\varepsilon_{2}) \! := \! \dfrac{\widehat{g}_{
12}(\varepsilon_{1},\varepsilon_{2})}{\widehat{g}_{22}(\varepsilon_{1},
\varepsilon_{2})},
\end{equation*}
and $\Gamma (\boldsymbol{\cdot})$ is the gamma function.

Let $\mathcal{H}(\tau)$ be the Hamiltonian function defined in 
Equation~{\rm (35)} corresponding to the function $u(\tau)$ given above. Then
\begin{align}
\mathcal{H}(\tau) \underset{\tau \to \infty \me^{\frac{\mi \pi \varepsilon_{
1}}{2}}}{=}& \, (-1)^{\frac{1+\varepsilon_{1}}{2}} \mi \! \left(3(\varepsilon 
b)^{2/3} \vert \tau \vert^{1/3} \! + \! 2 \vert \varepsilon b \vert^{1/3} 
\vert \tau \vert^{-1/3} \! \left(\! \left(a \! - \! (-1)^{\varepsilon_{2}} 
\mi/2 \right) \! - \! 2 \sqrt{\smash[b]{3}} \, \mi \! \left(\widehat{\nu}
(\varepsilon_{1},\varepsilon_{2}) \! + \! 1 \right) \right. \right. \nonumber 
\\
+&\left. \left. o \! \left(\tau^{-\delta} \right) \right) \! + \! \dfrac{
(a \! - \! (-1)^{\varepsilon_{2}} \mi/2)^{2}}{2 \tau} \right).
\end{align}
\end{ddddd}
\begin{ddddd}
Let $\varepsilon_{1} \! = \! \pm 1$, $\varepsilon_{2} \! = \! 0,\pm 1$, 
$\varepsilon b \! = \! \vert \varepsilon b \vert \me^{\mi \pi \varepsilon_{
2}}$, and $u(\tau)$ be a solution of Equation~{\rm (\ref{eq:dp3})} 
corresponding to the monodromy data 
$(a,s^{0}_{0},s^{\infty}_{0},s^{\infty}_{1},g_{11},g_{12},g_{21},g_{22})$. 
Suppose that
\begin{equation*}
\widehat{g}_{21}(\varepsilon_{1},\varepsilon_{2}) \! = \! 0, \qquad \quad 
\widehat{g}_{11}(\varepsilon_{1},\varepsilon_{2}) \widehat{g}_{22}
(\varepsilon_{1},\varepsilon_{2}) \! = \! 1.
\end{equation*}
Then $\exists \, \, \delta \! > \! 0$ such that $u(\tau)$ has the asymptotic 
expansion
\begin{align}
u(\tau) \underset{\tau \to \infty \me^{\frac{\mi \pi \varepsilon_{1}}{2}}}{=}& 
\, (-1)^{\frac{1+\varepsilon_{1}}{2}} \mi \! \left(\dfrac{\varepsilon 
(\varepsilon b)^{2/3}}{2} \vert \tau \vert^{1/3} \! + \! \dfrac{(-1)^{
\varepsilon_{1}} \varepsilon \sqrt{\smash[b]{\vert \varepsilon b \vert}} \, 
(s_{0}^{0} \! - \! \mi \me^{(-1)^{\varepsilon_{2}+1} \pi a})}{2^{3/2}3^{1/4} 
\sqrt{\smash[b]{\pi}}} \! \left(\! \dfrac{\sqrt{\smash[b]{3}}- \! 1}{\sqrt{
\smash[b]{3}}+ \! 1} \right)^{(-1)^{\varepsilon_{2}} \mi a} \right. \nonumber 
\\
\times&\left. \, \exp \! \left(-\mi \! \left(3 \sqrt{\smash[b]{3}} \, \vert 
\varepsilon b \vert^{1/3} \vert \tau \vert^{2/3} \! - \! \dfrac{\pi}{4} 
\right) \right) \right) \! \left(1 \! + \! o \! \left(\tau^{-\delta} \right) 
\right).
\end{align}
\end{ddddd}
\begin{ddddd}
Let $\varepsilon_{1} \! = \! \pm 1$, $\varepsilon_{2} \! = \! 0,\pm 1$, 
$\varepsilon b \! = \! \vert \varepsilon b \vert \me^{\mi \pi \varepsilon_{
2}}$, and $u(\tau)$ be a solution of Equation~{\rm (\ref{eq:dp3})} 
corresponding to the monodromy data 
$(a,s^{0}_{0},s^{\infty}_{0},s^{\infty}_{1},g_{11},g_{12},g_{21},g_{22})$. 
Suppose that
\begin{equation*}
\widehat{g}_{12}(\varepsilon_{1},\varepsilon_{2}) \! = \! 0, \qquad \quad 
\widehat{g}_{11}(\varepsilon_{1},\varepsilon_{2}) \widehat{g}_{22}
(\varepsilon_{1},\varepsilon_{2}) \! = \! 1.
\end{equation*}
Then $\exists \, \, \delta \! > \! 0$ such that $u(\tau)$ has the asymptotic 
expansion
\begin{align}
u(\tau) \underset{\tau \to \infty \me^{\frac{\mi \pi \varepsilon_{1}}{2}}}{=}& 
\, (-1)^{\frac{1+\varepsilon_{1}}{2}} \mi \! \left(\dfrac{\varepsilon 
(\varepsilon b)^{2/3}}{2} \vert \tau \vert^{1/3} \! + \! \dfrac{(-1)^{
\varepsilon_{1}} \varepsilon \sqrt{\smash[b]{\vert \varepsilon b \vert}} \, 
(s_{0}^{0} \! - \! \mi \me^{(-1)^{\varepsilon_{2}+1} \pi a})}{2^{3/2} 3^{1/4} 
\sqrt{\smash[b]{\pi}}} \! \left(\! \dfrac{\sqrt{\smash[b]{3}}+ \! 1}{\sqrt{
\smash[b]{3}}- \! 1} \right)^{(-1)^{\varepsilon_{2}} \mi a} \right. \nonumber 
\\
\times&\left. \, \exp \! \left(\mi \! \left(3 \sqrt{\smash[b]{3}} \, \vert 
\varepsilon b \vert^{1/3} \vert \tau \vert^{2/3} \! + \! \dfrac{3 \pi}{4} 
\right) \right) \right) \! \left(1 \! + \! o \! \left(\tau^{-\delta} \right) 
\right).
\end{align}
\end{ddddd}
\begin{ddddd}
Let $\varepsilon_{1} \! = \! \pm 1$, $\varepsilon_{2} \! = \! 0,\pm 1$, 
$\varepsilon b \! = \! \vert \varepsilon b \vert \me^{\mi \pi \varepsilon_{
2}}$, and $u(\tau)$ be a solution of Equation~{\rm (\ref{eq:dp3})} 
corresponding to the monodromy data 
$(a,s^{0}_{0},s^{\infty}_{0},s^{\infty}_{1},g_{11},g_{12},g_{21},g_{22})$. 
Suppose that
\begin{equation*}
\vert \Im (a) \vert \! < \! 1, \qquad \quad \widehat{g}_{11}(\varepsilon_{1},
\varepsilon_{2}) \widehat{g}_{22}(\varepsilon_{1},\varepsilon_{2}) \! \not= \! 
0, \qquad \quad \rho \! \not= \! 0, \qquad \quad \vert \Re (\rho) \vert \! < 
\! \dfrac{1}{2},
\end{equation*}
where
\begin{eqnarray*}
\cos (2 \pi \rho) \! := \! -\dfrac{\mi s^{0}_{0}}{2} \! = \cosh (\pi a) \! + 
\! \dfrac{1}{2} s_{0}^{\infty}s_{1}^{\infty} \me^{\pi a}.
\end{eqnarray*}
Then $\exists \, \, \delta \! > \! 0$ such that $u(\tau)$ has the asymptotic 
expansion
\begin{align}
u(\tau) \underset{\tau \to 0 \me^{\frac{\mi \pi \varepsilon_{1}}{2}}}{=}& \, 
\dfrac{(-1)^{\varepsilon_{2}} \tau b}{16 \pi} \exp \! \left((-1)^{\varepsilon_{
2}} \dfrac{\pi a}{2} \right) \! \left(\mathfrak{p}((-1)^{\varepsilon_{2}}a,
\rho) \chi_{1}(\vec{\widehat{g}}(\varepsilon_{1},\varepsilon_{2});\rho) \vert 
\tau \vert^{2 \rho} \! + \! \mathfrak{p}((-1)^{\varepsilon_{2}}a,-\rho) 
\right. \nonumber \\
\times&\left. \, \chi_{1}(\vec{\widehat{g}}(\varepsilon_{1},\varepsilon_{2});
-\rho) \vert \tau \vert^{-2 \rho} \right) \! \left(\mathfrak{p}((-1)^{
\varepsilon_{2}+1}a,\rho) \me^{-\mi \pi \rho} \chi_{2}(\vec{\widehat{g}}
(\varepsilon_{1},\varepsilon_{2});\rho) \vert \tau \vert^{2 \rho} \right. 
\nonumber \\
+&\left. \, \mathfrak{p}((-1)^{\varepsilon_{2}+1}a,-\rho) \me^{\mi \pi \rho} 
\chi_{2}(\vec{\widehat{g}}(\varepsilon_{1},\varepsilon_{2});-\rho) \vert \tau 
\vert^{-2 \rho} \right) \! \left(1 \! + \! \mathcal{O} \! \left(\tau^{\delta} 
\right) \right),
\end{align}
where
\begin{equation*}
\mathfrak{p}(z_{1},z_{2}) \! := \! \left(\dfrac{\vert \varepsilon b \vert}{32} 
\me^{\frac{\mi \pi}{2}} \right)^{z_{2}} \! \left(\dfrac{\Gamma (\frac{1}{2} \! 
- \! z_{2})}{\Gamma (1 \! + \! z_{2})} \right)^{2} \dfrac{\Gamma (1 \! + \! 
z_{2} \! + \! \frac{\mi z_{1}}{2})}{\tan (\pi z_{2})},
\end{equation*}
\begin{equation}
\begin{gathered}
\chi_{1}(\vec{\widehat{g}}(\varepsilon_{1},\varepsilon_{2});z_{3}) \! := \! 
\widehat{g}_{11}(\varepsilon_{1},\varepsilon_{2}) \me^{\mi \pi z_{3}} \me^{
\frac{\mi \pi}{4}} \! + \! \widehat{g}_{21}(\varepsilon_{1},\varepsilon_{2}) 
\me^{-\mi \pi z_{3}} \me^{-\frac{\mi\pi}{4}}, \\
\chi_{2}(\vec{\widehat{g}}(\varepsilon_{1},\varepsilon_{2});z_{4}) \! := \! 
\widehat{g}_{12}(\varepsilon_{1},\varepsilon_{2}) \me^{\mi \pi z_{4}} \me^{
\frac{\mi \pi}{4}} \! + \! \widehat{g}_{22}(\varepsilon_{1},\varepsilon_{2}) 
\me^{-\mi \pi z_{4}} \me^{-\frac{\mi \pi}{4}}.
\end{gathered}
\end{equation}

Let $\mathcal{H}(\tau)$ be the Hamiltonian function defined in 
Equation~{\rm (35)} corresponding to the function $u(\tau)$ given above. Then
\begin{align}
\mathcal{H}(\tau) \underset{\tau \to 0 \me^{\frac{\mi \pi \varepsilon_{1}}{2}
}}{=}& \, \dfrac{2 \rho}{\tau} \dfrac{\left(\mathfrak{p}((-1)^{\varepsilon_{2}}
a,\rho) \chi_{1}(\vec{\widehat{g}}(\varepsilon_{1},\varepsilon_{2});\rho) 
\vert \tau \vert^{2 \rho} \! - \! \mathfrak{p}((-1)^{\varepsilon_{2}}a,-\rho) 
\chi_{1}(\vec{\widehat{g}}(\varepsilon_{1},\varepsilon_{2});-\rho) \vert \tau 
\vert^{-2 \rho} \right)}{\left(\mathfrak{p}((-1)^{\varepsilon_{2}},\rho) \chi_{
1}(\vec{\widehat{g}}(\varepsilon_{1},\varepsilon_{2});\rho) \vert \tau \vert^{
2 \rho} \! + \! \mathfrak{p}((-1)^{\varepsilon_{2}}a,-\rho) \chi_{1}(\vec{
\widehat{g}}(\varepsilon_{1},\varepsilon_{2});-\rho) \vert \tau \vert^{-2 
\rho} \right)} \nonumber \\
+& \, \dfrac{1}{2 \tau} \! \left(a \! \left(a \! - \! (-1)^{\varepsilon_{2}} 
\mi \right) \! + \! \dfrac{1}{4} \! + \! 8 \rho^{2} \right) \! + \! o \! 
\left(\dfrac{1}{\tau} 
\right).
\end{align}
\end{ddddd}
\begin{ddddd}
Let $\varepsilon_{1} \! = \! \pm 1$, $\varepsilon_{2} \! = \! 0,\pm 1$, 
$\varepsilon b \! = \! \vert \varepsilon b \vert \me^{\mi \pi \varepsilon_{
2}}$, and $u(\tau)$ be a solution of Equation~{\rm (\ref{eq:dp3})} 
corresponding to the monodromy data 
$(a,s^{0}_{0},s^{\infty}_{0},s^{\infty}_{1},g_{11},g_{12},g_{21},g_{22})$. 
Suppose that
\begin{equation*}
\vert \Im (a) \vert \! < \! 1, \qquad \quad \widehat{g}_{11}(\varepsilon_{1},
\varepsilon_{2}) \widehat{g}_{22}(\varepsilon_{1},\varepsilon_{2}) \! \not= \! 
0, \qquad \quad s^{0}_{0} \! = \! 2 \mi.
\end{equation*}
Then $\exists \, \, \delta \! > \! 0$ such that $u(\tau)$ has the asymptotic 
expansion
\begin{align}
u(\tau) \underset{\tau \to 0 \me^{\frac{\mi \pi \varepsilon_{1}}{2}}}{=}& \, 
\dfrac{(-1)^{\varepsilon_{2}} \tau b \exp ((-1)^{\varepsilon_{2}} \frac{\pi 
a}{2})}{2a \sinh (\frac{\pi a}{2})} \! \left(\chi_{1}(\vec{\widehat{g}}
(\varepsilon_{1},\varepsilon_{2});0) \! \left(1 \! - \! \tfrac{(-1)^{
\varepsilon_{2}} \mi a}{2}Q((-1)^{\varepsilon_{2}}a) \right) \! + \! \dfrac{
(-1)^{\varepsilon_{2}} \pi a}{4} \right. \nonumber \\
\times&\left. \, (\widehat{g}_{21}(\varepsilon_{1},\varepsilon_{2}) \me^{-
\frac{\mi \pi}{4}} \! - \! 3 \widehat{g}_{11}(\varepsilon_{1},\varepsilon_{2}) 
\me^{\frac{\mi \pi}{4}}) \! + \! (-1)^{\varepsilon_{2}} \mi a \chi_{1}(\vec{
\widehat{g}}(\varepsilon_{1},\varepsilon_{2});0) \ln \vert \tau \vert \right) 
\nonumber \\
\times&\, \left(\chi_{2}(\vec{\widehat{g}}(\varepsilon_{1},\varepsilon_{2});0) 
\! \left(1 \! + \! \tfrac{(-1)^{\varepsilon_{2}} \mi a}{2}Q((-1)^{\varepsilon_{
2}+1}a) \right) \! + \! \dfrac{(-1)^{\varepsilon_{2}} \pi a}{4}(\widehat{g}_{1
2}(\varepsilon_{1},\varepsilon_{2}) \me^{\frac{\mi \pi}{4}} \right. \nonumber 
\\
-&\left. \, 3 \widehat{g}_{22}(\varepsilon_{1},\varepsilon_{2}) \me^{-\frac{
\mi \pi}{4}}) \! - \! (-1)^{\varepsilon_{2}} \mi a \chi_{2}(\vec{\widehat{g}}
(\varepsilon_{1},\varepsilon_{2});0) \ln \vert \tau \vert \right) \! \left(1 
\! + \! \mathcal{O} \! \left(\tau^{\delta} \right) \right),
\end{align}
where $\chi_{j}(\vec{\widehat{g}}(\varepsilon_{1},\varepsilon_{2});
\boldsymbol{\cdot})$, $j \! = \! 1,2$, are defined in Theorem~{\rm A.4}, 
Equations~{\rm (A.6)},
\begin{equation*}
Q(z) \! := \! 4 \psi (1) \! - \! \psi (\mi z/2) \! + \! \ln 2 \! - \! \ln 
(\vert \varepsilon b \vert),
\end{equation*}
$\psi (x) \! := \! \tfrac{\md}{\md x} \ln \Gamma (x)$ is the psi function, and 
$\psi (1) \! = \! -0.57721566490 \dotsc$.

Let $\mathcal{H}(\tau)$ be the Hamiltonian function defined in 
Equation~{\rm (35)} corresponding to the function $u(\tau)$ given above. 
Then
\begin{equation}
\mathcal{H}(\tau) \! \underset{\tau \to 0 \me^{\frac{\mi \pi \varepsilon_{1}}{
2}}}{=} \! \dfrac{1}{2 \tau} \! \left(a \! \left(a \! - \! (-1)^{\varepsilon_{
2}} \mi \right) \! + \! \dfrac{1}{4} \right) \! + \! \dfrac{\widehat{b}_{2}
(\varepsilon_{1},\varepsilon_{2})}{\tau (\widehat{a}_{2}(\varepsilon_{1},
\varepsilon_{2}) \! + \! \widehat{b}_{2}(\varepsilon_{1},\varepsilon_{2}) \ln 
\vert \tau \vert)} \! + \! o \! \left(\dfrac{1}{\tau} \right),
\end{equation}
where
\begin{align*}
\widehat{a}_{2}(\varepsilon_{1},\varepsilon_{2}) :=& \, \chi_{1}(\vec{\widehat{
g}}(\varepsilon_{1},\varepsilon_{2});0) \! \left(1 \! - \! \tfrac{(-1)^{
\varepsilon_{2}} \mi a}{2}Q((-1)^{\varepsilon_{2}}a) \right) \! + \! \dfrac{
(-1)^{\varepsilon_{2}} \pi a}{4}(\widehat{g}_{21}(\varepsilon_{1},\varepsilon_{
2}) \me^{-\frac{\mi \pi}{4}} \nonumber \\
-& \, 3 \widehat{g}_{11}(\varepsilon_{1},\varepsilon_{2}) \me^{\frac{\mi \pi}{
4}}),
\end{align*}
\begin{equation*}
\widehat{b}_{2}(\varepsilon_{1},\varepsilon_{2}) \! := \! (-1)^{\varepsilon_{
2}} \mi a \chi_{1}(\vec{\widehat{g}}(\varepsilon_{1},\varepsilon_{2});0).
\end{equation*}
\end{ddddd}
\begin{hhhhh}
For the conditions stated in Theorem~{\rm A.4},
\begin{align}
\pmb{\pmb{\boldsymbol{\tau}}}(\tau) \underset{\tau \to 0 \me^{\frac{\mi \pi 
\varepsilon_{1}}{2}}}{=}& \, \mathrm{const.} \, \tau^{\frac{1}{2} \left(a(a-
(-1)^{\varepsilon_{2}} \mi)+\frac{1}{4}+8\rho^{2} \right)} \! \left(\mathfrak{
p}((-1)^{\varepsilon_{2}}a,\rho) \chi_{1}(\vec{\widehat{g}}(\varepsilon_{1},
\varepsilon_{2});\rho) \vert \tau \vert^{2 \rho} \right. \nonumber \\
+&\left. \, \mathfrak{p}((-1)^{\varepsilon_{2}}a,-\rho) \chi_{1}(\vec{\widehat{
g}}(\varepsilon_{1},\varepsilon_{2});-\rho) \vert \tau \vert^{-2 \rho} \right) 
\! \left(1 \! + \! o \! \left(\tau^{\delta} \right) \right).
\end{align}
\end{hhhhh}
\begin{hhhhh}
For the conditions stated in Theorem~{\rm A.5},
\begin{equation}
\pmb{\pmb{\boldsymbol{\tau}}}(\tau) \! \underset{\tau \to 0 \me^{\frac{\mi \pi 
\varepsilon_{1}}{2}}}{=} \! \mathrm{const.} \, \tau^{\frac{1}{2} \left(a(a-
(-1)^{\varepsilon_{2}} \mi)+\frac{1}{4} \right)} \! \left(\widehat{a}_{2}
(\varepsilon_{1},\varepsilon_{2}) \! + \! \widehat{b}_{2}(\varepsilon_{1},
\varepsilon_{2}) \ln \vert \tau \vert \right) \! \left(1 \! + \! o \! \left(
\tau^{\delta} \right) \right).
\end{equation}
\end{hhhhh}

\end{document}